
\input amstex
\font\smallbf=ptmbo at 10pt

\define\bbR{{\bold R}}

\define\bbC{{\bold C}}

\define\bbCP{{\bold C}\text{\rm P}}

\define\rtr{\text{\bf R}\hskip-.7pt^3}

\define\spanr{\text{\rm Span}_{\hskip.4pt\bbR\hskip-1.2pt}}
\define\dimr{\dim_{\hskip.1pt\bbR\hskip-1.7pt}}
\define\dimc{\dim_{\hskip.4pt\bbC\hskip-1.2pt}}
\define\trc{\text{\rm Trace}\hs_\bbC}
\define\Lie{\text{\smallbf L}\hskip.5pt}
\define\sgn{\text{\rm sgn}\hs}

\define\y{c}
\define\bz{b\hs}

\define\fy{\sigma}
\define\vta{\varTheta}
\define\ax{A}
\define\bx{B}
\define\cx{C}
\define\fe{f}
\define\iy{\Cal I}
\define\iyp{\Cal I\hs'}
\define\jy{\Cal J}
\define\ly{\varPsi}
\define\dfe{d\hskip-.8ptf}

\define\e{\text{\rm exp}}

\define\hs{\hskip.7pt}
\define\nh{\hskip-1pt}
\define\ptmi{\phantom{i}}
\define\ptmii{\phantom{ii}}
\define\ns{\hskip-1.2pt}
\define\qd{\hskip7pt$\blacksquare$}
\define\dsq{d^{\hskip.6pt2}\hskip-1pt}
\define\dcu{d^{\hskip1.2pt3}\hskip-1pt}
\define\hrz{^{\hskip.5pt\text{\rm hrz}}}
\define\vrt{^{\hskip.2pt\text{\rm vrt}}}

\define\a{}
\define\f{\thetag}
\define\ff{\tag}

\define\vx{\varXi}

\define\vp{{\tau\hskip-4.55pt\iota\hskip.6pt}}
\define\vpsq{{\tau\hskip-4.55pt\iota\hskip1.3pt^2}}
\define\mgmt{M,g,m\hs,\hskip.2pt\vp}

\define\si{\phi}
\define\ta{\psi}
\define\la{\lambda}
\define\my{\mu}
\define\nk{{K}}

\define\sa{s}

\define\ps{t}

\define\zx{z}
\define\srf{S}
\define\gx{\gamma}

\define\ah{\alpha}
\define\ba{\beta}

\define\dv{\delta}
\define\diml{-di\-men\-sion\-al}
\define\skc{skew-sym\-met\-ric}
\define\sky{skew-sym\-me\-try}
\define\Sky{Skew-sym\-me\-try}
\define\kip{Killing potential}
\define\om{\omega}
\define\omh{\omega^{(h)}}
\define\ve{\varepsilon}
\define\kx{\kappa}

\define\ri{\text{\rm r}}
\define\rih{\ri^{(h)}}
\define\roh{\rho^{(h)}}
\define\sc{\text{\rm s}}
\define\navp{\nabla\hskip-.2pt\vp}
\define\ts{\eta}
\define\proj{\pi}

\define\krp{K\"ahler-Ricci potential}
\define\id{0}
\define\fn{1}
\define\pr{2}
\define\cc{3}
\define\kr{4}
\define\kp{5}
\define\ck{6}
\define\sr{7}
\define\xm{8}
\define\mc{9}
\define\de{10}
\define\qs{11}
\define\cm{12}
\define\rr{13}
\define\sh{14}
\define\hr{15}
\define\dx{16}
\define\lc{17}
\define\ct{18}
\define\sz{19}
\define\as{20}
\define\su{21}
\define\ty{22}
\define\ec{23}
\define\ls{24}
\define\tg{25}
\define\bb{26}
\define\is{27}

\define\acg{1}
\define\ber{2}
\define\bes{3}
\define\bon{4}
\define\bch{5}
\define\cal{6}
\define\ekm{7}
\define\cao{8}

\define\cvl{9}

\define\sdk{10}
\define\hem{11}
\define\dmr{12}
\define\dmg{13}
\define\grh{14}
\define\hws{15}
\define\hsi{16}
\define\kno{17}
\define\khn{18}
\define\leb{19}
\define\pag{20}
\define\ptv{21}
\define\pet{22}
\define\sak{23}
\define\snc{24}
\define\tiz{25}
\define\tfr{26}
\define\wan{27}
\documentstyle{amsppt}
\magnification=1200
\NoBlackBoxes
\topmatter
\title Local classification of conformally-Einstein\\
K\"ahler metrics in higher dimensions
\endtitle
\rightheadtext{Conformally-Einstein K\"ahler metrics}
\author A. Derdzinski and G. Maschler\endauthor
\address Dept. of Mathematics, The Ohio State University,
Columbus, OH 43210, USA\endaddress
\email andrzej\@math.ohio-state.edu\endemail
\address Department of Mathematics, University of Toronto, Canada M5S 3G3
\endaddress
\email maschler\@math.toronto.edu\endemail
\thanks 2000 {\it Mathematics Subject Classification.} Primary 53B35, 53B20;
Secondary 53C25, 53C55
\endthanks
\abstract{The requirement that a (non-Einstein) K\"ahler metric in any given
complex dimension $\,m>2\,$ be almost-everywhere conformally Einstein turns
out to be much more restrictive, even locally, than in the case of complex
surfaces. The local biholomorphic-isometry types of such metrics depend, for
each $\,m>2$, on three real parameters along with an arbitrary
K\"ahler-Einstein metric $\,h\,$ in complex dimension $\,m-1$. We provide an
explicit description of all these local-isometry types, for any
given $\,h$. That result is derived from a more general local classification
theorem for metrics admitting functions we call {\it special K\"ahler-Ricci
potentials}.}
\endabstract
\endtopmatter
\voffset=-35pt
\document
\head\S\id. Introduction\endhead
This is the first of three papers dealing with local and global properties of
conformally-Einstein K\"ahler metrics in higher dimensions. Specifically, we
study quadruples $\,\mgmt\,$ in which
\vskip4pt
\settabs\+\noindent&\f{\id.1}\hskip24pt&\cr
\+&&$(M,g)\hs$ is a K\"ahler manifold of complex dimension $\hs m\hs$
and $\hs\vp\hs$ is\cr
\+&\f{\id.1}&a nonconstant $\hs\,C^\infty\hs$ function on $\hs\,M\,\hs$ such
that the conformally\cr
\+&&related metric $\,\,\tilde g=g/\vpsq$,\hskip4ptdefined wherever
$\,\,\vp\ne0$,\hskip4ptis Einstein.\cr
\vskip4pt
\noindent As shown in Proposition \a\ck.4 below, if \f{\id.1} holds with
$\,m\ge3$, then also
$$\text{\rm$\mgmt\,$ satisfy \f{\id.1} and $\,d\vp\wedge\hs d\Delta\vp=0\,$
everywhere in $\,M$,}\ff\id.2$$
so that locally, at points with $\,d\vp\ne0$, the Laplacian of $\,\vp\,$ is a
function of $\,\vp$. As a result, \f{\id.2} is of independent interest only
for K\"ahler {\it surfaces\/} ($m=2$).

We obtain a complete local classification, at points in general position, of
all $\,\mgmt\,$ with \f{\id.1} and $\,m\ge3$, or \f{\id.2} and $\,m=2\,$ (see
Theorem \a\ls.1), and observe (in Remark \a\ls.2) that, for each fixed $\,m$,
their local biholomorphic-isometry types depend on three real constants plus
an additional ``parameter'' in the form of an arbitrary local
biholomorphic-isometry type of a K\"ahler-Einstein metric in complex dimension
$\,m-1$.

Our Theorem \a\ls.1 follows from a more general local classification (Theorem
\a\ct.1) of K\"ahler manifolds with {\it special\hskip5.5pt\krp s}, that
is, functions satisfying the second-order condition \f{\sr.1} in \S\sr, which
also involves the Ricci tensor.

A {\it global\/} classification of all $\,\mgmt\,$ with \f{\id.1} and
$\,m\ge3$, or \f{\id.2} and $\,m=2$, for which $\,M\,$ is compact, can
similarly be derived from a global classification of compact K\"ahler
manifolds with special \krp s. These classification theorems both require
extensive additional arguments, based on two different methods, and will
therefore appear in separate papers \cite{\dmr}, \cite{\dmg}.

The simplest examples of quadruples $\,\mgmt\,$ with \f{\id.1} for which
$\,M\,$ is compact, $\,m\ge2$, and $\,\vp=0\,$ somewhere in $\,M$, are
provided by some Riemannian products that have apparently been known for
decades; see \S\tg\ below. Another family of compact examples with \f{\id.1},
representing all dimensions $\,m>2$, and this time having $\,\vp\ne0\,$
everywhere in $\,M\,$ (so that $\,g\,$ is globally conformally Einstein), was
constructed by Lionel B\'erard Bergery \cite{\ber}; cf. also \S\bb\ below.
More recent extensions of the results of \cite{\ber} can be found in
\cite{\wan}.

Compact K\"ahler {\it surfaces\/} $\,(M,g)\,$ with \f{\id.1} have been studied
by many authors. For instance, Page \cite{\pag} found the only known example
of an Einstein metric on a compact complex surface (namely, on
$\,M=\bbCP^2\,\#\,\,\overline{\bbCP^2}$) which is globally
conformally-K\"ahler, but not K\"ahler; its conformally-related K\"ahler
metric was independently discovered by Calabi (\cite{\cal}, \cite{\cvl}).
Page's manifold thus satisfies \f{\id.1} (and, in fact, \f{\id.2}) with
$\,m=2$, and has $\,\vp\ne0\,$ everywhere in $\,M$. Examples of \f{\id.2} with
$\,m=2\,$ such that $\,M\,$ is compact and $\,\vp\,$ vanishes somewhere are
constructed, for minimal ruled surfaces $\,M\nh$, in \cite{\hws} and
\cite{\tfr}. Other results concerning \f{\id.1} for compact K\"ahler surfaces
$\,(M,g)\,$ include LeBrun's structure theorem \cite{\leb} for Hermitian
Einstein metrics on compact complex surfaces and the variational
characterization of such metrics in \cite{\snc}. For noncompact examples, see
\cite{\acg}, \cite{\hem}.

Condition \f{\id.1} with $\,m=2\,$ is much less restrictive than for
$\,m\ge3$, as it implies \f{\id.2} in the latter case, but not in the former.
This reflects the fact that, at points where the scalar curvature
$\,\,\sc\,\,$ of a K\"ahler-surface metric $\,g\,$ is nonzero, the metric
$\,\tilde g=g/\sc^2$ already satisfies a consequence of the Einstein condition
(namely, vanishing of the divergence of the self-dual Weyl tensor), and is, up
to a constant factor, the only metric conformal to $\,g\,$ with this property;
see \cite{\sdk, top of p. 417}.

\head\S\fn. A summary of contents\endhead
This paper is organized as follows. Sections \pr\ -- \kp\ contain preliminary
material. Basic facts on {\it conformally Einstein metrics\/} and the
definition of {\it special\hskip4pt\krp s}, along with their relation to each
other and examples of the latter, are presented in sections \ck\ -- \mc.
Sections \de\ and \qs\ deal with some ordinary differential equations
associated with the two italicized types of objects. In \S\cm\ we establish
``the $\,\hs\si\hs\,$ alternative'', stating that a certain function $\,\si\,$
naturally defined by a special \krp\ must be identically zero, or nonzero
everywhere. Next, in \S\rr\ we show that, essentially, $\,\si=0\,$ if and only
if the underlying K\"ahler metric is locally reducible. After some preparation
in \S\sh\ we introduce in \S\hr\ a simple generalization of Riemannian
submersions, needed to verify, in \S\dx, the claims we make about the examples
constructed in \S\xm, and then, in \S\ct, to prove our first classification
result, Theorem \a\ct.1, which also uses some lemmas from \S\lc. Theorem
\a\ct.1 asserts that locally, at points in general position, any K\"ahler
metric with a special \krp\ looks like one of the examples we describe in
\S\xm.

The remaining sections of the paper deal exclusively with conformally-Einstein
K\"ahler metrics or, more precisely, with quadruples $\,\mgmt\,$ satisfying
\f{\id.1} for $\,m\ge3$, or \f{\id.2} for $\,m=2$. First, in sections \sz\ --
\su, we explicitly solve the corresponding ordinary differential equations of
\S\de. We use the resulting three distinct families of solutions to introduce,
in \S\ty, three separate types of such quadruples. In \S\ec\ we use them again
to describe three families of examples, one for each type, obtained via a
special case of the construction in \S\xm. Then, in \S\ls\ we prove our main
classification result, stating that (locally, at points in general position)
any quadruple $\,\mgmt\,$ with \f{\id.1} for $\,m\ge3\,$ or \f{\id.2} for
$\,m=2\,$ belongs to one of those three families. Finally, sections \tg\ --
\is\ are appendices, in which we describe the well-known examples of compact
product manifolds with \f{\id.1}, explicitly realize B\'erard Bergery's
metrics \cite{\ber} as a special case of one of our three families, and
discuss some geometrically relevant integrals of the differential equations in
\S\de.

\head\S\pr. Preliminaries\endhead
For a $\,C^1$ vector field $\,v\,$ on a Riemannian manifold $\,(M,g)\,$ we
will write
$$\nabla v:TM\to TM\qquad\text{\rm with}\quad(\nabla v)w\,=\,\nabla_{\!w}v\,,
\ff\pr.1$$
i.e., treat the covariant derivative $\,\nabla v\,$ as a vector-bundle
morphism sending each $\,w\in T_xM$, $\,x\in M$, to $\,\nabla_{\!w}v\in T_xM$.
Except for \S\cc, the symbol $\,\nabla\,$ will denote both the Levi-Civita
connection of a given Riemannian metric $\,g\,$ on a manifold $\,M$, and the
$\,g$-gradient operator acting on $\,C^1$ functions $\,M\to\bbR\,$ (see
\f{\pr.4}). Then
$$\alignedat2
&\text{\rm\ptmi i)}\quad&&
2g(\nabla_{\!w}v,u)\,=\,d_w[g(v,u)]\,+\,d_v[g(w,u)]\,-\,d_u[g(w,v)]\hskip55pt\\
&&&\hskip49pt+\,g(v,[u,w])\,+\,g(u,[w,v])\,-\,g(w,[v,u])\,,\\
&\text{\rm ii)}\quad&&
\nabla_{\!v}w\,-\,\nabla_{\!w}v\,=\,[v,w]\,,\qquad\text{\rm i.e.,}\hskip5pt
\nabla\hskip5pt\text{\rm is\ torsion-free,}\endalignedat\ff\pr.2$$
for $\,C^2$ vector fields $\,v,w,u$. In fact, \f{\pr.2.i} follows from
\f{\pr.2.ii} and the Leibniz rule for the directional derivatives
$\,d_w,\hs d_v,\hs d_u$ (cf. \cite{\kno, Proposition 2.3 on p. 160}). Next,
for a mul\-ti\-ply-co\-var\-i\-ant tensor $\,\bz\,$ and a vector (field)
$\,v$, we set
$$\imath_v\bz\,=\,\bz(v,\,\cdot\,,\dots,\,\cdot\,)\hs.\ff\pr.3$$
We will apply $\,\imath_v$ to $\,1$-forms, twice-co\-var\-i\-ant tensors (such
as the metric $\,g$, the Ricci tensor $\,\,\ri$, or 
derivative $\,\nabla d\vp$, as in 
\f{\pr.5}), and differential $\,2$-forms. For instance,
$$\imath_vg\,=\,g(v,\,\cdot\,)\,=\,d\vp\,,\qquad\text{\rm where}\quad
v=\navp\quad\text{\rm is\ the\ gradient\ of}\hskip9pt\vp\,.\ff\pr.4$$
More generally, let $\,\xi\,$ be a $\,1$-form. Then $\,\imath_vg=\xi\,$
whenever $\,\xi=g(v,\,\cdot\,)\,$ is the $\,1$-form corresponding to $\,v\,$
via $\,g$, i.e., obtained from $\,v\,$ by ``lowering the index''.

The second covariant derivative $\,\nabla d\vp\,$ of a $\,C^2$ function
$\,\vp\,$ on a Riemannian manifold clearly has the property that, for any
vector fields $\,u,w$,
$$(\nabla d\vp)(u,w)\,=\,g(u,\nabla_{\!w}v)\,=\,g(\nabla_{\!u}v,w)\,,\qquad
\text{\rm where}\quad v\,=\,\navp\,.\ff\pr.5$$
\remark{Remark \a\pr.1}By \f{\pr.5}, for any $\,C^2$ function $\,\vp\,$ on a
Riemannian manifold, the eigenvalues and eigenvectors of the symmetric
$\,2$-tensor $\,\nabla d\vp\,$ at any point are the same as those of
$\,\nabla v\,$ with \f{\pr.1}, where $\,v=\navp$.
\endremark
\medskip Any $\,C^2$ function $\,\vp\,$ in a Riemannian manifold $\,(M,g)\,$
satisfies the relations
$$\alignedat2
&\text{\rm\ptmii i)}\quad&&
\imath_vd\vp\,=\,d_v\vp\,=\,Q\,,\quad\text{\rm with}\quad v\,=\,\navp\quad
\text{\rm and}\quad Q\,=\,g(\navp,\navp)\,,\\
&\text{\rm\ptmi ii)}\quad&&
2\hskip1pt\imath_v\bz\,=\,dQ\,,\quad\text{\rm where}\quad v
=\navp,\hskip6pt\bz
=\nabla d\vp\hskip6pt\text{\rm and}\hskip6ptQ=g(\navp,\navp)\,,\hskip5pt\\
&\text{\rm iii)}\quad&&
\imath_w\bz\,=\,\nabla_{\!w}d\vp\qquad\text{\rm for}\quad \bz=\nabla d\vp
\quad\text{\rm and\ any\ vector\ field}\quad w\,.\endalignedat\ff\pr.6$$
In fact, \f{\pr.6.i} is clear as $\,d_v\vp=g(v,\navp)$, \f{\pr.6.ii} is
immediate from the local-coordinate expression
$\,2\vp^{,k}\vp_{,jk}=[\vp^{,k}\vp_{,k}]_{,j}$, and \f{\pr.6.iii} is obvious
from \f{\pr.3}, \f{\pr.5}.

We denote $\,\langle\hskip1pt\bz,\bz'\rangle\,$ the natural inner product of
twice-co\-var\-i\-ant tensors $\,\bz,\bz'$ at any point $\,x\,$ of a
Riemannian manifold $\,(M,g)$, with
$$\langle\hskip1pt\bz,\bz'\rangle\,=\,\bz^{jk}\bz'_{jk}\,
=\,\,\text{\rm Trace}\hskip1.8ptB^*\hskip-1.1ptB'\,,\ff\pr.7$$
where the components $\,\bz^{jk}=g^{jp}g^{kq}\bz_{pq}$ and $\,\bz'_{jk}$ refer
to any local coordinate system, while $\,B:T_xM\to T_xM\,$ is related to
$\,\bz\,$ via $\,\bz(u,w)=g(Bu,w)\,$ for all $\,u,w\in T_xM\,$ (and similarly
for $\,B'\nh,\hs\bz'$). In particular, we have the {\it$\,g$-trace\/} of a
twice-co\-var\-i\-ant tensor $\,\bz$, given by
$\,\,\text{\rm Trace}_g\hs\bz=g^{jk}\bz_{jk}$, and so, for $\,B\,$ related to
$\,\bz\,$ as above,
$$\text{\rm Trace}_g\hs b\,=\,\langle\hskip1ptb,g\rangle\,=\,\,\text{\rm Trace}\,B\hs.
\ff\pr.8$$
For the tensor and exterior products and the exterior derivative of
$\,1$-forms $\,\xi$, $\,\xi'$, and $\,C^1$ tangent vector fields $\,u,v\,$ on
any manifold, we have
$$\alignedat2
&\text{\rm\ptmi i)}\quad&&
(\xi\hskip-.7pt\otimes\xi')(u,v)\,=\,\xi(u)\hs\xi'(v)\,,\qquad\quad
\xi\wedge\xi'\,=\,\xi\hskip-.7pt\otimes\xi'\,-\,\xi'\hskip-1.7pt\otimes\xi\,,\\
&\text{\rm ii)}\quad&&
(d\hs\xi)(u,v)\,=\,d_u[\xi(v)]\,-\,d_v[\xi(u)]\,-\,\xi([u,v])\,.
\hskip100pt\endalignedat\ff\pr.9$$
In any Riemannian manifold, the divergence operator $\,\dv\,$ acts on
mul\-ti\-ply-co\-var\-i\-ant $\,C^1$ tensor fields $\,\bz\,$ and $\,C^1$
vector fields $\,v\,$ by $\,\dv\bz=\sum_u\imath_u(\nabla_{\!u}\bz)\,$ and
$\,\dv v=\,\text{\rm Trace}\hskip1.8pt\nabla v$, cf. \f{\pr.3}, \f{\pr.1},
with summation over the vectors $\,u\,$ of any orthonormal basis of the
tangent space at any given point; the two formulae are consistent, as
$\,\dv v=\dv\xi\,$ whenever $\,\xi=g(v,\,\cdot\,)=\imath_vg$. Thus, for
$\,C^2$ functions $\,\vp$, $\,C^1$ vector fields $\,v,w$, $\,C^1$ $\,1$-forms
$\,\xi\,$ and twice-co\-var\-i\-ant $\,C^1$ tensor fields $\,\bz$,
$$\alignedat2
&\text{\rm\ptmii i)}\quad&&
\Delta\vp\,=\,\dv(d\vp)\,=\,\dv v\,,\hskip9pt\text{\rm where}\hskip6ptv\,
=\,\navp\hskip6pt\text{\rm and}\hskip6pt\Delta\hskip6pt\text{\rm is\ the\
Laplacian,}\\
&\text{\rm\ptmi ii)}\quad&&
\dv\hskip.6pt[\vp w]\,=\,d_w\vp\,+\,\vp\hskip1pt\dv w\,,\hskip9pt
\text{\rm and}\hskip9pt
\dv\hskip.6pt[\vp\bz]\,=\,\imath_v\bz\,+\,\vp\hskip1pt\dv\bz\hskip7pt
\text{\rm if}\hskip7ptv=\nabla\!\vp,\hskip0pt\\
&\text{\rm iii)}\quad&&
\imath_w[(\imath_vg)\otimes\xi]\,=\,g(v,w)\hs\xi\hskip7pt\text{\rm and}
\hskip7pt\dv\hs[(\imath_vg)\otimes\xi]\,=\,(\dv v)\hs\xi\,+\,\nabla_{\!v}\xi\,.
\endalignedat\ff\pr.10$$
In fact, \f{\pr.10.i}, \f{\pr.10.ii} are obvious, and the two equalities in
\f{\pr.10.iii} follow from \f{\pr.3}, \f{\pr.4}, \f{\pr.9.i} and,
respectively, the obvious local-coordinate identity
$\,(v^k\xi_j){}_{,k}=v^k{}_{,k}\hs\xi_j+v^k\xi_{j,k}$.

The symbols $\,R\hs,\,\ri\,\,$ and $\,\,\sc\,=\,\text{\rm Trace}_g\hs\ri\,\,$
(cf. \f{\pr.8}) usually stand for the curvature tensor, Ricci tensor and
scalar curvature of a given Riemannian metric $\,g$. Thus,
\vskip4pt
\settabs\+\noindent&\f{\pr.9}\hskip11pt&c)\hskip8pt&\cr
\+&&a)&$2\hs\dv\hs\ri\,\,=\,d\hs\sc\,$,\cr
\+&\f{\pr.11}&b)&$\imath_v\ri\,=\dv\bz-\hs d\hs\dv v\,$ for any $\,C^2$
vector field $\,v$, with $\,\bz(u,w)=g(u,\nabla_{\!w}v)$,\cr
\+&&c)&$\dv\bz\,=\,\hs\imath_v\ri\,\,+\,\hs d\hs\Delta\vp\,$ if $\,v=\navp\,$
and $\,\bz=\nabla d\vp\,$ for any $\,C^2$ function $\,\vp$,\cr
\vskip5pt
\noindent with $\,\dv\bz,\dv v\,$ as above (for $\,\bz=\,\ri\,\,$ in
\f{\pr.11.a}) and $\,\imath_v\,$ given by \f{\pr.3}. The first two equalities,
known as the {\it Bianchi identity for the Ricci tensor\/} and the {\it
Bochner formula\/} read, in local coordinates,
$\,2\hs\ri_{kj,}{}^k=\hskip1pt\sc\hskip1pt_{,j}$ and
$\,\,\ri_{jk}v^k=v^k{}_{,jk}-v^k{}_{,kj}$ (see, e.g., \cite{\bes, Proposition
1.94 on p. 43} and \cite{\bes, formula (2.51)}). The last equality in
\f{\pr.11.b}, with arbitrary tangent vectors $\,u,w$, {\it defines\/} the
twice-co\-var\-i\-ant tensor field $\,\bz\,$ appearing in the first equality.
Finally, \f{\pr.11.c} is a special case of \f{\pr.11.b} (cf. \f{\pr.10.i} and
\f{\pr.5}).

Finally, a $\,C^1$ vector field $\,v\,$ on a Riemannian manifold $\,(M,g)\,$
is locally a gradient if and only if $\,\nabla v\,$ with \f{\pr.1} is
self-adjoint at every point. This is clear from \f{\pr.4}, as \f{\pr.9.ii} for
$\,\xi=\imath_vg\,$ gives, by \f{\pr.2.ii},
$\,(d\hs\xi)(u,v)=g(\nabla_{\!u}v,w)-g(u,\nabla_{\!w}v)$.

\head\S\cc. Connections and curvature\endhead
In this section we depart from our usual conventions and let $\,\nabla\,$ and
$\,R\,$ stand for an arbitrary (linear) connection in any real/com\-plex
vector bundle $\,\Cal E\,$ over a manifold and, respectively, the curvature
tensor of $\,\nabla$, with the sign convention
$$R(u,v)w\,=\,\nabla_{\!v}\nabla_{\!u}w\,-\,
\nabla_{\!u}\nabla_{\!v}w\,+\,\nabla_{[u,v]}w\ff\cc.1$$
for $\,C^2$ vector fields $\,u,v\,$ tangent to the base and a $\,C^2$ section
$\,w\,$ of $\,\Cal E$. Denoting $\,m\,$ the real/com\-plex fibre dimension of
$\,\Cal E\,$ and choosing $\,C^\infty$ sections $\,w_a$, $\,a=1,\dots,m$,
which trivialize the bundle $\,\Cal E\,$ over some open subset $\,\,U\,$ of
the base, we have
$$\nabla_{\!v}w_a=\varGamma_{\!a}^{\hs b}(v)\hs w_b\qquad
(\text{\rm summed\ over}\quad b=1,\dots,m)\ff\cc.2$$
for any vector field $\,v\,$ on $\,\,U$, with some real/com\-plex valued
$\,1$-forms $\,\varGamma_{\!a}^{\hs b}$. From \f{\cc.1} and \f{\pr.9}, we
now obtain, with summation over $\,c=1,\dots,m$,
$$R(u,v)w_a\,=\,R_a^{\hs c}(u,v)\hs w_c\,,\quad\text{\rm where}\quad
R_a^{\hs b}\,=\,-\hs d\varGamma_{\!a}^{\hs b}\,
+\,\varGamma_{\!a}^{\hs c}\wedge\varGamma_{\!c}^{\hs b}\,.\ff\cc.3$$
\remark{Remark \a\cc.1}Given a connection $\,\nabla\,$ in a complex line
bundle over any manifold, let $\,\varOmega\,$ denote the {\it curvature
form\/} of $\,\nabla$, that is, the $\,2$-form characterized by
$\,R(u,v)w=i\hs\varOmega(u,v)w\,$ for $\,u,v,w,R\,$ as in \f{\cc.1}. Any local
$\,C^\infty$ section $\,w\,$ without zeros then gives rise to the {\it
connection\/} $\,1${\it-form\/} $\,\varGamma\,$ defined by
$\,\nabla_{\!v}w=\varGamma(v)\hs w$, i.e., \f{\cc.2} with $\,m=1$. Thus,
\f{\cc.3} with $\,m=1\,$ yields $\,\varOmega=i\hskip1ptd\varGamma$.
\endremark
\medskip
For any vector bundle $\,\Cal L\,$ over a manifold $\,N$, we will write
$$\Cal L\,=\,\{(y,\zx):y\in N,\hskip6pt\zx\in\Cal L_y\}\qquad\text{\rm and}
\qquad N\,\subset\,\Cal L\,,\ff\cc.4$$
using the same symbol $\,\Cal L\,$ for its total space, and identifying
$\,N\,$ with the zero section formed by all $\,(y,0)\,$ with $\,y\in N$. We
similarly treat each fibre $\,\Cal L_y$ as a subset of $\,\Cal L$, identifying
it with $\,\{y\}\times\Cal L_y$. Thus, $\,N\,$ and all $\,\Cal L_y$ are
submanifolds of $\,\Cal L\,$ with its obvious manifold structure. Being a
vector space, every fibre $\,\Cal L_y$ may also be identified with a subspace
of $\,T_{(y,\zx)}\Cal L$, for any $\,\zx\in\Cal L_y$ (namely, with the tangent
space of the submanifold $\,\{y\}\times\Cal L_y$ at $\,(y,\zx)$). The vertical
distribution $\,\Cal V\,$ on the total space $\,\Cal L\,$ is the subbundle of
$\,T\Cal L\,$ with the fibres
$\,\Cal V_{(y,\zx)}=\Cal L_y\subset T_{(y,\zx)}\Cal L$. In the case where
$\,\Cal L\,$ is a complex line bundle over $\,N$, using the notation of
\f{\cc.4} and a fixed real number $\,a\ne0$, we may define two vertical vector
fields (i.e., sections of the vertical distribution $\,\Cal V$) on the total
space $\,\Cal L\,$ by
$$v(y,\zx)\,=\,a\zx\,,\qquad u(y,\zx)\,=\,i\hs a\zx\,.\ff\cc.5$$
\remark{Remark \a\cc.2}Any Hermitian fibre metric $\,\langle\,,\rangle\,$ in
a complex line bundle $\,\Cal L\,$ over a manifold $\,N\,$ is determined by
its {\it norm function\/} $\,r:\Cal L\to[\hs0,\infty)$, assigning
$\,|\zx|=\langle\zx,\zx\rangle^{1/2}$ to each $\,(y,\zx)\in\Cal L$. As
$\,\Cal V_{(y,\zx)}=\Cal L_y$ (see above), $\,\langle\,,\rangle\,$ may also be
treated as a Hermitian fibre metric in the vertical subbundle $\,\Cal V\,$ of
$\,T\Cal L$, and then, for $\,v,u\,$ given by \f{\cc.5},
$\,\langle v,v\rangle=\langle u,u\rangle=a^2r^2$ and
$\,\,\text{\rm Re}\hskip1pt\langle v,u\rangle=0$, while, on
$\,\Cal L\smallsetminus N$,
$$d_vr=ar\hskip6pt\text{\rm and}\hskip6.5ptd_ur=\hs d_wr=0\hskip9pt
\text{\rm for\ all\ horizontal\ vectors}\hskip4.5ptw\hs,\ff\cc.6$$
``horizontality'' referring to any given connection in $\,\Cal L\,$ that makes
$\,\langle\,,\rangle\,$ parallel.
\endremark
\medskip
Let $\,\varOmega\,$ and $\,\Cal H\,$ be the curvature form (Remark \a\cc.1)
and the horizontal distribution of a given connection in a complex line bundle
$\,\Cal L\,$ over a manifold $\,N$. If we use an unnamed local
trivializing section of $\,\Cal L$, defined on an open set $\,N'\subset N$, to
introduce the $\,1$-form $\,\varGamma\,$ as in Remark \a\cc.1, and also to
identify the portion $\,\Cal L'$ of $\,\Cal L\,$ lying over $\,N'$ with
$\,N'\nh\times\bbC$, then, for any $\,(y,\zx)\in\Cal L'$ and
$\,(w,\zeta)\in T_{(y,\zx)}\Cal L'$,
$$(w,\zeta)\vrt\,=\,(0,\zeta+\varGamma(w)\zx)\,,\qquad
(w,\zeta)\hrz\,=\,(w,-\varGamma(w)\zx)\,,\ff\cc.7$$
with $\,...\vrt,\hs...\hrz$ standing for the $\,\Cal V\,$ and $\,\Cal H\,$
components relative to the decomposition $\,T\Cal L=\Cal H\oplus\Cal V$. (In
fact, the formula for $\,(w,\zeta)\hrz$ follows if one writes down the
parallel-transport equation using $\,\varGamma$.) Since
$\,\varOmega=i\hskip1ptd\varGamma\,$ (Remark \a\cc.1), we have
$$[\hs\tilde w\hs,\tilde w{}'\hs]\vrt\,=\,i\hs\varOmega(w,w')\zx\qquad
\text{\rm at\ any}\quad(y,\zx)\in\Cal L\,,\ff\cc.8$$
cf. \f{\cc.4}, $\,\tilde w,\tilde w{}'$ being the horizontal lifts to
$\,\Cal L\,$ of any $\,C^\infty$ vector fields $\,w,w'$ on $\,N$. (This
follows from \f{\pr.9.ii}, as \f{\cc.7} gives
$\,\tilde w=(w,-\varGamma(w)\zx)\,$ at $\,(y,\zx)$.) Since
$\,\Cal L_y\subset T_{(y,\zx)}\Cal L$, both sides of \f{\cc.8} may be treated
as vectors tangent to $\,\Cal L$.

\head\S\kr. K\"ahler metrics and the Ricci form\endhead
Except in a few cases (such as Lemma \a\kr.2), we deal with K\"ahler manifolds
in real terms, using real-val\-ued functions and differential forms, real
vector fields (rather than complexified ones), and components of tensors in
{\it real\/} coordinate systems.

For any complex manifold $\,M\,$ we denote $\,J\,$ its complex-structure
tensor, treated either as a real vector-bundle morphism $\,J:TM\to TM\,$ with
$\,J^2=-1$, or as the multiplication by $\,\,i\,\,$ in the complex vector
bundle $\,TM$.

As usual, a K\"ahler manifold $\,(M,g)\,$ is a complex manifold $\,M\,$ with a
Riemannian metric $\,g\,$ that makes the complex-structure tensor $\,J\,$
skew-adjoint and parallel. We denote $\,\hs\om\,$ and $\,\rho\,$ its
K\"ahler and Ricci forms, so that, if $\,\,\ri\,\,$ is the Ricci tensor,
$$\om(u,v)\,=\,g(Ju,v)\,,\hskip11pt\rho(u,v)\,=\,\ri\hs(Ju,v)\hskip11pt
\text{\rm for}\hskip6ptu,v\in T_xM,\hskip4ptx\in M.\ff\kr.1$$
\Sky\ of $\,\,\om\,$ and $\,\rho\,$
reflects the fact that $\,g\,$ and $\,\,\ri\,\,$ are Hermitian, in the sense
that both $\,\bz=g\,$ and $\,\bz=\,\ri\,\,$ are symmetric
twice-co\-var\-i\-ant tensors with
$$\bz(Jw,w')\,=\,-\hs\bz(w,Jw')\qquad\text{\rm for\ all}\quad w,w'\in T_xM\quad
\text{\rm and}\quad x\in M.\ff\kr.2$$
In any (real) local coordinates, $\,\rho(u,v)=\rho_{jk}u^jv^k$, where, by
\f{\kr.1} -- \f{\kr.2},
$$\rho_{jk}\,=-\,\,\ri_{jl}J_k^l\,.\ff\kr.3$$
The following well-known lemmas state that $\,\rho\,$ is the curvature form
(Remark \a\cc.1) of a natural connection in the highest complex exterior power
of $\,TM$.
\proclaim{Lemma \a\kr.1}Let\/ $\,u,v\,$ be\/ $\,C^2$ vector fields on a
K\"ahler manifold\/ $\,(M,g)$. Then\/
$\,\,\text{\rm Trace}_{\hskip.4pt\bbC}[R(u,v)]=i\hs\rho(u,v)$, where\/
$\,\rho\,$ is the Ricci form and\/ $\,R(u,v):TM\to TM\,$ denotes the
complex-linear bundle morphism\/ $\,w\mapsto R(u,v)w\,$ satisfying\/ \f{\cc.1}.
\endproclaim
\demo{Proof}Both $\,J\,$ and the unique Hermitian fibre metric
$\,\langle\,,\rangle\,$ in $\,TM\,$ whose real part is $\,g\,$ are
$\,\nabla$-parallel. Thus, $\,\nabla\,$ is a connection in the {\it complex\/}
bundle $\,TM$, making $\,\langle\,,\rangle\,$ parallel, and so $\,R(u,v)\,$ is
complex-linear, i.e., the morphisms $\,R(u,v)\,$ and $\,J\,$ commute. Since
they both are skew-adjoint relative to $\,\langle\,,\rangle$, their composite
$\,F$, with fixed vector fields $\,u,v$, is self-adjoint. Hence
$\,2i\,\text{\rm Trace}_{\hskip.4pt\bbC}[R(u,v)]
=\,2\,\text{\rm Trace}_{\hskip.4pt\bbC}F
=\,\text{\rm Trace}_{\hskip.4pt\bbR}F$, the last equality being clear if one
evaluates both traces using a complex basis $\,e_1,\dots,e_m$ of any
given tangent space which diagonalizes $\,F\,$ and, respectively, the real
basis $\,e_1,\,ie_1,\,\dots\,,e_m,\,ie_m$.

In terms of components relative to any real local-coordinate system,
self-adjoint\-ness of $\,F\,$ (or, skew-adjoint\-ness of $\,J$) means that
$\,R_{jpkq}J_l^q$ is symmetric in $\,k,l\,$ (or, respectively,
$\,\ah^{jk}=g^{jl}J_l^k$ is \skc\ in $\,j,k$). The Bianchi identity thus
gives $\,(R_{jpkq}-R_{jkpq}-R_{jqkp})\ah^{pq}=0$, i.e.,
$\,2\hs R_{jpkq}\ah^{pq}=R_{jkpq}\ah^{pq}=-\hs R_{jkq}{}^pJ_p^q$, while,
by \f{\kr.3}, $\,R_{jpkq}\ah^{pq}=g^{pl}R_{jpkq}J_l^q=g^{pl}R_{jplq}J_k^q
=-\hs\ri_{jq}J_k^q=\rho_{jk}$. Combining these equalities and summing against
$\,u^jv^k$, we get $\,2\hs\rho(u,v)=\,-\hs\text{\rm Trace}_{\hskip.4pt\bbR}F
=-\hs2i\,\text{\rm Trace}_{\hskip.4pt\bbC}[R(u,v)]$, which completes the
proof.{\hfill\qd}
\enddemo
\proclaim{Lemma \a\kr.2}Let\/ $\,\rho\,$ be the Ricci form of a K\"ahler
manifold\/ $\,(M,g)\,$ of complex dimension\/ $\,m\,$ with\/ $\,C^\infty$
vector fields\/ $\,w_a$, $\,a=1,\dots,m$, on an open set\/ $\,\,U\subset M\,$
which trivialize the tangent bundle\/ $\,TU\,$ treated as a complex vector
bundle, and let\/ $\,\varGamma_{\!a}^{\hs b}$ be the complex-val\-ued\/
$\,1$-forms with\/ \f{\cc.2} for the Levi-Civita connection\/ $\,\nabla\,$
of\/ $\,g$. Then\/ $\,\rho\,=\,i\hskip1ptd\varGamma_{\!a}^{\hs a}$ on\/
$\,\,U$, with summation over\/ $\,a=1,\dots,m$.
\endproclaim
In fact, the matrix $\,[R_a^{\hs b}(u,v)]\,$ (see \f{\cc.3}) represents the
operator $\,R(u,v)\,$ of Lemma \a\kr.1 in the complex basis $\,w_a(x)$ at any
point $\,x\in U$. Now Lemma \a\kr.1 and \f{\cc.3} give
$\,i\rho=R_a^{\hs a}=-\hs d\varGamma_{\!a}^{\hs a}$ (with
$\,\varGamma_{\!a}^{\hs c}\wedge\varGamma_{\!c}^{\hs a}=0\,$ due to obvious
cancellations).{\hfill\qd}

\head\S\kp. Killing potentials\endhead
A real-val\-ued $\,C^\infty$ function $\,\vp\,$ on a K\"ahler manifold
$\,(M,g)\,$ is called a {\it Killing potential\/} if $\,u=J(\navp)\,$ is a
Killing vector field on $\,(M,g)$.
\remark{Remark \a\kp.1}As usual, a {\it Killing field\/} on a Riemannian
manifold $\,(M,g)\,$ is any $\,C^\infty$ vector field $\,u\,$ such that
$\,\nabla u\,$ is skew-adjoint at every point (cf. \f{\pr.1}). Then $\,u\,$
is uniquely determined by $\,u(x)\,$ and $\,(\nabla u)(x)\,$ at any given
point $\,x\in M$. In fact, $\,u\,$ is a Jacobi field along any geodesic, since
its local flow, applied to a geodesic segment, generates a variation of
geodesics. This implies a {\it unique continuation property\/}: a Killing
field is uniquely determined by its restriction to any nonempty open set.
\endremark
\medskip
We call a (real) $\,C^\infty$ vector field $\,v\,$ on a complex manifold {\it
hol\-o\-mor\-phic\/} if $\,\Lie_vJ=0$, where $\,\Lie\,$ is the Lie derivative.
For a $\,C^\infty$ vector field $\,v\,$ on a {\it K\"ahler\/} manifold,
\vskip4pt
\settabs\+\noindent&\f{\kp.1}\hskip11pt&c)\hskip8pt&\cr
\+&&a)&$v\,$ is holomorphic if and only if $\,[J,\nabla v]=0$,\cr
\+&\f{\kp.1}&b)&$[v,u]=0\,$ and $\,u\,$ is holomorphic if $\,v\,$ is
holomorphic and $\,u=Jv$,\cr
\+&&c)&$\nabla u=\hs J\circ(\nabla v)$, with the convention \f{\pr.1}, if
$\,u=Jv$.\cr
\vskip5pt
\noindent In fact, \f{\kp.1.a} (or, $\,[v,u]=0\,$ in \f{\kp.1.b}) follows if
we apply the relation $\,(\Lie_vJ)w=\Lie_v(Jw)-J(\Lie_vw)$, valid for any
$\,C^2$ vector fields $\,v,w$, to any $\,w$, obtaining
$\,\Lie_vJ=[J,\nabla v]\,$ due to \f{\pr.2.ii} and \f{\pr.1} with
$\,\nabla J=0\,$ (or, to $\,w=v$, noting that $\,\Lie_vw=[v,w]$). Now
\f{\kp.1.c} is clear as $\,\nabla_{\!w}u=\hs J(\nabla_{\!w}v)\,$ for every
tangent vector $\,w\,$ (since $\,\nabla J=0$), while \f{\kp.1.a} and
\f{\kp.1.c} yield the remainder of \f{\kp.1.b}.
\proclaim{Lemma \a\kp.2}For a\/ $\,C^\infty$ function\/ $\,\vp\,$ on a
K\"ahler manifold\/ $\,(M,g)$, the following three conditions are
equivalent\/{\rm:}\hskip7pt{\rm(i)}\hskip5pt$\vp\,$ is a \kip\/{\rm;}
\hskip7pt{\rm(ii)}\hskip5ptThe gradient\/ $\,v=\navp\,$ is a holomorphic
vector field\/{\rm;}\hskip7pt{\rm(iii)}\hskip5pt$\bz=\nabla d\vp\,$ is
Hermitian, as in\/ \f{\kr.2}.
\endproclaim
In fact, let $\,v=\navp\,$ and $\,u=Jv$. By \f{\kp.1.c},
$\,\nabla u=J\circ(\nabla v)\,$ and $\,[\nabla v]^*=\nabla v$, that is,
$\,\nabla v\,$ in \f{\pr.1} is self-adjoint at every point. Thus,
$\,\nabla u+[\nabla u]^*=\,J\circ(\nabla v)+[\nabla v]^*\circ J^*
=[J,\,\nabla v]\,$ and, by \f{\kp.1.a}, (i) is equivalent to (ii). Also, by
\f{\kr.2} for $\,\bz=g\,$ and \f{\pr.5},
$\,(\nabla d\vp)(Jw',w)=g(Jw',\nabla_{\!w}v)=-\hs g(w',\nabla_{\!w}u)$, so
that Hermitian symmetry of $\,\nabla d\vp\,$ amounts to skew-adjoint\-ness of
$\,\nabla u$, i.e., to (i).{\hfill\qd}
\medskip
Next, there is a well-known local one-to-one correspondence between \kip s
defined up to an additive constant, and holomorphic Killing vector fields:
\proclaim{Lemma \a\kp.3}Let\/ $\,(M,g)\,$ be a K\"ahler manifold. For every
\kip\/ $\,\vp\,$ on\/ $\,(M,g)$, the Killing field\/ $\,J(\navp)\,$ is
hol\-o\-mor\-phic. Conversely, if\/ $\,H^1(M,\bbR)=\{0\}$, then every
holomorphic Killing vector field on\/ $\,(M,g)\,$ has the form\/
$\,J(\navp)\,$ for a \kip\/ $\,\vp$, which is unique up to an additive
constant.
\endproclaim
The first assertion is clear from Lemma \a\kp.2(ii) and \f{\kp.1.b}. Next,
setting $\,v=-\hs Ju\,$ for a holomorphic Killing field $\,u$, we get
$\,\nabla v=-\hs J\circ(\nabla u)\,$ (by \f{\kp.1.c}). As $\,J\,$ and
$\,\nabla u\,$ commute (cf. \f{\kp.1.a}) and are skew-adjoint, their composite
$\,\nabla v\,$ is self-adjoint, i.e., locally, $\,v=\navp\,$ for some function
$\,\vp\,$ (see end of \S\pr).{\hfill\qd}
\remark{Remark \a\kp.4}Let $\,\vp\,$ be a nonconstant \kip\ on a K\"ahler
manifold $\,(M,g)$. Then $\,\nabla d\vp\ne0\,$ wherever $\,d\vp=0\,$ (and
hence $\,\,d\vp\ne0\,$ on a dense open subset of $\,M$, which also follows
from Lemma \a\kp.2(ii)). In fact, if $\,\nabla d\vp\,$ and $\,d\vp\,$ both
vanished at some point, so would $\,v=\navp\,$ and $\,\nabla v\,$ (by
\f{\pr.5}), as well as $\,u=Jv\,$ and $\,\nabla u\,$ (since
$\,\nabla u=J\circ\nabla v\,$ by \f{\kp.1.c}). The Killing
field $\,u=J(\navp)\,$ thus would vanish identically on $\,M\,$ (see Remark
\a\kp.1), contradicting nonconstancy of $\,\vp$.
\endremark
\proclaim{Lemma \a\kp.5}Let\/ $\,\vp:M\to\bbR\,$ be a \kip\ on a K\"ahler
manifold\/ $\,(M,g)$. If we set\/ $\,v=\navp$, $\,u=Jv=J(\navp)$,
$\,\bz=\nabla d\vp\,$ and
$$\ba(w,w')\,=\,\bz(Jw,w')\qquad\text{\rm for\ all}\quad
w,w'\in T_xM\quad\text{\rm and}\quad x\in M\,,\ff\kp.2$$
then\/ $\,\ba\,$ is a differential\/ $\,2$-form on\/ $\,M$, and
$$2\hskip.2pt\ba\,=\,d\hs\xi\,,\qquad\text{\rm where}\quad
\xi\,=\,\imath_v\om\,=\,\imath_ug\,.\ff\kp.3$$
\endproclaim
In fact, $\,\imath_v\om=\imath_ug\,$ by \f{\pr.3}, \f{\kr.1}, while Lemma
\a\kp.2(iii) implies \sky\ of $\,\ba(w,w')\,$ in $\,w,w'$. Also, \f{\pr.1} and
\f{\pr.9.ii} give $\,(d\hs\xi)(w,w')=g(Aw,w')\,$ with
$\,A=\nabla u-[\nabla u]^*$, for any $\,C^1$ vector field $\,u,w,w'$ and
$\,\xi=\imath_ug=g(u,\,\cdot\,)\,$ (cf. \f{\pr.3}). As $\,\vp\,$ is a \kip,
$\,u=Jv$, with $\,v=\navp$, is a Killing field, so that $\,A=2\hs\nabla u\,$
and $\,(d\hs\xi)(w,w')=2g(\nabla_{\!w}u,w')=2g(J\nabla_{\!w}v,w')
=-\,2g(\nabla_{\!w}v,Jw')=-\,2\bz(w,Jw')=2\hskip.2pt\ba(w,w')\,$ in view of
\f{\pr.5} with $\,\bz=\nabla d\vp\,$ and \f{\kp.2}.{\hfill\qd}
\medskip
For any \kip\ $\,\vp\,$ on a K\"ahler manifold $\,(M,g)$, we have (see
\cite{\cal})
$$2\hs\imath_v\ri\,\,=\,-\hs dY\,,\qquad
\text{\rm where}\quad v\,=\,\navp\quad\text{\rm and}\quad
Y\,=\,\dv v\,=\,\Delta\vp\hs,\ff\kp.4$$
$\ri\,\,$ being the Ricci tensor. In fact,
$\,J_k^pu^k=-\hs v^p$ (i.e., $\,Ju=-\hs v$) for the Killing field $\,u=Jv$,
while $\,J_l^pu_p=\vp_{,l}$ since
$\,w^lJ_l^pu_p=g(Jw,u)=-\hs g(w,Ju)=g(v,w)=d_w\vp=w^l\vp_{,l}$ for any vector
$\,w$. As $\,\nabla J=0\,$ and $\,u_{p,k}+u_{k,p}=0\,$ (see Remark \a\kp.1),
we have $\,u^k{}_{,k}=\dv u=0\,$ and
$\,\vp_{,lk}{}^k=(J_l^pu_p)_{,k}{}^k=J_l^pu_{p,k}{}^k=-\hs J_l^pu_{k,p}{}^k$.
This in turn equals $\,-\hs J_l^p\ri_{pk}u^k=\rho_{kl}u^k
=-\hs\rho_{lk}u^k=\ri_{lp}J_k^pu^k=-\hs\ri_{lp}v^p$ (cf. \f{\kr.3}), since the
coordinate form of identity \f{\pr.11.b}, with $\,u^k{}_{,k}=0$, gives
$\,\,\ri_{jk}u^k=u^k{}_{,jk}$. Hence $\,\dv\bz\,=\,-\hs\imath_v\ri$, with
$\,v,\bz\,$ as in \f{\pr.11.c}, which, by \f{\pr.11.c}, yields \f{\kp.4}.

\head\S\ck. Conformally-Einstein K\"ahler metrics\endhead
Let $\,\,\ri,\hskip1.5pt\tilde{\ri}\,\,$ and
$\,\,\sc,\hskip1.5pt\tilde{\sc},\,\,$ be the Ricci tensors and scalar
curvatures of conformally related Riemannian metrics $\,g\,$ and
$\,\tilde g=g/\vpsq$ in real dimension $\,n$. Then (see, e.g., \cite{\sdk, p.
411}), with $\,\nabla\,$ and $\,\Delta=g^{jk}\nabla_{\!j}\nabla_{\!k}$
denoting the $\,g$-gradient and $\,g$-Laplacian,
$$\aligned
\tilde{\ri}\,&=\,\,\ri\,\,+\,(n-2)\,\vp^{-1}\nabla d\vp\,+\,
\left[\vp^{-1}\Delta\vp\,-\,(n-1)\vp^{-2}Q\right]g\,,\\
\tilde{\sc}\,&=\,\,\sc\hs\vpsq\,+\,2(n-1)\,\vp\hskip.4pt\Delta\vp\,
-\,n(n-1)Q\,,\qquad\text{\rm where}\quad Q\,=\,g(\navp,\navp)\hs.
\endaligned\ff\ck.1$$
Thus, for $\,n\ge3$, the metric $\,\tilde g=g/\vpsq$ is Einstein if and only if
$$\nabla d\vp\,+\,(n-2)^{-1}\vp\hskip1pt\ri\,\,=\,\fy g\,,\qquad\text{\rm with}
\quad n\fy\,=\,\Delta\vp\,+\,(n-2)^{-1}\hs\sc\hs\vp/n\,.\ff\ck.2$$
\proclaim{Lemma \a\ck.1}Let\/ $\vp\,$ be a\/ $\,C^\infty$ function on a
K\"ahler manifold\/ $\,(M,g)\,$ such that, denoting\/ $\,\,\ri\,\,$ the Ricci
tensor, we have
$$\nabla d\vp\,\,+\,\,\chi\,\ri\,\,=\,\fy g\qquad\text{\rm for\hskip6ptsome}
\hskip9pt
C^\infty\hskip7pt\text{\rm functions}\hskip8pt\chi,\fy\hs.\ff\ck.3$$
Then\/ $\,\vp\,$ is a \kip, i.e.,
$\,u=J(\navp)\,$ is a Killing field on\/ $\,(M,g)$, and
$$\alignedat2
&\text{\rm\ptmi i)}\quad&&
d\hs\xi\,=\,2\,[\fy\hs\om\,
-\,\chi\hskip.5pt\rho\hskip.4pt]\,,\qquad\text{\rm where}\quad
\xi\,=\,\imath_v\om\quad\text{\rm and}\quad v\,=\,\navp\,,\hskip20pt\\
&\text{\rm ii)}\quad&&
d\fy\wedge\hs\om\,=\,d\chi\wedge\rho\hs.\endalignedat\ff\ck.4$$
\endproclaim
In fact, $\,\vp\,$ is a \kip\ in view of Lemma \a\kp.2(iii), as Hermitian
symmetry of $\,\nabla d\vp\,$ follows from that of $\,g\,$ and $\,\,\ri\,\,$
via \f{\ck.3} (cf. \f{\kr.1} -- \f{\kr.2}); now \f{\ck.4.i} is clear from
\f{\kp.3} since, for $\,\bz=\nabla d\vp$, the $\,2$-form $\,\ba\,$ with
\f{\kp.2} equals $\,\fy\hs\om\,-\,\chi\hskip.5pt\rho\,$ (by \f{\ck.3} and
\f{\kr.1}); and, applying $\,d\,$ to \f{\ck.4.i}, we obtain
\f{\ck.4.ii}, as $\,d\hs\om\hs=\hs d\rho=0$.{\hfill\qd}
\medskip
According to Lemma \a\ck.1, condition \f{\ck.3} imposed on a $\,C^\infty$
function $\vp\,$ on a K\"ahler manifold $\,(M,g)\,$ guarantees that $\,\vp\,$
is a \kip. A few special cases of \f{\ck.3} are of independent interest and
have been studied extensively. First, \f{\ck.3} holds if $\,\vp\,$ is a
constant function on a Riemannian manifold $\,(M,g)\,$ which is Einstein or
$\,2$\diml; more generally, for any $\,g,\vp\,$ for which $\,\tilde g=g/\vpsq$
is Einstein, we have \f{\ck.3} in the form of \f{\ck.2}; next, K\"ahler
metrics on compact manifolds with functions $\,\vp\,$ satisfying \f{\ck.3} for
{\it constants\/} $\,\chi,\fy\,$ such that $\,\chi\fy>0\,$ are known as {\it
K\"ahler-Ricci solitons\/} (see \cite{\ptv}, \cite{\tiz}, \cite{\cao});
finally, Riemannian (or pseu\-do-Riem\-ann\-i\-an) manifolds $\,(M,g)\,$
admitting nonconstant $\,C^\infty$ functions $\,\vp\,$ that satisfy \f{\ck.3}
with $\,\chi=0\,$ have been studied extensively, and their local structure is
completely understood (cf. \cite{\khn} and the references therein). See also
the comment following \f{\sr.1} in \S\sr.
\remark{Remark \a\ck.2}Discussing conformally-Einstein metrics, we always
assume, as in \f{\id.1}, that $\,\vp\,$ is a $\,C^\infty$ function on a given
Riemannian manifold $\,(M,g)\,$ such that the metric $\,\tilde g=g/\vpsq$,
{\it\hskip1ptdefined wherever\/} $\,\vp\ne0$, is Einstein. Thus, $\,\vp\,$ may
still vanish somewhere in $\,M$. Although relation \f{\ck.2} and its
consequences, Lemma \a\ck.3 and \f{\ck.5} below, are directly established only
in the open set where $\,\vp\ne0$, they are automatically valid {\it at every
point of\/} $\,M$. In fact, the equalities in question hold both in the set
where $\,\vp\ne0$, and in the interior of the pre\-im\-age $\,\vp^{-1}(0)$,
while the union of these two open sets is clearly dense in $\,M$.
\endremark
\proclaim{Lemma \a\ck.3}Let\/ $\,\mgmt\,$ satisfy\/ \f{\id.1} with\/
$\,m\ge2$. Then\/ $\,\vp\,$ is a \kip\ on\/ $\,(M,g)$, that is,
$\,u=J(\navp)\,$ is a Killing field.
\endproclaim
This is clear from Lemma \a\ck.1, since \f{\ck.2} is a special case of
\f{\ck.3}.{\hfill\qd}
\medskip
Assuming \f{\id.1} with $\,m\ge2\,$ and setting $\,Q=g(\navp,\navp)$,
$\,Y=\Delta\vp$, we have
$$\alignedat2
&\text{\rm a)}\quad&&
(2m-1)(m-2)\,dY\,+\,(m-1)\vp\hs d\hs\sc\,\,
-\,\,\sc\,d\vp\,=\,0\hs,\hskip24pt\\
&\text{\rm b)}\quad&&
\vpsq\hskip.4ptd\hs\sc\hs+2Y d\vp+2(m-1)\vp\hs dY-2m\hs dQ=0.
\endalignedat\ff\ck.5$$
In fact, applying $\,\imath_v$ for $\,v=\navp\,$ (or, $\,d$) to \f{\ck.2} (or,
to $\,\,\tilde{\sc}\,\,$ in \f{\ck.1}), with $\,n=2m$, and using \f{\pr.6.ii},
\f{\kp.4} and \f{\pr.4} (or, constancy of $\,\,\tilde{\sc}\hs$, which is the
scalar curvature the Einstein metric $\,\tilde g$), we obtain
$\,\,\sc\hs\vp\hskip.3ptd\vp+2(m-1)Y d\vp
+m\hs\vp\hs dY-2m(m-1)\hs dQ=0\,$ or, respectively,
$\,2\hskip1.2pt\sc\hs\vp\hskip.3ptd\vp\hs+\hs\vpsq\hskip.4ptd\hs\sc\,+
\hs2(2m-1)\hs[\hs Y\hs d\vp\hs+\hs\vp\hs dY\hs]\hs-\hs2m(2m-1)\hs dQ=0$.
Adding $\,1-2m\,$ times the former equality to $\,m-1\,$ times the latter (or,
subtracting twice the former from the latter), we now get \f{\ck.5.a} (or,
\f{\ck.5.b}). (To obtain \f{\ck.5.a}, we cancelled a $\,\,\vp\,$ factor; this
is allowed as \f{\ck.5.a} still holds in the interior of the pre\-im\-age
$\,\vp^{-1}(0)$, cf. Remark \a\ck.2.)
\proclaim{Proposition \a\ck.4}Suppose that\/ $\,\mgmt\,$ satisfy\/ \f{\id.1}
with\/ $\,m\ge3$. Then, locally in the open set where\/ $\,d\vp\ne0$, each of
the three functions\/ $\,\,\sc\hs$, $\,\,g(\navp,\navp)$, $\,\,\Delta\vp\,$ is
a function of\/ $\,\vp$. Here\/ $\,\,\sc\hs\,$ stands, as usual, for the
scalar curvature of\/ $\,g$.
\endproclaim
In fact, $\,d\,$ applied to \f{\ck.5.a} (or, \f{\ck.5.b}) gives
$\,d\vp\wedge\hs d\hs\sc\,=0\,$ (or, $\,d\vp\wedge\hs dY=0$, as
$\,d\vp\wedge\hs d\hs\sc\,=0$). Hence, applying $\,d\vp\wedge\,...\,\,$ to
\f{\ck.5.b} we get $\,d\vp\wedge\hs dQ=0$.{\hfill\qd}
\medskip
A weaker assertion for $\,m=2\,$ is well-known (see \cite{\sdk, Proposition 3
on p. 416}):
\proclaim{Proposition \a\ck.5}Let\/ $\,\mgmt\,$ satisfy\/ \f{\id.1} with\/
$\,m=2$. Then the scalar curvature\/ $\,\,\sc\,\,$ is a constant multiple of\/
$\,\vp$. Thus, either\/ $\,\,\sc\,\,$ is identically zero, or $\,\hs\vp\,$ is
a constant multiple of\/ $\,\,\sc\hs$.
\endproclaim
This is clear as \f{\ck.5.a} with $\,m=2\,$ gives
$\,\vp\hs d\hs\sc\,-\,\sc\,d\vp=0$, so that $\,\vp^{-1}\hs\sc\,\,$ is constant
on each connected component $\,\,U'$ of the open set $\,\,U\subset M\,$ on
which $\,\vp\ne0$. Thus, in view of Lemma \a\ck.3 and denseness of $\,\,U\,$
in $\,M\,$ (immediate from Remark \a\kp.4), $\,J(\navp)\,$ and
$\,J(\nabla\sc\hs)\,$ are Killing fields; due to their unique continuation
property (Remark \a\kp.1), the proportionality relation
$\,J(\nabla\sc\hs)=pJ(\navp)\,$ with $\,p=\vp^{-1}\hs\sc\hs$, valid on
$\,\,U'$, must hold everywhere in $\,M$. Hence $\,\vp^{-1}\hs\sc\,\,$ is
constant on $\,\,U$.{\hfill\qd}
\medskip
Using well-known analogues of \f{\ck.1} for $\,\dv W\,$
(the divergence of the Weyl tensor of a given metric $\,\,g$) and
$\,\tilde\dv\tilde W\,$ (its counterpart for $\,\,\tilde g=g/\vpsq$), one
verifies that, for a K\"ahler metric $\,g\,$ in complex dimension $\,m$, with
the K\"ahler form $\,\om$, we have $\,(\tilde\dv\tilde W)\hskip.3pt\om=0\,$ if
and only if $\,2(2m-1)(m-2)\hs\imath_v\ri\,
=(m-1)\vp\hs d\hs\sc\,-\,\sc\,d\vp$. This is a generalization of
equation \f{\ck.5.a}: if $\,\tilde g\,$ is Einstein,
$\,\tilde\dv\tilde W\,$ must vanish, while, by \f{\kp.4},
$\,2(2m-1)(m-2)\hs\imath_v\ri\,=-\hs(2m-1)(m-2)\,dY$.

\head\S\sr. Special K\"ahler-Ricci potentials\endhead
We call $\,\vp\,$ a {\it special K\"ahler-Ricci potential\/} on a K\"ahler
manifold $\,(M,g)\,$ if
\vskip4pt
\settabs\+\noindent&\f{\sr.1}\hskip11pt&\cr
\+&&$\vp\,$ is a nonconstant Killing potential on $\hs(M,g)\hs$ (\S\kp) such
that, \hskip.6ptat every\cr
\+&\f{\sr.1}&point with $\,\,d\vp\ne0$, \ all nonzero tangent vectors
orthogonal to $\,\,v=\navp$\cr
\+&&and $\,\,u\hs=\hs Jv\,\,$ are\hs\hs\ eigenvectors\hs\hs\ of\hs\hs\ both
$\,\,\nabla d\vp\,\,$ and\hs\hs\ the\hs\hs\ Ricci\hs\hs\ tensor
$\hs\,\,\ri\hs$.\cr
\vskip4pt
\noindent Examples are described in \S\xm\ and Corollary \a\mc.3 below.
Moreover, \f{\sr.1} is closely related to \f{\ck.3}: assuming \f{\sr.1} we get
\f{\ck.3} on a suitable open set (Remark \a\sr.4), while \f{\ck.3} plus some
strong additional conditions gives \f{\sr.1} (Corollary \a\mc.2).

It is because of this relation with \f{\ck.3} that we use the phrase
`K\"ahler-Ricci potential', adding the word `special' to indicate an extra
assumption (the ``eigenvector clause'') made in \f{\sr.1}. For more on special
K\"ahler-Ricci potentials, see \cite{\dmr}.
\remark{Remark \a\sr.1}For a distribution $\,\Cal V\,$ on a Riemannian
manifold $\,(M,g)\,$ and a symmetric twice-co\-var\-i\-ant tensor $\,\bz\,$
at a point $\,x\in M$, consider the conditions
\vskip4pt
\hbox{\hskip-.8pt
\vbox{\hbox{\f{\sr.2}}\vskip3.2pt}
\hskip16pt
\vbox{
\hbox{a)\hskip9ptAll nonzero vectors in $\,\Cal V_x$ and
$\,\Cal H_x=\Cal V_x^\perp$ are eigenvectors of $\,\bz$.}
\vskip1pt
\hbox{b)\hskip9ptAll nonzero vectors in $\,\Cal H_x=\Cal V_x^\perp$ are
eigenvectors of $\,\bz$.}}}
\vskip4pt
\noindent Let $\,\,\Cal V\,$ now be a $\,J${\it-invariant\/} distribution of
complex dimension one on a K\"ahler manifold $\,(M,g)$. For those symmetric
twice-co\-var\-i\-ant tensors $\,\bz\,$ at a point $\,x\in M\,$ which are also
Hermitian (see \f{\kr.2}), condition \f{\sr.2.b} then implies \f{\sr.2.a}. In
fact, the operator $\,B:T_xM\to T_xM\,$ with $\,\bz(w,w')=g(Bw,w')\,$ for all
$\,w,w'\in T_xM\,$ is self-adjoint, commutes with $\,J$, and
$\,B\Cal V_x^\perp\subset\Cal V_x^\perp$. Hence $\,B\Cal V_x\subset\Cal V_x$.
Choosing $\,v\in\Cal V_x\smallsetminus\{0\}\,$ and $\,\la\in\bbR\,$ with
$\,Bv=\la v$, we thus have $\,Bu=\la u\,$ for $\,u=Jv\,$ (as $\,BJv=JBv$),
which yields \f{\sr.2.a} since $\,\,\dimr\Cal V_x=2$.
\endremark
\definition{Definition \a\sr.2}Given a special \krp\ $\,\vp:M\to\bbR\,$ on a
K\"ahler manifold $\,(M,g)$, as in \f{\sr.1}, let $\,M'\subset M\,$
be the open set on which $\,d\vp\ne0$, and let the vector fields $\,v,u\,$ on
$\,M\,$ and distributions $\,\Cal H,\Cal V\,$ on $\,M'$ be given by
$$\Cal V\,=\,\,\text{\rm Span}\,\{v,u\}\hskip9pt\text{\rm and}\hskip9pt
\Cal H\,=\,\Cal V^\perp\,,\hskip6pt\text{\rm with}\hskip6pt v\,=\,\navp
\hskip6pt\text{\rm and}\hskip6pt u\,=\,Jv\,.\ff\sr.3$$
The {\it septuple $\,(Q,Y,\hs\sc\hs,\si,\ta,\la,\my)\,$ of functions on\/
$\,M'$ associated with\/} $\,\vp\,$ consists of $\,Q=g(\navp,\navp)$,
$\,Y=\Delta\vp$, the scalar curvature $\,\,\sc$, and the ``eigenvalue
functions'' $\,\si,\ta,\la,\my\,$ with $\,\si=\la=0\,$ when $\,\dimc M=1$,
and, in general,
$$\alignedat2
&\ri\,\,=\,\la\hs g\quad\text{\rm on}\quad\Cal H\hs,\qquad
&&\ri\,\,=\,\my\hs g\quad\text{\rm on}\quad\Cal V\hs,\\
&\nabla d\vp\,=\,\si g\quad\text{\rm on}\quad\Cal H\hs,\qquad
&&\nabla d\vp\,=\,\ta g\quad\text{\rm on}\quad\Cal V\hs,\hskip45pt\\
&\ri\hskip1pt(\Cal H,\Cal V)\,=\,(\nabla d\vp)(\Cal H,\Cal V)\,=\,&&\{0\}
\qquad\text{\rm for}\quad\Cal H,\Cal V\quad\text{\rm as\ in}\quad
\text{\rm\f{\sr.3}.}\endalignedat\ff\sr.4$$
The last line states that $\,\Cal H,\Cal V\,$ are $\,\,\ri$-or\-thog\-o\-nal
and $\,\nabla d\vp$-or\-thog\-o\-nal to each other (cf. the eigenvector clause
of \f{\sr.1}). The rest of \f{\sr.4} now follows, for some $\,C^\infty$
functions $\,\la,\my,\si,\ta\,$ on $\,M'$, as \f{\sr.2.b} gives \f{\sr.2.a}
(Remark \a\sr.1). Also, by \f{\sr.3},
$$g(v,v)\,=\,g(u,u)\,=\,Q\,,\qquad g(v,u)\,=\,0\qquad\text{\rm everywhere\ in}
\quad M.\ff\sr.5$$
\enddefinition
\remark{Remark \a\sr.3}Conversely, \f{\sr.1} obviously holds for any
nonconstant Killing potential $\,\vp\,$ on a K\"ahler manifold which satisfies
\f{\sr.4} on the open set where $\,d\vp\ne0$, with $\,\Cal V,\Cal H\,$ given
by \f{\sr.3} and some $\,\si,\ta,\la,\my$.
\endremark
\remark{Remark \a\sr.4}Let $\,\vp\,$ be a\/ $\,C^\infty$ function on a
Riemannian manifold $\,(M,g)\,$ such that some distribution $\,\Cal V\,$
on $\,M\,$ satisfies \f{\sr.2.a} at every $\,x\in M$, both for
$\,\bz=\nabla d\vp\,$ and the Ricci tensor $\,\bz=\,\ri\hs$. Then \f{\ck.3}
holds in the open set of all points at which $\,\,\ri\,\,$ is not a multiple
of $\,g\,$ (i.e., $\,\,\ri\,\ne\,\sc\hs g/n$, with $\,n=\,\dimr M$).

In fact, the space of all symmetric twice-co\-var\-i\-ant tensors $\,\bz\,$
with \f{\sr.2.a} at a fixed $\,x\,$ is two\diml, with the eigenvalues serving
as parameters. Hence $\,\nabla d\vp\,$ is a combination of $\,\,\ri\,\,$ and
$\,g\,$ at points where they are linearly independent.

Thus, every special \krp\ $\,\vp\,$ on a K\"ahler manifold satisfies \f{\ck.3}
on the set where $\,\,\ri\,\ne\,\sc\hs g/n\,$ and $\,d\vp\ne0\,$ (by
\f{\sr.1}, as \f{\sr.4} then implies \f{\sr.2.a}).
\endremark
\proclaim{Lemma \a\sr.5}Let\/ $\,(Q,Y,\hs\sc\hs,\si,\ta,\la,\my)\,$ be the
septuple associated, as above, with a special \krp\/ $\,\vp\,$ on a K\"ahler
manifold $\,(M,g)\,$ of complex dimension\/ $\,m\ge1$. Then, in the open set\/
$\,M'$ where\/ $\,d\vp\ne0$, we have\/ $\,dQ=2\ta\,d\vp$,
$\,dY=-\hs2\my\,d\vp\,$ and\/
$\,(m-2)\hs\nabla\la=-\hs[\nabla\la]\vrt+2(m-1)(\my-\la)\hs\si\hs v/Q$,
where\/ $\,v=\navp\,$ and\/ $\,...\vrt$ stands for the\/ $\,\Cal V\,$
component relative to the decomposition\/ $\,TM'=\Cal H\oplus\Cal V$, cf.
\f{\sr.3}. Furthermore, if\/ $\,m\ne2\,$ then\/
$\,Q\,d\la=2(\my-\la)\hs\si\,d\vp$, while, for any\/ $\,m\ge1$,
\widestnumber\item{(ii)}\roster
\item"(i)"$\,\,\sc\,\,=\,2\my\,+\,2(m-1)\la\,,\qquad
Y\,=\,2\ta\,+\,2(m-1)\hs\si\,,$
\item"(ii)"$\,\nabla_{\!v}v\,=\,-\hs\nabla_{\!u}u\,=\,\ta\hskip.4ptv\,,
\hskip13pt\nabla_{\!v}u\,=\,\nabla_{\!u}v\,=\,\ta\hskip.4ptu\,,$\hskip13pt
with\/ $\,u=Jv$.
\endroster
\endproclaim
\demo{Proof}As $\,\,\sc\,=\,\text{\rm Trace}_g\hs\ri\,\,$ and
$\,Y=\Delta\vp=\,\text{\rm Trace}_g(\nabla d\vp)\,$ (cf. \f{\pr.8}), \f{\sr.4}
gives (i). Since $\,v\,$ is a section of $\,\Cal V$, \f{\sr.4}, \f{\pr.3} and
\f{\pr.4} give $\,\imath_v\ri\,=\my\,\imath_vg=\my\,d\vp\,$ and
$\,\imath_v(\nabla d\vp)=\ta\hs d\vp$. Now, by \f{\kp.4} and \f{\pr.6.ii},
$\,dY=-\hs2\my\,d\vp\,$ and $\,dQ=2\ta\,d\vp$.

By \f{\sr.4} and Remark \a\pr.1, $\,\nabla_{\!v}v=\ta\hskip.4ptv\,$ and
$\,\nabla_{\!u}v=\ta\hskip.4ptu$, As $\,u=Jv\,$ and $\,\nabla J=0$, this gives
$\,\nabla_{\!u}u=-\hs\ta\hskip.4ptv$. Now (ii) follows since \f{\kp.1.b} and
\f{\pr.2.ii} yield $\,\nabla_{\!v}u=\nabla_{\!u}v$.

Setting $\,\bz=\,d\vp\otimes d\vp+\xi\hskip-.7pt\otimes\xi\,$ on $\,M'$,
with $\,\xi=\imath_ug$, we have $\,\dv\bz=Yd\vp$, by \f{\pr.10.iii} with
$\,\imath_vg=\hs d\vp\,$ (see \f{\pr.4}), or $\,\xi=\imath_ug$, as
$\,\dv v=Y\,$ by \f{\pr.10.i}, while $\,\dv u=0\,$ since $\,u\,$ is a Killing
field (cf. Remark \a\kp.1) and, finally,
$\,\nabla_{\!v}(d\vp)=-\hs\nabla_{\!u}\xi\,=\,\ta\,d\vp\,$ (due to the first
relation in (ii), with $\,v,u\,$ replaced by the corresponding $\,1$-forms
$\,d\vp,\xi$). Also, by \f{\pr.10.iii} and \f{\sr.5},
$\,\imath_v\bz=Q\,d\vp=Q\hs\imath_vg$,
$\,\imath_u\bz=Q\hs\xi=Q\hs\imath_ug$, and $\,\imath_w\bz=0\,$ for vectors
$\,w\,$ orthogonal to both $\,v\,$ and $\,u$. Thus,
$\,\,\ri\,-\la\hs g=\bz'$, with $\,\bz'=(\my-\la)\hs\bz/Q\,$ since, by
\f{\sr.4}, both sides have the same image under $\,\imath_w$ for any tangent
vector $\,w$. Combined with (i), \f{\pr.10.ii} and the relations
$\,\dv\bz=Yd\vp$ and $\,dQ=2\ta\,d\vp$, i.e., $\,\nabla Q=2\ta\hs v$, the
above expressions for $\,\imath_v\bz$, $\,\imath_u\bz\,$ and $\,\imath_w\bz\,$
also show that $\,\dv\bz'=\imath_wg\,$ for the vector field
$\,w=[\nabla\my]\vrt-[\nabla\la]\vrt+2(m-1)(\my-\la)\hs\si\hs v/Q$. However,
$\,\imath_wg=\dv\bz'=\hs\dv(\ri\,-\la\hs g)\,$ (since
$\,\bz'=\,\ri\,-\la\hs g$) which, as
$\,2\hs\dv(\ri\,-\la\hs g)=\hs d\hs\sc\,-\hs2\hs d\la\,$ (by \f{\pr.11.a},
\f{\pr.10.ii} and \f{\pr.4}), gives $\,\imath_wg=d\my+(m-2)\hs d\la\,$ (by
(i)), and so $\,w=\nabla\my+(m-2)\hs\nabla\la\,$ (cf. \f{\pr.4}).

As $\,dY=-\hs2\my\,d\vp$, we have $\,0=\hs ddY=-\hs2\hs d\my\wedge d\vp$, so
that $\,\nabla\my\,$ equals a function times $\,\nabla\vp\,$ and, in
particular, $\,[\nabla\my]\vrt=\nabla\my$. Comparing the above two expressions
for $\,w$, we thus obtain the required formula for
$\,(m-2)\hs\nabla\la$. When $\,m>2$, this shows that $\,\nabla\la\,$ is a
section of $\,\Cal V=\,\text{\rm Span}\,\{v,u\}$, i.e.,
$\,[\nabla\la]\vrt=\nabla\la$, and hence
$\,\nabla\la=2(\my-\la)\hs\si\hs v/Q$, that is,
$\,Q\,d\la=2(\my-\la)\hs\si\,d\vp$. Since both sides of the last equality
vanish when $\,m=1\,$ (Definition \a\sr.2), this completes the
proof.{\hfill\qd}
\enddemo

\head\S\xm. Examples\endhead
To describe examples of special \krp s $\,\vp\,$ on K\"ahler manifolds
$\,(M'\!,g)\,$ (see \f{\sr.1}), we assume that the following data are given:
$$\jy,r,\theta,\vp,\fe\hs;\qquad a,\ve,\y\hs;\qquad m,N,h\hs;\qquad
\Cal L,\proj,\Cal V\hs;\qquad\Cal H,\langle\,,\rangle,M'\nh.\hskip10pt
\ff\xm.1$$
Here $\,\jy\subset(0,\infty)\,$ is an open interval, $\,r\in\jy\,$ is a real
variable, $\,a\in\bbR\smallsetminus\{0\}\,$ and $\,\ve\in\{-\hs1,0,1\}\,$ are
constants, while $\,\theta,\vp,\fe:\jy\to\bbR\,$ are $\,C^\infty$ functions
such that $\,\theta>0$, $\,\fe>0$, $\,d\vp/dr=ar\theta$, and either
$\,\ve=0\,$ and $\,\fe=1$, or $\,\ve=\pm\hs1\,$ and $\,\fe=2\ve\hs(\vp-\y)\,$
with a constant $\,\y$. Furthermore, $\,m\ge1\,$ is an integer and $\,(N,h)\,$
is a K\"ahler manifold of complex dimension $\,m-1$, which we assume to be
Einstein unless $\,m=2$. Also, $\,\Cal L\,$ is a holomorphic line bundle with
a Hermitian fibre metric $\,\langle\,,\rangle\,$ over $\,N$, while $\,\proj\,$
and $\,\Cal V\,$ denote the bundle projection $\,\Cal L\to N\,$ and the
vertical distribution in $\,\Cal L\,$ (\S\cc). Next, $\,\Cal H\,$ is the
horizontal distribution of a linear connection in $\,\Cal L\,$ making
$\,\langle\,,\rangle\,$ parallel, whose curvature form $\,\varOmega\,$ (Remark
\a\cc.1) equals $\,-\hs2\ve a\,\omh$, where $\,\omh$ is the K\"ahler form of
$\,(N,h)$. We also assume that $\,\Cal H\,$ is $\,J$-invariant as a subbundle
of $\,T\Cal L$, where $\,J:T\Cal L\to T\Cal L\,$ denotes the complex structure
tensor of $\,\Cal L$. (Cf. Remark \a\lc.4.) Finally, $\,M'$ is a connected
open subset of $\,\Cal L\smallsetminus N\,$ contained in $\,r^{-1}(\jy)\,$
(the $\,r$-pre\-im\-age of $\,\jy$), where the symbol $\,r\hs$ stands for
the norm function of $\,\langle\,,\rangle\,$ as well (see Remark \a\cc.2), so
that $\,\fe,\theta,\vp\,$ and other $\,C^\infty$ functions of $\,r\in\jy\,$
may be treated as $\,C^\infty$ functions $\,M'\to\bbR\hs$.

A metric $\,g\,$ on $\,M'$ now is defined by requiring that
$\,g(\Cal H,\Cal V)=\{0\}$, i.e., $\,\Cal H\,$ be $\,g$-or\-thog\-o\-nal to
$\,\Cal V$, while $\,g=\fe\hs\proj^*\nh h\,$ on $\,\Cal H\,$ and
$\,g=\theta\,\text{\rm Re}\hskip1pt\langle\,,\rangle\,$ on $\,\Cal V$, where
$\,\,\text{\rm Re}\hskip1pt\langle\,,\rangle\,$ is the standard Euclidean
metric on each fibre of $\,\Cal L$.

Our $\,M'$, being an open submanifold of $\,\Cal L$, is a complex manifold of
complex dimension $\,m\ge1$. We will verify later, in \S\dx, that $\,g\,$ is a
K\"ahler metric on $\,M'$, and $\,\vp\,$ is a special \krp\ on $\,(M'\!,g)$,
as defined in \f{\sr.1}.

\head\S\mc. More on the conformally-Einstein case\endhead
\proclaim{Lemma \a\mc.1}Let a nonconstant\/ $\,C^\infty$ function\/
$\,\vp:M\to\bbR\,$ on a K\"ahler manifold\/ $\,(M,g)\,$ of any complex
dimension $\,m\ge2\,$ satisfy\/ \f{\ck.3} with\/ $\,\chi,\fy:M\to\bbR\,$ such
that\/ $\,d\fy\wedge\hs d\vp=\hs d\chi\wedge\hs d\vp=0$, i.e., any point
with\/ $\,d\vp\ne0\,$ has a neighborhood on which both\/ $\,\chi,\fy\,$ are
functions of\/ $\,\vp$. Also, let\/ $\,\Cal V,\Cal H\,$ be given by\/
\f{\sr.3} on the open set\/ $\,\,M'\subset M\,$ where\/ $\,d\vp\ne0$.
Then, in the open subset of\/ $\,M'$ on which\/ $\,d\chi\ne0$, we have\/
\f{\sr.4} for some\/ $\,C^\infty$ functions\/ $\,\la,\my,\si,\ta\,$ with\/
$\,d\fy\,=\,\la\,d\chi$.
\endproclaim
\demo{Proof}Let $\,\,\bz\,=\,(d\chi/d\vp)\,\ri\,\,-\,(d\fy/d\vp)\hs g\,\,$ on
$\,M'$, with $\,\imath_w$ as in \f{\pr.3}. (Note that $\,\bz\,$ is
well-defined: $\,d\chi/d\vp\,$ and $\,d\fy/d\vp\,$ make sense as
$\,d\fy\wedge\hs d\vp=\hs d\chi\wedge\hs d\vp=0$.) Then
$\,\Cal H_x\subset\{w\in T_xM:\imath_w\bz=0\}\,$ at every $\,x\in M'$.

In fact, we may fix $\,x\in M'$ and assume that $\,\bz(x)\ne0$. For
$\,\ba=(d\chi/d\vp)\,\rho\,-\,(d\fy/d\vp)\,\om\hs$,
\f{\ck.4.ii} gives $\,\ba\wedge\hs d\vp=0$, and so, locally,
$\,\ba=\xi\wedge\hs d\vp\,$ for some $\,1$-form $\,\xi$. As $\,\ba\,$ and
our $\,\bz\,$ satisfy \f{\kp.2} (cf. \f{\kr.1}), $\,\xi\,$ and $\,d\vp\,$ are
linearly independent at $\,x$. By \f{\pr.3} and \f{\pr.9.i},
$\,\imath_w\ba=\xi(w)\hs d\vp-\hs(d_w\vp)\xi\,$ for any vector $\,w$. Thus,
the nullspace $\,\,\text{\rm Ker}\,\ba(x)=\{w\in T_xM:\imath_w\ba=0\}\,$
is the intersection of the kernels of $\,\xi\,$ and $\,d\vp\,$ at $\,x$, i.e.,
$\,[\hs\text{\rm Ker}\,\ba(x)]^\perp=\,\text{\rm Span}\,\{v(x),u'(x)\}$,
where $\,v=\navp\,$ and $\,u'$ is the vector field with
$\,\xi=g(u',\,\cdot\,)$. However, \f{\kr.2} (which holds for our $\,\bz$, as
it does for $\,g,\,\ri$), along with \f{\pr.3} and \f{\kp.2}, implies that
$\,\,\text{\rm Ker}\,\ba(x)=\,\text{\rm Ker}\,\bz(x)\,$ (the nullspace of
$\,\bz(x)$). As $\,\bz\,$ is Hermitian (cf. \f{\kr.2}), both
$\,\,\text{\rm Ker}\,\bz(x)\,$ and
$\,[\hs\text{\rm Ker}\,\bz(x)]^\perp$ are $\,J$-invariant. Since $\,v(x)
\in[\hs\text{\rm Ker}\,\bz(x)]^\perp=\,\text{\rm Span}\,\{v(x),u'(x)\}$, we
have $\,[\hs\text{\rm Ker}\,\bz(x)]^\perp=\,\text{\rm Span}\,\{v(x),u(x)\}\,$
with $\,u=Jv$, as required.

As all nonzero vectors in $\,\Cal H\,$ are eigenvectors of
$\,\bz=(d\chi/d\vp)\,\ri\,\,-\,(d\fy/d\vp)\hs g\,$ for the eigenvalue $\,0$,
at points $\,x\in M'$ where $\,d\chi\ne0\,$ they are eigenvectors of
$\,\,\ri\,\,$ for the eigenvalue $\,\la(x)$, with
$\,\la=(d\fy/d\vp)/(d\chi/d\vp)=\,d\fy/d\chi\,$ (i.e.,
$\,d\fy\,=\,\la\,d\chi$). Now \f{\ck.3} implies the same for
$\,\,\nabla d\vp\,$ (with some eigenvalue $\,\si(x)$) instead of $\,\,\ri\,\,$
and $\,\la(x)$. Consequently, $\,\bz=\,\ri\,$ and $\,\bz=\nabla d\vp\,$ both
satisfy \f{\sr.2.b}, and hence (by Remark \a\sr.1) also \f{\sr.2.a}, at every
such $\,x$, which completes the proof.{\hfill\qd}
\enddemo
\proclaim{Corollary \a\mc.2}Let a nonconstant\/ $\,C^\infty$ function\/
$\,\vp:M\to\bbR\,$ on a K\"ahler manifold\/ $\,(M,g)\,$ of complex dimension\/
$\,m\ge2\,$ satisfy\/ \f{\ck.3} with\/ $\,\chi,\fy:M\to\bbR\,$
such that\/ $\,d\fy\wedge\hs d\vp=\hs d\chi\wedge\hs d\vp=0\,$ and\/
$\,d\chi\ne0\,$ wherever\/ $\,d\vp\ne0$. Then\/ $\,\vp\,$ is a special \krp\
on\/ $\,(M,g)$, i.e., satisfies\/ \f{\sr.1} as well.
\endproclaim
This is clear from Lemmas \a\ck.1 and \a\mc.1 along with Remark
\a\sr.3.{\hfill\qd}
\proclaim{Corollary \a\mc.3}Condition\/ \f{\id.1} with\/ $\,m\ge3$, or\/
\f{\id.2} with\/ $\,m=2$, implies\/ \f{\sr.1}.
\endproclaim
In fact, \f{\ck.2} and Propositions \a\ck.4, \a\ck.5 yield the hypotheses of
Corollary \a\mc.2.{\hfill\qd}
\remark{Remark \a\mc.4}In terms of the septuple
$\,(Q,Y,\hs\sc\hs,\si,\ta,\la,\my)\,$ associated with a special \krp\
$\,\vp\,$ on a K\"ahler manifold of complex dimension $\,m\ge2\,$ (Definition
\a\sr.2), condition \f{\de.3.iii} below is necessary and sufficient for
$\,\mgmt\,$ to satisfy \f{\id.1}. In fact, by \f{\sr.4}, the first equality in
\f{\ck.2}, with $\,n=2m$, is equivalent to
$\,2(m-1)\hs\si+\vp\la=(n-2)\fy=2(m-1)\ta+\vp\my$.
\endremark

\head\S\de. Differential equations related to \f{\sr.1} and \f{\id.1}\endhead
For a fixed integer $\,m\ge2$, let the system of ordinary differential
equations
$$\aligned
dQ\,&=\,2\ta\,d\vp\,,\qquad dY\,=\,-\hs2\my\,d\vp\,,\qquad
Q\,d\si\,=\,2(\ta-\si)\hs\si\,d\vp\,,\\
Q\,d\ta\,&=\,[2(m-1)(\si-\ta)\hs\si\,-\,\my Q]\,d\vp\,,\qquad
Q\,d\la\,=\,2(\my-\la)\hs\si\,d\vp\,,
\endaligned\ff\de.1$$
$$\alignedat2
&\text{\rm\ptmi i)}\quad&&
2(m-1)(\ta-\si)\hs Q\,d\my\\
&\quad&&\qquad
=\,(\la-\my)\left[\la Q\,+\,(2m-3)\my Q\,
+\,4(m-1)^2(\ta-\si)\hs\si\right]d\vp\,,\\
&\text{\rm ii)}\quad&&
(m-1)(\ta-\si)\,d\hs\sc\,=\,(\la-\my)\hs[\la\,+\,(2m-3)\my]\,d\vp\,,
\hskip20pt\endalignedat\ff\de.2$$
be imposed on seven unknown $\,C^1$ functions
$\,Q,Y,\hs\sc\hs,\si,\ta,\la,\my\,$ of the real variable $\,\vp$, defined on
an unspecified interval. We will also require that, on this interval,
$$\alignedat2
&\text{\rm\ptmii i)}\quad&&
Q\,\ne\,0\quad\text{\rm everywhere,}\\
&\text{\rm\ptmi ii)}\quad&&
Y\,=\,2\ta\,+\,2(m-1)\hs\si\,,\qquad\sc\,\,=\,2\my\,+\,2(m-1)\la\,,\hskip70pt\\
&\text{\rm iii)}\quad&&
2(m-1)(\ta-\si)\,=\,(\la-\my)\vp.\endalignedat\ff\de.3$$
The constraints \f{\de.3.ii}, \f{\de.3.iii} amount to vanishing of specific
{\it integrals\/} of the system \f{\de.1} -- \f{\de.2}; see Lemma \a\de.1
below.

Conditions \f{\de.1} -- \f{\de.3} arise in our discussion of \f{\sr.1} and
\f{\id.1} (see Lemma \a\qs.1 and Corollary \a\qs.2). Certain particularly
simple $\,C^\infty$ solutions to \f{\de.1} -- \f{\de.3} have $\,\si=0\,$
identically (cf. Corollary \a\rr.2(ii), (iii)). These special solutions will
be classified in \S\sz, while in  \S\as\ -- \S\su\ we deal with the opposite
case:
$$\si\,\ne\,0\quad\text{\rm everywhere\ in\ the\ interval.}\ff\de.4$$
\proclaim{Lemma \a\de.1}Let\/ $\,m\ge1\,$ be a fixed integer. Then both\/
$\,2\y=2\hskip.5pt\vp\,-\,Q/\si\,$ and\/ $\,2\ta+2(m-1)\hs\si-Y\,$ are
constant, as functions of\/ $\,\vp$, on any interval on which\/
$\,Q,Y,\hs\sc\hs,\si,\ta,\la,\my\,$ are\/ $\,C^1$ functions
satisfying\/ \f{\de.3.i}, \f{\de.4} and the first four equations in\/
\f{\de.1}. Similarly, relations\/ \f{\de.1} -- \f{\de.3.i}, with\/ $\,m\ge2$,
imply constancy of\/ $\,\hs2\my+2(m-1)\la-\,\sc\,\,$ on any interval on
which\/ $\,\si\ne\ta$, and constancy of\/ $\,Z-\vp$, with\/
$\,Z=2(m-1)(\ta-\si)/(\la-\my)$, on any interval on which\/ $\,\la\ne\my$.
\endproclaim
In fact, the second assertion follows since \f{\de.1} -- \f{\de.3.i} yield,
for $\,m\ge2$,
\vskip4pt
\noindent($*$)\hskip10pt
$Z\hs d\hs[\my+(m-1)\la]/d\vp\hs=\hs\la\hs+\hs(2m-3)\my\,\,$ and
$\,\,dZ/d\vp\hs=\hs1\,\,$ wherever $\,\,\la\ne\my$,
\vskip4pt
\noindent
while the first four equations in \f{\de.1} with $\,m\ge1\,$ give
$\,d\hs[Q/\si]/d\vp=2\,$ wherever $\,\si\ne0\,$ and
$\,d\hs[\ta+(m-1)\hs\si]/d\vp=-\hs\my\,$ wherever $\,Q\ne0$, proving the first
claim.{\hfill\qd}
\medskip
For further integrals of the system \f{\de.1} -- \f{\de.4}, see \S\is.
\remark{Remark \a\de.2}For a septuple $\,(Q,Y,\hs\sc\hs,\si,\ta,\la,\my)\,$ of
$\,C^1$ functions of a variable $\,\vp$, \f{\de.1} and \f{\de.3} with
$\,m\ge1\,$ imply \f{\de.2}. In fact, applying $\,d\,$ to \f{\de.3.iii}, then
multiplying by $\,(\la-\my)Q\,$ and, finally, replacing: $\,(\la-\my)\vp\,$
with $\,2(m-1)(\ta-\si)\,$ (cf. \f{\de.3.iii}), and
$\,d\la,\hs d\si,\hs d\ta\,$ with what is provided by \f{\de.1}, we obtain
\f{\de.2.i}. Now \f{\de.2.ii} is clear from \f{\de.2.i} and the formulae for
$\,d\la\,$ and $\,\,\sc\,\,$ in \f{\de.1}, \f{\de.3.ii}.
\endremark

\head\S\qs. Equations \ \f{\de.1} \ for special K\"ahler-Ricci
potentials\endhead
\proclaim{Lemma \a\qs.1}Let\/ $\,(Q,Y,\hs\sc\hs,\si,\ta,\la,\my)\,$ be the
septuple associated with a special \krp\/ $\,\vp\,$ on a K\"ahler manifold\/
$\,(M,g)\,$ of complex dimension\/ $\,m\ge1$, as in\/ {\rm Definition
\a\sr.2}. Then, in the open set\/ $\,M'\subset M\,$ where\/ $\,d\vp\ne0$,
\widestnumber\item{(b)}\roster
\item"(a)"Each of the five functions\/ $\,Q,Y,\si,\ta,\my\,$ is, locally, a\/
$\,C^\infty$ function of\/ $\,\vp$.
\item"(b)"$\,Q,Y,\hs\sc\hs,\si,\ta,\la,\my\,$ satisfy\/ \f{\de.1.i} --
\f{\de.1.iv}, \f{\de.3.i} and\/ \f{\de.3.ii},
\endroster
where\/ \f{\de.1.i} -- \f{\de.1.v} are the five equations forming\/
\f{\de.1}. If, in addition, $\,m\ne2$, then\/ $\,\,\sc\hs,\la\,$ also are,
locally in\/ $\,M'$, functions of\/ $\,\vp$, and\/ \f{\de.1.v} holds.
\endproclaim
\demo{Proof}On $\,M'$, Lemma \a\sr.5 gives \f{\de.1.i}, \f{\de.1.ii}, which
implies (a) for $\,Q,Y$, and
hence also for $\,\ta,\my\,$ (as $\,2\ta=\hs dQ/d\vp$, $\,2\my=-\hs dY/d\vp$),
and $\,\si\,$ (cf. Lemma \a\sr.5(i); by Definition \a\sr.2, $\,\si=\la=0\,$
when $\,m=1$.) This proves (a) for all five functions, and (combined with
Lemma \a\sr.5) also the final clause.

Let us now set $\,v=\navp\,$ and $\,\bz=\nabla d\vp$,
As \f{\pr.11.c}, \f{\kp.4} and \f{\de.1.ii} give $\,\dv\bz=-\hs\my\,d\vp\,$
for $\,\bz=\nabla d\vp$, \f{\pr.6.i} yields
$\,\imath_v(\dv\bz)=-\hs\my\hs\imath_vd\vp=-\hs\my Q$, while, by \f{\sr.4},
$\,\,\langle\bz,\bz\rangle=2\ta^2+2(m-1)\hs\si^2$, with
$\,\langle\,,\rangle\,$
as in \f{\pr.7}. However,
$\,\dv(\nabla_{\!v}v)=\imath_v(\dv\bz)\,+\,\langle \bz,\bz\rangle\,$ due to
the obvious local-coordinate identity
$\,[\vp^{,k}\vp_{,jk}]^{,j}=\vp^{,k}\vp_{,jk}{}^j+\vp_{,jk}\vp^{,jk}$. Thus,
$\,\dv(\nabla_{\!v}v)=-\my Q+2\ta^2+2(m-1)\hs\si^2$ while, as $\,Y=\Delta\vp\,$
and $\,\nabla_{\!v}v=\ta\hskip.4ptv\,$ (Lemma \a\sr.5(ii)),
$\,\dv(\nabla_{\!v}v)=\dv(\ta\hskip.4ptv)=\,d_v\ta\,+\,\ta\hs Y
=\,d_v\ta+[2\ta+2(m-1)\hs\si]\hs\ta\,$ (by \f{\pr.10.ii}, \f{\pr.10.i} and Lemma
\a\sr.5(i)). Equating the two formulae for $\,\dv(\nabla_{\!v}v)$, we obtain
$\,d_v\ta=2(m-1)(\si-\ta)\hs\si-\my Q$. Both sides of \f{\de.1.iv} thus yield the
same value under $\,\imath_v$ (cf. \f{\pr.6.i}) and, since they both are
functional multiples of $\,d\vp\,$ (by (a)), while $\,\imath_vd\vp=Q\ne0\,$ by
\f{\pr.6.i}, they must coincide. This proves \f{\de.1.iv}. Since
$\,(1-m)\,d\si=\hs d\ta+\my\,d\vp\,$ (as one sees evaluating $\,dY\,$ from
Lemma \a\sr.5(i) and using \f{\de.1.ii}), \f{\de.1.iii} now is immediate from
\f{\de.1.iv}. Finally, \f{\de.3.i} follows since $\,Q=g(\navp,\navp)>0\,$ on
$\,M'$, while \f{\de.3.ii} is obvious from Lemma \a\sr.5(i), which completes
the proof.{\hfill\qd}
\enddemo
\proclaim{Corollary \a\qs.2}Let\/ $\,\mgmt\,$ satisfy\/ \f{\id.1} with\/
$\,m\ge3$, or\/ \f{\id.2} with\/ $\,m=2$, and let\/
$\,(Q,Y,\hs\sc\hs,\si,\ta,\la,\my)\,$ be the septuple associated with\/
$\,\vp$, cf. {\rm Definition \a\sr.2} and\/ {\rm Corollary \a\mc.3}. Then\/
\f{\de.1} -- \f{\de.3} hold in the open set on which\/ $\,d\vp\ne0$.
\endproclaim
In fact, $\,\vp\,$ satisfies \f{\sr.1} (Corollary \a\mc.3). Now Lemma \a\qs.1
and Remark \a\mc.4 give \f{\de.1} and \f{\de.3}, except for the last
formula in \f{\de.1} when $\,m=2$. That formula then follows, however, from
the expression for $\,(m-2)\hs\nabla\la\,$ in Lemma \a\sr.5, in which we now
have $\,[\nabla\la]\vrt=\nabla\la$, as $\,\la\,$ is (locally, at points with
$\,d\vp\ne0$) a function of $\,\vp$, since so are $\,\,\sc\,\,$ (Proposition
\a\ck.5) and $\,\my\,$ (Lemma \a\qs.1), while $\,2(m-1)\la=\,\sc\,-2\my\,$ by
Lemma \a\sr.5(i). Now Remark \a\de.2 yields \f{\de.2}.{\hfill\qd}

\head\S\cm. The $\,\hs\si\hs\,$ alternative\endhead
The following two lemmas are well-known.
\proclaim{Lemma \a\cm.1}Let\/ $\,u\,$ be a Killing field on a Riemannian
manifold\/ $\,(M,g)$, and let\/ $\,y\in M\,$ be a point such that\/
$\,u(y)=0$.
\widestnumber\item{(ii)}\roster
\item"(i)"For every sufficiently small\/ $\,\,\text{\rm d}>0$, the flow
of\/ $\,u\,$ restricted to the radius\/ $\,\,\text{\rm d}\,$ open ball\/
$\,\,U\,$ centered at\/ $\,y\,$ consists of ``global'' isometries\/
$\,\,U\to\,U$.
\item"(ii)"If\/ $\,\,\text{\rm d}\,$ and\/ $\,\,U\,$ in\/ {\rm(i)} are chosen
so that, in addition, the exponential mapping\/ $\,\,\e_y:U'\to\,U\,$ is a
diffeomorphism, where\/ $\,\,U'$ is the radius\/ $\,\,\text{\rm d}\,$ open
ball in\/ $\,T_yM\,$ centered at\/ $\,0$, then\/ $\,u\,$
restricted to\/ $\,\,U\,$ is the $\,\,\e_y$-image of the linear vector field
on\/ $\,\,U'$ given by\/ $\,w\mapsto\nabla_{\!w}u$.
\endroster
\endproclaim
This is clear as the diffeomorphism $\,\,\e_y:U'\to\,U$, for any small
$\,\,\text{\rm d}\hs$, makes every isometry $\,\Phi:U\to\,U\,$ with
$\,\Phi(y)=y\,$ correspond to the linear mapping $\,d\Phi_y:U'\to U'$. Now
(ii) is immediate if we apply this to the isometries forming the flow of
$\,u\,$ and note that their differentials at $\,y\,$ form a one-parameter
group of linear isometries with the infinitesimal generator
$\,(\nabla u)(y)$.{\hfill\qd}
\proclaim{Lemma \a\cm.2}For a Killing vector field\/ $\,u\,$ on a Riemannian
manifold\/ $\,(M,g)$, let\/ $\,N(u)=\{y\in M:u(y)=0\}\,$ be the set of all
zeros of\/ $\,u$. If\/ $\,u\ne0\,$ somewhere in\/ $\,M$, then, for every
connected component\/ $\,N\,$ of\/ $\,N(u)$, with\/ $\,\nabla u\,$ as in\/
\f{\pr.1},
\widestnumber\item{(d)}\roster
\item"(a)"$N\,$ is contained in an open set that does not intersect any other
component.
\item"(b)"$N\subset M\,$ is a closed set and a submanifold with the subset
topology.
\item"(c)"The submanifold\/ $\,N\,$ is totally geodesic in\/ $\,(M,g)\,$ and\/
$\,\dim M-\hs\dim N\ge2$.
\item"(d)"For any\/ $\,y\in N\,$ we have
$\,T_yN=\,\text{\rm Ker}\,[(\nabla u)(y)]=\{w\in T_yM:\nabla_{\!w}u=0\}$.
\endroster
Furthermore, the set\/ $\,M'=M\smallsetminus N(u)\,$ is connected, open and
dense in\/ $\,M$.
\endproclaim
In fact, (a) -- (d) (except for the inequality in (c)) are immediate from
Lemma \a\cm.1(ii). Now let us fix $\,y\in M\,$ with $\,u(y)=0\,$ and
$\,\,\text{\rm d}\hs,\,U'\!,\,U\,$ with (i), (ii) of Lemma \a\cm.1. The space
$\,V=\,\text{\rm Ker}\,[(\nabla u)(y)]\subset T_yM\,$ has $\,\dim V<\dim M$,
or else $\,u\,$ and $\,\nabla u\,$ would both be zero at $\,y$, i.e., $\,u\,$
would vanish identically (Remark \a\kp.1). Next, choosing
$\,w\in V^\perp\subset T_yM\,$ with $\,w\ne0$, we have
$\,\nabla_{\!w}u\in V^\perp\cap w^\perp$ (as $\,(\nabla u)(y)\,$ is
skew-adjoint) and $\,\nabla_{\!w}u\ne0\,$ (since $\,w\notin V$). Hence
$\,w,\nabla_{\!w}u\in V^\perp$ are nonzero and orthogonal to each other, which
yields $\,\dim M-\hs\dim V\ge2$. As $\,\,U\cap N(u)=\,\e_y(U'\cap V)\,$ (by
Lemma \a\cm.1(ii)), we thus obtain $\,\dim M-\hs\dim N\ge2$, which (along with
(a)) gives connectedness and denseness of
$\,M'=M\smallsetminus N(u)$.{\hfill\qd}
\proclaim{Corollary \a\cm.3}Let\/ $\,\vp\,$ be a nonconstant \kip\ on a
K\"ahler manifold, cf. {\rm\S\kp}. Then the open subset on which\/
$\,d\vp\ne0\,$ is connected and dense in\/ $\,M$.
\endproclaim
This is obvious from Lemma \a\cm.2, since $\,u=J(\navp)\,$ is a Killing
field.{\hfill\qd}
\proclaim{Lemma \a\cm.4}Let\/ $\,\vp\,$ satisfy\/ \f{\sr.1} on a K\"ahler
manifold\/ $\,(M,g)\,$ of complex dimension\/ $\,m\ge1$, and let\/
$\,(0,L)\ni\sa\mapsto x(\sa)\in M'$ be an integral curve of the vector field\/
$\,\pm\hs v/|v|$, where\/ $\,\pm\,$ is some sign, $\,v=\navp\,$ and\/
$\,M'\subset M\,$ is the open set on which\/ $\,d\vp\ne0$. Then, for the
eigenvalue function\/ $\,\si\,$ appearing in\/ \f{\sr.4}, either\/ $\,\si=0\,$
at\/ $\,x(\sa)\,$ for all\/ $\,\sa\in(0,L)$, or\/ $\,\si\ne0\,$ at all\/
$\,x(\sa)$, $\,\sa\in(0,L)$.
\endproclaim
In fact, let $\,I'$ be a maximal (nonempty) subinterval of $\,(0,L)\,$ with
$\,\si\ne0\,$ at all $\,x(\sa)$, $\,\sa\in I'$. By Lemmas \a\qs.1(b) and
\a\de.1, there exists a constant $\,\y\,$ such that $\,Q=2(\vp-\y)\hs\si\,$ at
every $\,x(\sa)$, $\,\sa\in I'$, where $\,Q=g(v,v)$. Hence, if such a
subinterval exists, it must coincide with $\,(0,L)$, or else it would have an
endpoint $\,L'\in(0,L)$, with $\,Q=\si=0\,$ at $\,x(L')$, contrary to our
assumption that $\,x(\sa)\in M'$, i.e., $\,v\ne0\,$ at $\,x(\sa)$, for all
$\,\sa\in(0,L)$.{\hfill\qd}
\proclaim{Lemma \a\cm.5}Given\/ $\,\vp\,$ with\/ \f{\sr.1} on a K\"ahler
manifold\/ $\,(M,g)\,$ of complex dimension\/ $\,m\ge1$, let\/
$\,(Q,Y,\hs\sc\hs,\si,\ta,\la,\my)\,$ be as in\/ {\rm Definition \a\sr.2}, and
let\/ $\,M'\subset M\,$ be the open set on which\/ $\,d\vp\ne0$. Then either\/
$\,\si=0\,$ identically on\/ $\,M'$, or\/ $\,\si\ne0\,$ everywhere in\/
$\,M'$. In the latter case, there exists a constant\/ $\,\,\y\,$ such that
relations\/ $\,Q/\si=2(\vp-\y)\,$ and\/ $\,\vp\ne\y\,$ hold everywhere in\/
$\,M'$.
\endproclaim
In fact, given $\,x\in M'$ with $\,\si(x)=0$, let $\,P\,$ be the connected
component, containing $\,x$, of the pre\-im\-age of $\,\vp(x)\,$
under $\,\vp:M'\to\bbR\hs$. Thus, $\,\si=0\,$ on $\,P\,$ since $\,\si(x)=0\,$
and, by Lemma \a\qs.1(a), $\,\si\,$ restricted to $\,P\,$ is constant. Lemma
\a\cm.4(b) now implies that $\,\si=0\,$ along every integral curve of
$\,v=\navp\,$ in $\,M'$ which intersects $\,P$. As $\,P\,$ is a submanifold of
$\,M'$, while $\,v=\navp$ is normal to $\,P\,$ and nonzero everywhere in
$\,M'$, the union of these integral curves contains a neighborhood of $\,x\,$
in $\,M'$. Thus, the set of all $\,x\in M'$ with $\,\si(x)=0\,$ is not only
(relatively) closed, but also open in $\,M'$. Our either-or claim now follows
from connectedness of $\,M'$ (Corollary \a\cm.3). Finally, if $\,\si\ne0\,$
everywhere in $\,M'$, Lemmas \a\qs.1(b) and \a\de.1 show that the required
constant $\,\y\,$ exists, locally, in $\,M'$. As $\,M'$ is connected
(Corollary \a\cm.3), $\,\y\,$ has a unique value throughout $\,M'$.{\hfill\qd}
\medskip

\head\S\rr. Special K\"ahler-Ricci potentials and reducibility\endhead
\proclaim{Lemma \a\rr.1}Suppose that\/ $\,(Q,Y,\hs\sc\hs,\si,\ta,\la,\my)\,$
is the septuple associated, as in\/ {\rm Definition \a\sr.2}, with a special
\krp\/ $\,\vp\,$ on a K\"ahler manifold\/ $\,(M,g)\,$ of complex dimension\/
$\,m\ge2$, and\/ $\,v,u,\Cal V,\Cal H\,$ are given by\/ \f{\sr.3}, while\/
$\,w,w'$ are\/ $\,C^1$ vector fields, defined in the open set\/ $\,\,M'$ on
which\/ $\,d\vp\ne0\,$ and orthogonal to\/ $\,v\,$ and\/ $\,u$. Denoting\/
$\,[\nabla_{\!w}w']\vrt$ and\/ $\,[\hs w\hs,w'\hs]\vrt$ the\/ $\,\Cal V\,$
components of\/ $\,\nabla_{\!w}w'$ and the Lie bracket\/ $\,[w,w']\,$ relative
to the decomposition\/ $\,TM'=\Cal H\oplus\Cal V$, we have
$$\alignedat2
&\text{\rm a)}\hskip4.5pt&&
Q\hs[\nabla_{\!w}w']\vrt=-\hs\si\hs[g(w,w')v\hs+\hs\om(w,w')u]\quad
\text{\rm and}\\
&\text{\rm b)}\hskip4.5pt&&
Q\hs[\hs w\hs,w'\hs]\vrt=-\hs2\si\,\om(w,w')\hs u\hs,\hskip5pt
\text{\rm where}\hskip4pt\om\hskip2.7pt\text{\rm is\ the\ K\"ahler\ form\ of}
\hskip3.5pt(M,g),\hskip-4pt\\
&\text{\rm c)}\hskip4.5pt&&
g(\nabla_{\!w}v,w')\,=\,\si\hs g(w,w')\,,\quad
g(\nabla_{\!w}u,w')\,=\,\si\hskip1.3pt\om(w,w')\,.
\endalignedat\ff\rr.1$$
If, in addition, $\,w,w'$ commute with\/ $\,v\,$ and\/ $\,u$, then
$$\alignedat2
&\text{\rm\ptmi i)}\quad&&
\nabla_{\!v}w\,=\,\nabla_{\!w}v\,=\,\si\hs w\,,\quad
\nabla_{\!u}w\,=\,\nabla_{\!w}u\,=\,\si\hs Jw\,,\hskip95pt\\
&\text{\rm ii)}\quad&&
d_v\varPsi\,=\,\,d_u\varPsi\,=\,0\qquad\text{\rm for}\quad
\varPsi\,=\,\si\hs g(w,w')/Q\,.\endalignedat\ff\rr.2$$
\endproclaim
\demo{Proof}By \f{\pr.5} and \f{\sr.4},
$\,g(\nabla_{\!w}v,w')=(\nabla d\vp)(w,w')=\si\hs g(w,w')$, while
\f{\pr.5}, \f{\sr.4} and \f{\kr.1} yield
$\,g(\nabla_{\!w}u,w')=g(\nabla_{\!w}(Jv),w')=g(J(\nabla_{\!w}v),w')
=-\hs g(\nabla_{\!w}v,Jw')$ $=-\hs(\nabla d\vp)(w,Jw')=-\hs\si\hs g(w,Jw')
=\si\hskip1.3pt\om(w,w')$, which proves \f{\rr.1.c}.

Since \f{\rr.2.ii} is obvious when $\,\si=0\,$ identically, let us assume (cf.
Lemma \a\cm.5) that $\,\si\ne0\,$ everywhere in $\,M'$. By \f{\pr.2.ii},
$\,g(\nabla_{\!v}w,w')=g(\nabla_{\!w}v,w')\,$ and
$\,g(\nabla_{\!u}w,w')=g(\nabla_{\!w}u,w')\,$ whenever $\,w,w'$ commute
with $\,v\,$ and $\,u$, so that \f{\rr.1.c} implies
$\,d_v\varUpsilon=2\si\varUpsilon\,$ and $\,d_u\varUpsilon=0\,$ for the
function $\,\varUpsilon=g(w,w')$. Also, by Lemma \a\cm.5,
$\,d(\si/Q)=-\hs(\si/Q)^2\hs d(Q/\si)=-\hs2(\si/Q)^2\hs d\vp$, so that
$\,d_v(\si/Q)=-\hs2\hs\si^2/Q\,$ and $\,d_u(\si/Q)=0\,$ since $\,d_v\vp=Q\,$
and $\,d_u\vp=g(u,v)=0\,$ (see \f{\pr.6.i}, \f{\sr.5}). As
$\,\varPsi=(\si/Q)\varUpsilon\,$ in \f{\rr.2.ii}, our formulae for
$\,d_v\varUpsilon,\hs d_u\varUpsilon\,$ give \f{\rr.2.ii}.

Next, $\,g(v,\nabla_{\!w}w')=-\hs g(\nabla_{\!w}v,w')\,$ since $\,g(v,w')=0\,$
(and similarly for $\,u\,$ instead of $\,v$). This proves \f{\rr.1.a}; in
fact, \f{\rr.1.c} and \f{\sr.5} now show that both sides of \f{\rr.1.a} have
the same inner product with $\,u\,$ and $\,v$, while
$\,\Cal V=\,\text{\rm Span}\,\{v,u\}$. Thus, \f{\rr.1.b} is immediate from
\f{\rr.1.a} and \f{\pr.2.ii}.

Finally, for $\,w,w'$ orthogonal to $\,v,u\,$ and commuting with them,
$\,\nabla_{\!v}w\,$ and $\,\nabla_{\!u}w\,$ are orthogonal to $\,v\,$ and
$\,u\,$ (by the Leibniz rule, cf. the last paragraph,  and Lemma \a\sr.5(ii)),
while $\,\nabla_{\!v}w=\nabla_{\!w}v$, $\,\nabla_{\!u}w=\nabla_{\!w}u\,$ by
\f{\pr.2.ii}. Hence \f{\rr.2.i} follows, as one sees taking the inner products
of the quantities involved with any vector field $\,w'$ orthogonal to $\,v\,$
and $\,u\,$ and using \f{\rr.1.c}. This completes the proof.{\hfill\qd}
\enddemo
\proclaim{Corollary \a\rr.2}Let\/ $\,\vp\,$ satisfy\/ \f{\sr.1} on a K\"ahler
manifold\/ $\,(M,g)\,$ of complex dimension\/ $\,m\ge2$, and let\/
$\,M'\subset M\,$ be the open set of points at which\/ $\,d\vp\ne0$. For the
distributions\/ $\,\Cal V,\Cal H\,$ and the function\/ $\,\si\,$ on\/ $\,M'$
given by\/ \f{\sr.3} -- \f{\sr.4},
\widestnumber\item{(iii)}\roster
\item"(i)"$\si=0\,$ identically on\/ $\,M'$ if and only if\/ $\,\Cal H\,$ is
integrable.
\item"(ii)"If\/ $\,\si=0\,$ identically on\/ $\,M'$, then, locally,
$\,(M'\!,g)\,$ is a Riem\-ann\-i\-an product with\/ $\,\Cal H,\Cal V\,$
serving as the factor distributions.
\item"(iii)"Conversely, if there exists a nonempty open connected set\/
$\,\,U\subset M'$ such that\/ $\,(U,g)\,$ is a non-Einstein,
Riem\-ann\-i\-an-and-K\"ahler product of two lower\diml\ K\"ahler manifolds,
then\/ $\,\si=0\,$ identically on\/ $\,M'$.
\endroster
\endproclaim
In fact, \f{\rr.1.b} implies (i). If $\,\si=0\,$ on $\,M'$, then $\,\Cal V\,$
is parallel since, by Lemma \a\sr.5(ii) and \f{\rr.2.i}, the covariant
derivatives of $\,v=\navp\,$ and $\,u=Jv\,$ in all directions lie in
$\,\Cal V$. This yields (ii). Finally, making $\,\,U\,$ in (iii) smaller, we
may assume \f{\sr.4} with $\,\la\ne\my\,$ everywhere in $\,\,U\,$ (so that
$\,\,\ri\,\,$ is not a multiple of $\,g\,$ at any point of $\,\,U$). Then
$\,\Cal V\,$ restricted to $\,\,U\,$ is an eigenspace distribution of complex
dimension one for $\,\,\ri$, and so $\,\Cal V\,$ must be contained in one
factor distribution of the product decomposition in (iii). Now \f{\rr.1.a}
with $\,w=w'\ne0\,$ tangent to the other factor gives $\,\si=0\,$ on $\,\,U$,
and hence on $\,M'$ (see Lemma \a\cm.5).{\hfill\qd}

\head\S\sh. Submersions and horizontal vectors\endhead
Given manifolds $\,M,N\,$ and a $\,C^\infty$ submersion $\,\proj:M\to N$, we
will refer to the subbundle $\,\Cal V=\,\text{\rm Ker}\,d\proj\,$ of $\,TM\,$
as the {\it vertical distribution\/} of $\,\proj$. A fixed Riemannian metric
$\,g\,$ on $\,M\,$ then gives rise to the {\it horizontal distribution\/}
$\,\Cal H=\Cal V^\perp$, and we denote $\,w\hrz,\hs w\vrt$ the $\,\Cal H\,$
and $\,\Cal V\,$ components of any vector (field) $\,w\,$ tangent to $\,M\,$
relative to the decomposition $\,TM=\Cal H\oplus\Cal V$. Also, the
Levi-Civita connection $\,\nabla\,$ of $\,g\,$ induces connections
$\,\nabla\hrz,\hs\nabla\vrt$ in the vector bundles $\,\Cal H,\hs\Cal V\,$ over
$\,M$, with
$$\nabla\hrz_{\!w}w'\,=\,[\nabla_{\!w}w']\hrz\,,\qquad
\nabla\vrt_{\!w}w'\,=\,[\nabla_{\!w}w']\vrt\ff\sh.1$$
for vector fields $\,w\,$ and horizontal/vertical $\,C^1$ vector fields
$\,w'$ on $\,M$.
\remark{Remark \a\sh.1}Given a Riemannian manifold $\,(M,g)\,$ and a {\it
surjective\/} submersion $\,\proj:M\to N$, we will use the horizontal-lift
operation to identify vector fields on $\,N\,$ with those horizontal vector
fields on $\,M\,$ (i.e., sections of $\,\Cal H$) which are
$\,\proj$-pro\-jecta\-ble. This does not usually lead to notational confusion,
with one notable exception. Namely, the Lie bracket of two
$\,\proj$-pro\-jecta\-ble vector
fields is $\,\proj$-pro\-jecta\-ble onto the Lie bracket of their
$\,\proj$-im\-ages (see \cite{\kno, p. 10}), but Lie brackets of
horizontal fields need not be horizontal. Also, all vertical vector fields
are $\,\proj$-pro\-jecta\-ble (onto $\,0$). Thus, if
$\,[\hs w\hs,w'\hs]\,$ is the Lie bracket in $\,M\,$ of any
$\,\proj$-pro\-jecta\-ble horizontal $\,C^\infty$ vector fields $\,w,w'\,$ on
$\,M$, then, under our identification,
$$[\hs w\hs,w'\hs]\hrz\hskip7pt\text{\rm corresponds\ to\ the\ Lie\ bracket\
of}\hskip6ptw\hskip4pt\text{\rm and}\hskip4ptw'\hskip5pt\text{\rm in}
\hskip5ptN.\ff\sh.2$$
\endremark
We will encounter cases in which $\,(M,g)\,$ is a Riemannian manifold with a
submersion $\,\proj:M\to N\,$ and, in addition, there are given objects
$\,u,\ba,\vta\,$ such that
\vskip4pt
\hbox{\hskip-.8pt
\vbox{\hbox{\f{\sh.3}}\vskip15pt}
\hskip7.5pt
\vbox{
\hbox{a)\hskip7pt$u\,$ is a vertical vector field on $\hs M\hs$ and
$\hs\ba\hs$ is a differential $\hs2$-form on $\hs M$,}
\hbox{b)\hskip7pt$\vta:\Cal H\to\Cal H\,$ is a vector-bundle morphism, with
$\,\Cal H\,$ as above, while}
\hbox{c)\hskip7pt$[w,w']\vrt=\,\ba(w,w')\hs u\,$ for all horizontal $\,C^1$
vector fields $\,w,w'$, and}
\hbox{d)\hskip7pt$[\nabla_{\!u}w]\hrz=\,\vta w\,$ for all
$\,\proj$-pro\-jecta\-ble horizontal $\,C^1$ vector fields $\,w$.}
}}
\vskip5pt
\noindent Conditions \f{\sh.3.c} -- \f{\sh.3.d} are geometrically natural. In
\f{\sh.3.c}, the dependence of $\,[w,w']\vrt$ on $\,w,w'$ is pointwise, i.e.,
the value of $\,[w,w']\vrt$ at any $\,x\in M\,$ depends only on $\,w(x)\,$ and
$\,w'(x)$, as one sees multiplying $\,w\,$ or $\,w'$ by a function (cf.
\cite{\kno, Proposition 3.1 on p. 26}). In \f{\sh.3.d},
$\,\proj$-pro\-jecta\-ble horizontal vector fields form a trivializing space
of sections in the vector bundle $\,\Cal H\,$ restricted to any of the
submanifolds $\,\proj^{-1}(y)$, $\,y\in N$, and this is why, for
$\,\nabla\hrz$ as in \f{\sh.1}, $\,\nabla\hrz_{\!u}$ makes sense as a {\it
bundle morphism}. (For details, see Remark \a\sh.2 below.)
\remark{Remark \a\sh.2}Let $\,V\,$ be a {\it trivializing space of\/}
$\,C^\infty$ {\it sections\/} of a vector bundle $\,\Cal E\,$ over any base
manifold, that is, a vector space of sections such that, for every point
$\,y\,$ of the base and every $\,w_0\in\Cal E_y$, there exists a {\it
unique\/}  $\,w\in V\,$ with $\,w(y)=w_0$. For any fixed connection
$\,\nabla\,$ in $\,\Cal E\,$ and any vector field $\,u\,$ on the base, we then
may treat $\,\nabla_{\!u}$ as a vector-bundle morphism $\,\Cal E\to\Cal E$,
sending any $\,w_0\in\Cal E_y$ to $\,(\nabla_{\!u}w)(y)$, with $\,w\,$ chosen
as above for $\,y\,$ and $\,w_0$.
\endremark

\head\S\hr. Horizontally homothetic submersions\endhead
We say that a $\,C^\infty$ submersion $\,\proj:M\to N\,$ between Riemannian
manifolds $\,(M,g),\hs(N,h)\,$ is {\it horizontally homothetic\/} if,
for some function $\,\fe:M\to(0,\infty)$,
$$\proj^*h\,=\,g/\!\fe\quad\text{\rm on}\quad\Cal H\quad\text{\rm and}\quad
[\nabla\!\fe]\hrz\,=\,0\,,\ff\hr.1$$
with $\,...\hrz,\hs\Cal H\,$ as in \S\sh, i.e., if the pullback tensor
$\,\proj^*h\,$ agrees with $\,g/\!\fe\,$ on the horizontal distribution, and
$\,\fe\,$ has a vertical gradient.

A submersion $\,\proj\,$ with \f{\hr.1} obviously becomes a {\it Riemannian
submersion\/} if one replaces $\,g\,$ with the conformally related metric
$\,\tilde g=g/\!\fe$. Conversely, any Riemannian submersion $\,\proj:M\to N\,$
between $\,(M,\tilde g)\,$ and $\,(N,h)\,$ is a horizontally homothetic
submersion between $\,(M,g)\,$ and $\,(N,h)$, where $\,g=\fe\tilde g\,$ for
any fixed $\,C^\infty$ function $\,\fe:M\to(0,\infty)\,$ with
$\,[\nabla\!\fe]\hrz=0$. Further examples of \f{\hr.1} are provided by
projections $\,\proj:M\to\srf\,$ of a product manifold $\,M=N\times\srf\,$
endowed with a {\it warped product\/} metric $\,g$, obtained from a metric
$\,h\,$ on $\,N\,$ along with some metric on $\,\srf\,$ and a function
$\,\srf\to(0,\infty)$. See, e.g., \cite{\bes, p. 237}.
\proclaim{Lemma \a\hr.1}Let\/ $\,\proj:M\to N\,$ be a horizontally
homothetic, surjective submersion between Riemannian manifolds\/
$\,(M,g),\hs(N,h)$, and let\/ $\,\nabla,\hskip1pt\text{\rm D}\,\,$ stand for
the Levi-Civita connections of\/ $\,g\,$ and\/ $\,h$, respectively. For any\/
$\,\proj$-pro\-jecta\-ble horizontal\/ $\,C^\infty$ vector fields\/ $\,w,w'$
on\/ $\,\,M\,$ treated as vector fields on\/ $\,N$, cf. {\rm Remark \a\sh.1},
\widestnumber\item{(ii)}\roster
\item"(i)"The vector field\/ $\,\,[\nabla_{\!w}w']\hrz$ on\/ $\,\,M\,$ is\/
$\,\proj$-pro\-jecta\-ble.
\item"(ii)"$\,\text{\rm D}_ww'$ equals\/ $\,[\nabla_{\!w}w']\hrz$ treated as
a vector field on\/ $\,N$.
\endroster
\endproclaim
\demo{Proof}Our $\,g$, $\,\nabla\,$ and any $\,\proj$-pro\-jecta\-ble
horizontal $\,C^\infty$ vector fields $\,u,v,w\,$ in $\,M\,$ satisfy
\f{\pr.2.i} and, since $\,\Cal V\,$ and $\,\Cal H\,$ are
$\,g$-or\-tho\-gon\-al, \f{\pr.2.i} will remain valid when $\,\nabla_{\!w}v\,$
and all Lie brackets are replaced by their horizontal components. Since
$\,\fe\,$ in \f{\hr.1} is constant along $\,u,v\,$ and $\,w$, \f{\pr.2.i} will
still hold if we further replace every occurrence of $\,g\,$ by $\,g/\!\fe$.
Finally, using \f{\sh.2} and \f{\hr.1} we see that this last modified version
of \f{\pr.2.i} is also valid with $\,\,\text{\rm D}_wv\,$ instead of
$\,[\nabla_{\!w}v]\hrz$ (since \f{\pr.2.i} holds for $\,(N,h)\,$ as well).
Both $\,\,\text{\rm D}_wv,\hskip3.5pt[\nabla_{\!w}v]\hrz$ thus are
horizontal vector fields on $\,M\,$ having the same inner product
with every $\,\proj$-pro\-jecta\-ble horizontal field, and hence must
coincide. This in turn implies (i), and completes the proof.{\hfill\qd}
\enddemo
\proclaim{Lemma \a\hr.2}Let\/ $\,\,\Cal V,\Cal H\,$ be the vertical and
horizontal distributions of a holomorphic, surjective, horizontally homothetic
submersion\/ $\,\proj:M\to N\,$ between K\"ahler manifolds\/ $\,(M,g)\,$ and\/
$\,(N,h)$, and let\/ \f{\sh.3} hold for some\/ $\,u$, $\,\ba\,$ and a
complex-linear morphism\/ $\,\vta:\Cal H\to\Cal H$. Denoting\/
$\,\rho$, $\,\roh$ and\/ $\,\varOmega\,$ the Ricci forms of\/ $\,g,\hs h\,$
and, respectively, the curvature form, cf. {\rm Remark \a\cc.1}, of the
connection in the highest complex exterior power of\/ $\,\Cal V$, induced by\/
$\,\nabla\vrt$ with\/ \f{\sh.1}, and letting\/ $\,\,\trc\vta:M\to\bbC\,$ be
the pointwise trace of\/ $\,\vta$, we have
$$\rho\,\,=\,\,\varOmega\,\,+\,\,\proj^*\roh\,
-\,\,i\hs[\trc\vta]\hs\ba
\qquad\text{\rm on}\quad\Cal H\,.\ff\hr.2$$
\endproclaim
\demo{Proof}Let $\,m,q\,$ be the complex dimensions of $\,M\,$ and $\,N$, and
let the ranges of indices be: $\,a,b\in\{1,\dots,m\}$,
$\,j,k\in\{1,\dots,q\}$, $\,\lambda,\mu\in\{q+1,\dots,m\}$. Locally, we can
choose $\,\proj$-pro\-jecta\-ble horizontal vector fields $\,w_j$ and vertical
vector fields $\,w_\lambda$, which together trivialize the complex vector
bundle $\,TM\,$ on an open set $\,\,U\subset M$. Choosing the corresponding
$\,1$-forms $\,\varGamma_{\!a}^b$ as in \f{\cc.2}, we now have
$\,[\nabla_{\!v}w_j]\hrz=\varGamma_{\!j}^k(v)\hs w_k$ and
$\,[\nabla_{\!v}w_\lambda]\vrt=\varGamma_{\!\lambda}^\mu(v)\hs w_\mu$ for
any vector field $\,v\,$ on $\,\,U\,$ (due to \f{\cc.2} with
$\,[w_j]\vrt=[w_\lambda]\hrz=0$). Now \f{\sh.1} gives
$\,\nabla\vrt_{\!v}w_\lambda=\varGamma_{\!\lambda}^\mu(v)\hs w_\mu$, and so,
for $\,\hat w=w_{q+1}\wedge\ldots\wedge w_m$ we have
$\,\nabla\vrt_{\!v}\hat w=\varGamma(v)\hs\hat w\,$ with
$\,\varGamma=\varGamma_{\!\lambda}^\lambda$. Thus,
$\,\varOmega\,=\,i\hskip1ptd\varGamma_{\!\lambda}^\lambda$, since
$\,\varOmega=i\hskip1ptd\varGamma\,$ (Remark \a\cc.1). Note that we sum
over repeated indices.

As $\,[\nabla_{\!v}w_j]\hrz=\varGamma_{\!j}^k(v)\hs w_k$, Lemma \a\hr.1(i)
shows that $\,\varGamma_{\!j}^k(v)\,$ is $\,\proj$-pro\-jecta\-ble, i.e.,
descends to a function in $\,N$, whenever $\,v\,$ is a
$\,\proj$-pro\-jecta\-ble horizontal vector field in $\,M$, so that the
$\,1$-forms $\,\varGamma_{\!j}^k$ in $\,M$, restricted to $\,\Cal H$, may be
viewed as $\,1$-forms in $\,N\,$ (cf. Remark \a\sh.1). Lemma \a\hr.1(ii) then
gives $\,\,\text{\rm D}_{v}w_j=\varGamma_{\!j}^k(v)\hs w_k$, with the
$\,\proj$-pro\-jecta\-ble horizontal vector fields $\,v\,$ and $\,w_k$ in
$\,M\,$ treated, simultaneously, as vector fields in $\,N$. Now, by Lemma
\a\kr.2 (for $\,(N,h)$), $\,\roh=i\,d_N\hskip-1pt\varGamma_{\!j}^j$, where
$\,d_N$ is the exterior derivative in $\,N$. Note that $\,d_N$ differs from
$\,d$, the exterior derivative in $\,M$, due to the Lie-bracket term in
\f{\pr.9.ii}, as the $\,M\,$ and $\,N\,$ Lie brackets in \f{\sh.2} differ by
$\,[w,w']\vrt$, i.e., by $\,\ba(w,w')\hs u$, cf. \f{\sh.3.c}. Thus,
$\,d\varGamma_{\!j}^j=\hs d_N\hskip-1pt\varGamma_{\!j}^j
-\varGamma_{\!j}^j(u)\hs\ba\,$ on $\,\Cal H$, while
$\,\varGamma_{\!j}^j(u)=\,\trc\vta$, since
$\,\vta w_j=[\nabla_{\!u}w_j]\hrz=\varGamma_{\!j}^k(u)\hs w_k$. Consequently,
$\,d\varGamma_{\!j}^j=\,d_N\hskip-1pt\varGamma_{\!j}^j
-\,[\trc\vta]\hs\ba$. However, Lemma \a\kr.2 for $\,(M,g)\,$ gives
$\,\rho\,=\,i\hskip1ptd\varGamma_{\!a}^a
=i\hskip1ptd\varGamma_{\!j}^j+\hs i\hskip1ptd\varGamma_{\!\lambda}^\lambda$.
Hence the last relation, along with
$\,\varOmega\,=\,i\hskip1ptd\varGamma_{\!\lambda}^\lambda$ and
$\,\roh=i\,d_N\hskip-1pt\varGamma_{\!j}^j$ (see above), yields \f{\hr.2}. This
completes the proof.{\hfill\qd}
\enddemo
\proclaim{Corollary \a\hr.3}Suppose that a surjective holomorphic submersion\/
$\,\proj:M'\to N\,$ between K\"ahler manifolds\/ $\,(M'\!,g),\hs(N,h)\,$ of
complex dimensions\/ $\,m\ge2\,$ and\/ $\,m-1\,$ satisfies\/ \f{\hr.1} for
some\/ $\,\fe:M'\to(0,\infty)$. Also, let\/ $\,v\,$ be a\/ $\,C^\infty$
vertical vector field without zeros on\/ $\,M'$ such that, for all\/
$\,\proj$-pro\-jecta\-ble horizontal\/ $\,C^\infty$ vector fields\/ $\,w,w'$
on\/ $\,M'$,
$$\nabla_{\!v}v\hs=\hs\ta\hskip.4ptv\hs,\hskip9pt
\nabla_{\!u}v\hs=\hs\ta\hskip.4ptu\hs,\hskip9pt
\nabla_{\!w}v\hs=\hs\si\hs w\hs,\hskip9pt\nabla_{\!u}w\hs=\hs\si\hs Jw\hs,
\hskip9pt\text{\rm with}\hskip6ptu\hs=\hs Jv\hs,\ff\hr.3$$
with some functions\/ $\,\si,\ta:M'\to\bbR$. Finally, let\/ \f{\sh.3.c} hold
for the\/ $\,2$-form\/ $\,\ba=-\hs2\ve\hskip1pt\om/\!\fe$, where\/ $\,\ve\,$
is a real constant, $\,u=Jv$, and\/ $\,\om\,$ is the K\"ahler form of
$\,(M'\!,g)$, with\/ $\,...\vrt,\hs...\hrz,\hs\Cal V,\Cal H\,$ as in\/
{\rm\S\sh}.

Denoting\/ $\,\ri\,\,$ and\/ $\,\,\rih$ the Ricci tensors of\/ $\,g\,$ and\/
$\,h$, we have
\widestnumber\item{(ii)}\roster
\item"(i)"If\/ $\,\,\rih=\kx\hs h\,$ for a function\/
$\,\kx:N\to\bbR\hs$, then $\,\,\ri\,=\la\hs g\,$ on $\,\Cal H$, with\/
$\,\la=(\kx-\ve Y)/\!f\,$ and\/ $\,Y=2\ta+2(m-1)\hs\si$.
\item"(ii)"If\/ $\,\,\ri\,=\la\hs g\,$ on $\,\Cal H\,$ for a function\/
$\,\la:M'\to\bbR\hs$, then $\,\,\rih=\kx\hs h\,$ with $\,\kx:N\to\bbR\,$
given by\/ $\,\kx=\la\fe+\ve Y$, where\/ $\,Y=2\ta+2(m-1)\hs\si$.
\endroster
\endproclaim
\demo{Proof}Our $\,v\,$ is a global trivializing section of the complex line
bundle $\,\Cal V\,$ over $\,M'$, while \f{\sh.1} and \f{\hr.3} give
$\,\nabla\vrt_{\!w}v=\varGamma(w)\hs v\,$ for all vectors $\,w$, where
$\,\varGamma\,$ is the complex-val\-ued $\,1$-form with $\,\varGamma(v)=\ta$,
$\,\varGamma(u)=i\ta$, and $\,\varGamma(w)=0\,$ if $\,w\,$ is horizontal.
(Note that the dependence of $\,\nabla\vrt_{\!w}v\,$ on $\,w\,$ is pointwise.)
Using \f{\pr.9.ii} for $\,\xi=\varGamma\,$ and \f{\sh.3.c} we now get
$\,(d\varGamma)(w,w')=2\ve\hs i\ta\,\om(w,w')/\!\fe\,$ for all horizontal
vectors $\,w,w'$. Also, verticality of $\,v\,$ and $\,J$-invariance of
$\,\Cal H$, along with our assumptions \f{\sh.3.c} and \f{\hr.3} show that
conditions \f{\sh.3} are satisfied by by our $\,\ba\,$ with $\,u=Jv\,$ and
$\,\vta=\si J\,$ (so that $\,\vta\,$ multiplies sections of $\,\Cal H\,$ by
the function $\,i\hs\si$). Lemma \a\hr.2 now gives \f{\hr.2}, which here reads
$\,\rho\,-\,\proj^*\roh=-\hs\ve Y\om/\!\fe\,$ on $\,\Cal H$, for
$\,Y=2\ta+2(m-1)\hs\si$. In fact, $\,\varOmega=i\hskip1ptd\varGamma\,$ (see
Remark \a\cc.1; the `highest complex exterior power' of $\,\Cal V\,$ required
in Lemma \a\hr.2 now is $\,\Cal V\,$ itself), so that
$\,\varOmega=-\hs2\ve\ta\hskip1.2pt\om/\!\fe\,$ on $\,\Cal H\,$ (by the above
formula for $\,(d\varGamma)(w,w')$), while $\,\,\trc\vta=(m-1)\hs i\hs\si$,
i.e., $\,-\hs i\hs[\trc\vta]\hs\ba=-\hs2(m-1)\ve\si\hskip1.2pt\om/\!\fe$.
Hence $\,\rho\,-\,\proj^*\roh=-\hs\ve\hs Y\hs\proj^*\omh\,$ on $\,\Cal H$,
where $\,\omh\,$ is the K\"ahler form of $\,(N,h)$. (This is clear since
$\,\proj^*\omh=\hs\om/\!\fe\,$ on $\,\Cal H$, due to \f{\kr.1} and \f{\hr.1}.)

Assertions (i), (ii) now follow from the two formulae for
$\,\rho\,-\,\proj^*\roh$, as $\,\proj^*\omh=\hs\om/\!\fe\,$ on $\,\Cal H$,
while conditions $\,\,\ri\,=\la\hs g\,$ (on $\,\Cal H$) and
$\,\,\rih=\kx\hs h\,$ are equivalent to
$\,\rho=\la\hskip1.2pt\om\hs\,$ (on $\,\Cal H$) and, respectively,
$\,\roh=\kx\hskip1.2pt\omh$. This completes the proof.{\hfill\qd}
\enddemo

\head\S\dx. Proof of the claim made in \S\xm\endhead
Let the data \f{\xm.1} satisfy the conditions listed in \S\xm. Then $\,g\,$
defined in \S\xm\ is a K\"ahler metric on the complex manifold $\,M'$, and
$\,\vp$, as a function $\,M'\to\bbR\hs$, is a special \krp\ on $\,(M'\!,g)\,$
(see \f{\sr.1}). Also, for the septuple $\,Q,Y,\hs\sc\hs,\si,\ta,\la,\my\,$
associated with $\,\vp\,$ as in Definition \a\sr.2, and the Ricci tensor
$\,\,\rih$ of $\,h$,
\widestnumber\item{(c)}\roster
\item"(a)"$Q=a^2r^2\theta$, and $\,d\vp\ne0\,$ everywhere in $\,M'$.
\item"(b)"$\la=(\kx-\ve Y)/\!f$, i.e., $\,\kx=\la\fe+\ve Y$, with
$\,\kx:N\to\bbR\,$ such that $\,\,\rih=\kx\hs h$. (We identify $\,\kx\,$ with
the composite $\,\kx\circ\proj$, i.e., view it as a function on $\,M'$.)
\item"(c)"Either $\,\si=0\,$ identically and $\,\ve=0$, or $\,\si\ne0\,$
everywhere, $\,\ve=\pm\hs1\,$ and $\,\y\,$ in \f{\sr.1} is the same as in
Lemma \a\cm.5, so that $\,\hs Q/\si=2(\vp-\y)$.
\endroster
This can be verified as follows. We denote $\,v,u\,$ the vector fields
\f{\cc.5} on $\,\Cal L$, with our fixed $\,a\ne0$. Also, throughout this
section, $\,w,w'$ stand for
any two $\,C^\infty$ vector fields in $\,N\,$ and, simultaneously, for their
horizontal lifts to $\,\Cal L$, while the symbol $\,J\,$ is used for the
complex structure tensors of both $\,N\,$ and $\,\Cal L$. As $\,\Cal H\,$ is
$\,J$-invariant and $\,\proj:\Cal L\to N\,$ is holomorphic, $\,Jw\,$ means the
same, whether treated as a vector field in $\,N$, or as a
$\,\proj$-pro\-jecta\-ble horizontal field in $\,\Cal L$.

It follows now that $\,v=\navp$, i.e., $\,v\,$ is the $\,g$-gradient of
$\,\vp$, and
$$\aligned
&\nabla_{\!v}v\,=\,-\hs\nabla_{\!u}u\,=\,\ta\hskip.4ptv\,,\quad
\nabla_{\!v}u\,=\,\nabla_{\!u}v\,=\,\ta\hskip.4ptu\,,\quad
\nabla_{\!v}w\,=\,\nabla_{\!w}v\,=\,\si\hskip.4ptw\,,\\
&\nabla_{\!u}w\,=\,\nabla_{\!w}u\,=\,\si\hskip.4ptJw\,,\quad
\nabla_{\!w}w'\,=\,\,\text{\rm D}_ww'\,-\,\ve\hs
[h(w,w')v\,+\,h(Jw,w')u]\hs,\endaligned\ff\dx.1$$
where $\,\nabla,\,\text{\rm D}\,\,$ are the Levi-Civita connections
of $\,g\,$ and $\,h$, while $\,\ta,\si:M'\to\bbR\,$ are given by
$\,2\hs\theta\ta/a=2\hs\theta+r\hs d\theta/dr$,
$\,\hs2\hs\fe\si/a=r\hs\dfe/dr$. In fact, as $\,d\vp/dr=ar\theta\,$ in
\S\xm, \f{\cc.6} gives $\,d_v\vp=Q\,$ with $\,Q=a^2r^2\theta\,$ and
$\,\hs d_u\vp=\hs d_w\vp=0$.
Formulae $\,g(v,v)=g(u,u)=Q\,$ and $\,g(w,w')=\fe h(w,w')$,
$\,g(v,u)=g(v,w)=g(u,w)=0\,$ (due to
$\,\langle v,v\rangle=\langle u,u\rangle=a^2r^2$ and
$\,\,\text{\rm Re}\hskip1pt\langle v,u\rangle=0\,$ in Remark \a\cc.2) now give
$\,v=\navp\,$ (by showing that $\,v-\navp\,$ is $\,g$-or\-thog\-o\-nal to
$\,v,u\,$ and all $\,w$), and hence assertion (a). As $\,J\,$ is skew-adjoint
and $\,\fe\si=\ve Q\,$ (since $\,\dfe/d\vp=2\ve\,$ and $\,d\vp/dr=ar\theta$),
they also imply that the connection $\,\nabla\,$ {\it defined\/} by \f{\dx.1}
makes $\,g\,$ parallel. Next,
$$[w,w']\vrt\,=\,-\hs2\hs\ve\hskip1pt\omh(w,w')u\,
=\,-\hs2\hs\ve\hs h(Jw,w')u\,,\ff\dx.2$$
$\omh\,$ being the K\"ahler form of $\,h$. (This is clear from \f{\cc.8},
\f{\cc.5}, the assumption about the curvature form in \S\xm, and
\f{\kr.1}.) Also, $\,[v,u]=[v,w]=[u,w]=0$, since $\,v,u\,$ and all $\,w\,$ are
obviously invariant under fibrewise multiplications by complex constants,
which constitute the flows of $\,v\,$ and $\,u\,$ in $\,\Cal L$. Therefore, by
\f{\pr.2.ii}, \f{\sh.2} and \f{\dx.2}, $\,\nabla\,$ defined by \f{\dx.1} is
torsion-free and so, as $\,\nabla g=0$, it must coincide with the Levi-Civita
connection of $\,g$.

Furthermore, the complex structure tensor $\,J\,$ in $\,M'$ is
$\,\nabla$-parallel, i.e., commutes with $\,\nabla_{\!v}$, $\,\nabla_{\!u}$
and all $\,\nabla_{\!w}$, as one sees combining \f{\dx.1} with the relation
$\,u=Jv$, immediate from \f{\cc.5}, and \f{\kr.2} for $\,\bz=h$. Hence
$\,g\,$ is a K\"ahler metric on $\,M'$.

Remark \a\pr.1 and \f{\dx.1} now imply the part of \f{\sr.4} involving
$\,\nabla d\vp$, with our $\,\Cal V,\Cal H\,$ and $\,\ta,\si\,$ as in
\f{\dx.1}. The eigenspaces of $\,\nabla d\vp\,$ at every point thus are
$\,J$-invariant and $\,g$-or\-thog\-o\-nal to one another, so that
$\,\nabla d\vp\,$ is Hermitian, i.e., $\,\vp\,$ is a nonconstant \kip\ on
$\,(M'\!,g)\,$ (Lemma \a\kp.2(iii)).

Since our $\,v,u,\Cal V,\Cal H\,$ obviously satisfy \f{\sr.3}, all that
remains to be verified is the part of \f{\sr.4} involving $\,\,\ri\,\,$ (cf.
Remark \a\sr.3), that is, \f{\sr.2.a} for $\,\bz=\,\ri$, at every point
$\,x\in M'$, as well as relations (b), (c) at the beginning of this section.

First, setting $\,Y=\Delta\vp\,$ we have
$\,Y=\dv v=\,\text{\rm Trace}\hskip1.8pt\nabla v\,$ (see \f{\pr.10.i}), and
so, by \f{\dx.1}, $\,Y=2\ta+2(m-1)\hs\si$. Thus, $\,Y\,$ is a function of
$\,r$, i.e., of $\,\vp\,$ (as $\,d\vp/dr=ar\theta\ne0\,$ in \S\xm) and, by
\f{\kp.4}, $\,\imath_v\ri\,=\my\,d\vp=\my\,\imath_vg$, with $\,\my\,$ given by
$\,2\my=-\hs dY/d\vp$. Consequently, $\,v(x)\,$ is, at every point
$\,x\in M'$, an eigenvector of $\,\,\ri(x)\,$ for the eigenvalue $\,\my(x)$.
Hermitian symmetry of $\,\,\ri\,\hs$ (cf. \f{\kr.2}) now implies that
$\,u(x)=Jv(x)\,$ is also an eigenvector of $\,\,\ri(x)\,$ for the same
eigenvalue $\,\my(x)$.

The hypotheses of Corollary \a\hr.3 now are satisfied by $\,\fe,\si,\ta\,$
defined above. In fact, \f{\hr.3} (or,
\f{\hr.1}) is a part of \f{\dx.1} (or, of the definition of $\,g$),
$\,\nabla\!\fe\,$ being vertical since so is $\,\navp=v\,$ and $\,\fe\,$ is a
(linear) function of $\,\vp$. Also, \f{\sh.3.c} with
$\,\ba=-\hs2\ve\hskip1pt\om/\!\fe\,$ is clear from \f{\dx.2}, since relation
$\,g=\fe\hs\proj^*\nh h\,$ on $\,\Cal H\,$ and \f{\kr.1} give
$\,\proj^*\omh=\hs\om/\!\fe\,$ on $\,\Cal H$. Moreover, the premise of
Corollary \a\hr.3(i) is satisfied (cf. \S\xm), and hence so is its
conclusion, i.e., $\,\,\ri\,=\la\hs g\,$ on $\,\Cal H\,$ for
$\,\la=(\kx-\ve Y)/\!f$. Combined with the last paragraph, this yields
\f{\sr.2.a} for $\,\bz=\,\ri$, as well as (b). Finally, by Remark \a\pr.1,
$\,\si\,$ in \f{\sr.4} coincides with $\,\si\,$ in \f{\dx.1}, so that
(as $\,\dfe/d\vp=2\ve\,$ and $\,d\vp/dr=ar\theta$) $\,\si,\hs\ve\,$ are
related as stated in (c), and the remainder of (c) follows from (a) with
$\,\fe\si=\ve Q\,$ and $\,\fe=2\ve\hs(\vp-\y)$.{\hfill\qd}
\remark{Remark \a\dx.1}The construction in \S\xm\ may be conveniently
rewritten using a
modified version of the data \f{\xm.1}, in which the initial quintuple
$\,\jy,r,\theta,\vp,\fe\,$ is replaced by a triple $\,\iy,\vp,Q\,$ formed by
an open interval $\,\iy$, a variable $\,\vp\in\iy$, and a positive
$\,C^\infty$ function $\,Q\,$ of the variable $\,\vp$. This triple is subject
to just one condition, reflecting positivity of $\,\fe$, and stating that
$\,\ve\hs(\vp-\y)>0\,$ for all $\,\vp\in\iy$, unless $\,\ve=0$, or,
equivalently, either $\,\ve=0$, or $\,\y\notin\iy\,$ and
$\,\ve=\,\text{\rm sgn}\,(\vp-\y)\,$ whenever $\,\vp\in\iy$.

In addition, $\,Q\,$ treated as a function $\,M'\to\bbR\,$ is given by
$\,Q=g(\navp,\navp)$.

In fact, let us set $\,Q=a^2r^2\theta$. As $\,d\vp/dr=ar\theta\ne0\,$ for all
$\,r\in\jy$, the assignment $\,r\mapsto\vp\,$ constitutes a $\,C^\infty$
diffeomorphism of $\,\jy\,$ onto some interval $\,\iy$, allowing us to replace
$\,r\hs$ by the new variable $\,\vp\in\iy\,$ and view $\,Q\,$ as a function of
$\,\vp$. Now $\,Q=g(\navp,\navp)$, since $\,\langle v,v\rangle=a^2r^2$ (Remark
\a\cc.2) and $\,v=\navp\,$ (see above).

Conversely, given $\,\iy,Q,\vp\,$ with $\,\ve\hs(\vp-\y)>0\,$ for all
$\,\vp\in\iy\,$ unless $\,\ve=0$, we may choose a diffeomorphism
$\,\vp\mapsto r\,$ of $\,\iy\,$ onto some interval $\,\jy\subset(0,\infty)\,$
such that $\,d\hs[\log r]/d\vp=a/Q$, and treat $\,Q\,$ (originally, a function
of $\,\vp\in\iy\hs$), as well as $\,\vp\,$ itself, as functions of the new
variable $\,r\in\jy$. We then define $\,\theta,\fe:\jy\to(0,\infty)\,$ by
$\,\theta=Q/(ar)^2$ (so that $\,d\vp/dr=ar\theta$), and $\,\fe=1\,$ (when
$\,\ve=0$) or $\,\fe=2\ve\hs(\vp-\y)\,$ with our given constant $\,\y\,$ (when
$\,\ve=\pm\hs1$).
\endremark

\head\S\lc. Linearity of invariant connections\endhead
\proclaim{Lemma \a\lc.1}Let there be given a holomorphic real vector field\/
$\,\,v\,$ on a complex manifold\/ $\,M\,$ of complex dimension\/ $\,m\ge1$,
cf. {\rm\S\kp}, a point\/ $\,x\in M\,$ with\/ $\,v(x)\ne0$, and a real
constant\/ $\,a\ne0$. Then there exists a biholomorphic identification of a
neighborhood\/ $\,\,U\,$ of\/ $\,x\,$ in\/ $\,M\,$ with an open set in\/
$\,\Cal L\smallsetminus N$, for some holomorphic
line bundle\/ $\,\Cal L\,$ over a complex manifold\/ $\,N$, which makes\/
$\,v\,$ and\/ $\,u=Jv\,$ appear as the vector fields\/ \f{\cc.5}. Note that\/
$\,N\subset\Cal L\,$ according to the convention\/ \f{\cc.4}.
\endproclaim
In fact, using the flow of $\,v\,$ and the holomorphic inverse mapping theorem
we can find holomorphic local coordinates $\,y_1,\dots,y_m$ for $\,M\,$ at
$\,x\,$ for which $\,v\,$ and $\,u\,$ are the coordinate vector fields for the
real coordinate directions of $\,\,\text{\rm Re}\,y_m$ and
$\,\,\text{\rm Im}\,y_m$. The new coordinates $\,y_1,\dots,y_{m-1},z\,$ with
$\,z=\e\hskip2ptay_m$ represent points near $\,x\,$ by pairs $\,(y,z)$, where
$\,y=(y_1,\dots,y_{m-1})$. We have thus biholomorphically identified a
neighborhood of $\,x\,$ in $\,M\,$ with $\,N\times D\,$ for some open sets
$\,N\subset\bbC^{m-1}$ and $\,D\subset\bbC\smallsetminus\{0\}$. Treating
$\,N\times D\,$ as an open subset of the product line bundle
$\,\Cal L=N\times\bbC$, we easily obtain \f{\cc.5} for our $\,v,u\,$ (for
instance, expressing their flows first in terms of $\,(y,y_m)$, and then with
the aid of $\,(y,z)$).{\hfill\qd}
\medskip
For vector spaces $\,V,W\,$ with $\,\dim V=1$, linearity of mappings
$\,V\to W\,$ follows from their homogeneity, which is the underlying principle
of the next (well-known)
\proclaim{Lemma \a\lc.2}Suppose that\/ $\,\Cal L\,$ is the total space, with\/
\f{\cc.4}, of a\/ $\,C^\infty$ complex line bundle over a manifold\/ $\,N$,
and\/ $\,\,M'\subset\Cal L\,$ is an open connected set having a nonempty
connected intersection with every fibre\/ $\,\Cal L_y$, $\,y\in N$. Also,
let\/ $\,v,u\,$ be the vector fields\/ \f{\cc.5} on\/ $\,M'$ with a fixed real
constant\/ $\,a\ne0$. A\/ $\,C^\infty$ distribution\/ $\,\Cal H\,$ on\/ $\,M'$
is the restriction to\/ $\,M'$ of the horizontal distribution of a\/
$\,C^\infty$ linear connection in\/ $\,\Cal L\,$ admitting a parallel
Hermitian fibre metric\/ $\,\langle\,,\rangle\,$ if and only if
\widestnumber\item{(c)}\roster
\item"(a)"$\Cal H\,$ is invariant under the local flows of both\/ $\,v\,$
and\/ $\,u$,
\item"(b)"$TM'=\Cal H\oplus\Cal V\,$ for the vertical distribution\/
$\,\Cal V$, cf. {\rm\S\cc},
\item"(c)"The real span of\/ $\,\Cal H\,$ and\/ $\,u\,$ is an integrable
distribution on\/ $\,M'\smallsetminus N$.
\endroster
\endproclaim
\demo{Proof}Necessity of (a) -- (c) is clear: the flows of $\,v\,$ and $\,u\,$
consist of multiplications by complex scalars in $\,\Cal L\,$ (which gives
(a)), (b) is obvious, and, finally, the norm function $\,r\hs$ of
$\,\langle\,,\rangle\,$ is constant along both $\,\Cal H\,$ and\/ $\,u\,$
(which yields (c)).

Conversely, for $\,\Cal H\,$ with (a) -- (c) and any $\,x\in M'$, let us
use a local trivialization of $\,\Cal L\,$ to identify a neighborhood of
$\,x\,$ in $\,M'$ with $\,N'\times D\,$ for some open sets $\,N'\subset N\,$
and $\,D\subset\bbC\hs$, and let $\,\varGamma\,$ be the complex-val\-ued
$\,1$-form on $\,N'$ assigning, to any vector $\,w\in T_yN$, $\,y\in N'$, the
complex number $\,\varGamma(w)\,$ such that the vector
$\,(w,0)\in T_{(y,z)}(N'\times D)$, for some (or any)
$\,z\in D\smallsetminus\{0\}$, has the $\,\Cal V\,$ component
$\,(0,\varGamma(w)z)\,$ (or, equivalently, its $\,\Cal H\,$ component is
$\,(w,-\varGamma(w)z)$). Note that $\,\varGamma\,$ is well-defined, i.e.,
independent of the choice of $\,z$, since (a) implies invariance of
$\,\Cal H$, in addition to that of $\,\Cal V$, under multiplications by
complex scalars in $\,\Cal L$.

As $\,(0,\zeta)\hrz=0\,$ whenever $\,(0,\zeta)\in T_{(y,z)}(N'\times D)$,
the second formula in \f{\cc.7} now shows that $\,\Cal H\,$ coincides, on
$\,N'\times D$, with the horizontal distribution of the unique linear
connection in $\,\Cal L'$ (the restriction of $\,\Cal L\,$ to $\,N'$) whose
connection $\,1$-form in our fixed local trivialization (cf. Remark \a\cc.1)
is $\,\varGamma$. Since the intersections $\,M'\cap\Cal L_y$ are nonempty and
connected, these ``local'' connections fit together to form a single
``global'' linear connection in $\,\Cal L\,$ with the horizontal distribution
$\,\Cal H$.

Finally, by (c), locally in $\,M'\smallsetminus N\,$ there exist $\,C^\infty$
functions $\,\varphi\,$ with $\,d_v\varphi>0\,$ and
$\,d_u\varphi=d_w\varphi=0\,$ for all $\,w\,$ in $\,\Cal H$. As $\,[v,u]=0$,
(a) implies constancy of $\,d_v\varphi\,$ along $\,u\,$ and $\,\Cal H\,$ as
well. (For instance, $\,d_wd_v\varphi=0\,$ for all
$\,C^\infty$ sections
$\,w\,$ of $\,\Cal H$, since $\,d_vd_w\varphi=0\,$ and $\,\,\Lie_vw=[v,w]\,$
is, by (a), a section of $\,\Cal H$.) Thus, $\,d_v\varphi=\Upsilon(\varphi)\,$
for some positive $\,C^\infty$ function $\,\Upsilon\,$ of a real variable.
Choosing a positive function $\,r\hs$ of the variable $\,\varphi\,$ with
$\,\Upsilon(\varphi)\,d\hs[\hs\log\hs r\hs]/d\varphi\,=\,a$, and treating it
as a function on an open subset of $\,M'\smallsetminus N$, we now get
$\,d_vr=ar\,$ and $\,d_ur=0$, so that, by \f{\cc.6}, $\,r\hs$ is the
restriction to an open subset of $\,\Cal L\,$ of the norm function of some
Hermitian fibre metric $\,\langle\,,\rangle\,$ in $\,\Cal L\,$ (Remark
\a\cc.2). The linear connection we found makes $\,\langle\,,\rangle\,$
parallel (i.e., its parallel transport is norm-preserving), since $\,r\hs$ is
constant along the horizontal distribution $\,\Cal H$. This completes the
proof.{\hfill\qd}
\enddemo
\remark{Remark \a\lc.3}Given a special \krp\ $\,\vp\,$ (see \f{\sr.1}) on a
K\"ahler manifold $\,(M,g)\,$ of complex dimension $\,m\ge1$, a point
$\,x\in M\,$ with $\,(d\vp)(x)\ne0$, and a fixed real constant $\,a\ne0$, let
$\,v,u,\Cal V,\Cal H\,$ be as in \f{\sr.3}. Then there exists a holomorphic
line bundle $\,\Cal L\,$ over a complex manifold $\,N\,$ along with a
biholomorphic identification of a neighborhood of $\,x\,$ in $\,M\,$ with an
open subset of $\,\Cal L\smallsetminus N\,$ (notation of \f{\cc.4}), under
which $\,v,u\,$ and $\,\Cal V\,$ become the vector fields \f{\cc.5} and the
vertical distribution in $\,\Cal L\,$ (\S\cc), while $\,\Cal H\,$ appears as
the horizontal distribution of some $\,C^\infty$ linear connection in
$\,\Cal L\,$ that admits a parallel Hermitian fibre metric
$\,\langle\,,\rangle$.

In fact, $\,v\,$ is holomorphic (by \f{\sr.1} and Lemma \a\kp.2(ii)). Applying
Lemma \a\lc.1 we may choose $\,N,\Cal L\,$ and a neighborhood
$\,\,U\,$ of $\,x\,$ in $\,M$, biholomorphically identified with an open set
in $\,\Cal L\smallsetminus N\,$ so that $\,v,u\,$ are given by \f{\cc.5} (and
$\,\Cal V=\,\text{\rm Span}\,\{v,u\}\,$ is the vertical distribution). Our
claim will now follow from Lemma \a\lc.2, once we have shown that conditions
(a) -- (c) in Lemma \a\lc.2 are satisfied. (The ``nonempty connected
intersections'' property of $\,\,U\,$ will hold once we make $\,\,U\,$ smaller
and replace $\,N\,$ by a suitable open submanifold.)

First, (b) is obvious. To establish (a), let us consider any vector field
$\,w\,$ commuting with
$\,v\,$ and $\,u$. Such fields $\,w\,$ exist, locally, and
realize all vectors at any given point $\,x$, as one verifies by prescribing
$\,w\,$ along a submanifold with the tangent space $\,\Cal H_x$ at $\,x$, and
then spreading it over a neighborhood of $\,x\,$ using the flows of $\,v\,$
and $\,u\,$ (which commute, cf. \f{\kp.1.b}). Setting $\,\vartheta=g(v,w)$,
$\,\vartheta'=g(u,w)$, and using the Leibniz rule along with the equalities
$\,\nabla_{\!v}w=\nabla_{\!w}v$, $\,\nabla_{\!u}w=\nabla_{\!w}u\,$ (see
\f{\pr.2.ii}), the fact that $\,\nabla v\,$ is self-adjoint (as $\,v=\navp$,
cf. end of \S\pr), while $\,\nabla u\,$ is skew-adjoint (as $\,u\,$ is
Killing; see \f{\sr.1} and Remark \a\kp.1) and, finally, Lemma \a\sr.5(ii), we
obtain $\,d_v\vartheta=2\ta\vartheta$, $\,d_v\vartheta'=2\ta\vartheta'$,
$\,d_u\vartheta=d_u\vartheta'=0$. (For instance,
$\,d_v\vartheta=g(\nabla_{\!v}v,w)+g(v,\nabla_{\!v}w)$, while
$\,g(v,\nabla_{\!v}w)=g(v,\nabla_{\!w}v)=g(\nabla_{\!v}v,w)=\ta g(v,w)
=\ta\vartheta$.) Thus, if $\,\vartheta=\vartheta'=0\,$ at one point of a leaf
of $\,\Cal V=\,\text{\rm Span}\,\{v,u\}$, we have $\,\vartheta=\vartheta'=0\,$
everywhere in the leaf (due to uniqueness of solutions for ordinary
differential equations), which proves (a). Finally, integrability of
$\,\,\spanr\hs(\Cal H,u)\,$ (i.e., (c)) is clear as
$\,\,\spanr\hs(\Cal H,u)=v^\perp$ and $\,v=\navp\,$ (by \f{\sr.5}, \f{\sr.3}).
\endremark
\remark{Remark \a\lc.4}Let $\,\Cal L\,$ be a  holomorphic line bundle over a
complex manifold $\,N$. Every Hermitian fibre metric $\,\langle\,,\rangle\,$
in $\,\Cal L\,$ then has a {\it Chern connection\/} (\cite{\bch},
\cite{\grh}), i.e., a unique linear connection in $\,\Cal L\,$ making
$\,\langle\,,\rangle\,$ parallel and having a $\,J$-invariant horizontal
distribution, $\,J\,$ being the complex structure tensor on the total space.

In fact, every real subspace $\,V\,$ of real codimension one in a complex
vector space $\hs\,T\,$ with $\,\hs\dim\hs T< \infty\,$ contains a unique
complex subspace $\,H\,$ of complex codimension one (in $\hs\,T$), as one sees
fixing a Hermitian inner product $\,(\hskip1pt,)\,$ in $\hs\,T\,$ and a
nonzero vector $\,w\,$ which is
$\,\,\text{\rm Re}\hskip1pt(\hskip1pt,)$-or\-thog\-o\-nal to $\,V$, and noting
that all complex subspaces
$\,\,\text{\rm Re}\hskip1pt(\hskip1pt,)$-or\-thog\-o\-nal to $\,w\,$ are also
$\,(\hskip1pt,)$-or\-thog\-o\-nal to $\,w$, so that, for dimensional reasons,
our $\,H\,$ must be the $\,(\hskip1pt,)$-or\-thog\-o\-nal complement of $\,w$.

Applying this, at any $\,x\in\Cal L\smallsetminus N\,$ (cf. \f{\cc.4}), to
$\,T=T_x\Cal L\,$ and $\,V\,$ which is the tangent space at $\,x\,$ of the
submanifold given by $\,r=r(x)$, where $\,r\hs$ is the norm function of
$\,\langle\,,\rangle\,$ (Remark \a\cc.2), we obtain a space
$\,H=\Cal H_x\subset V\subset T_x\Cal L$. Along with $\,\Cal H_y=T_yN\,$ for
$\,y\in N\subset\Cal L$, this defines the horizontal distribution
$\,\Cal L\ni x\mapsto\Cal H_x$ of a linear connection making
$\,\langle\,,\rangle\,$ parallel. In fact, assumptions (a) -- (c) in Lemma
\a\lc.2 are satisfied: (a) due to uniqueness of each $\,H=\Cal H_x$, and (c)
since $\,r\hs$ is constant along $\,\,\spanr\hs(\Cal H,u)$.
\endremark

\head\S\ct. A local classification theorem\endhead
In \S\xm\ we constructed a family of special \krp s $\,\vp\,$ on
K\"ahler manifolds $\,(M,g)\,$ of any complex dimension $\,m\ge1$.
We now show that those examples are, locally, the only ones possible. Although
we establish this here just at points in general position, a similar assertion
holds at all points (see \cite{\dmr}).
\proclaim{Theorem \a\ct.1}Let\/ $\,\vp\,$ be a special \krp, with\/ \f{\sr.1},
on a K\"ahler manifold\/ $\,(M,g)\,$ of complex dimension\/ $\,m\ge1$, and
let\/ $\,M'\subset M\,$ be the open set where\/ $\,d\vp\ne0$. Then\/ $\,M'$ is
dense in\/ $\,M\,$ and every point of\/ $\,M'$ has a neighborhood on which, up
to a biholomorphic isometry, $\,g\,$ and\/ $\,\vp\,$ are obtained as in\/ {\rm
\S\xm}.
\endproclaim
\demo{Proof}Denseness of $\,M'$ is clear from Remark \a\kp.4. Applying Remark
\a\lc.3 to $\,M,g,\vp\,$ with any fixed constant $\,a\ne0\,$ and a given point
$\,x\in M'$, and then replacing $\,M\,$ by a sufficiently small neighborhood
of $\,x$, we may assume that $\,M=M'$ and identify $\,M'$ biholomorphically
with an open connected subset of $\,\Cal L\smallsetminus N$, in such a way
that this identification and
$\,\Cal L,N,v,u,\Cal V,\Cal H,\langle\,,\rangle\,$ have all the properties
listed in Remark \a\lc.3 plus
the ``nonempty connected intersections'' property of Lemma \a\lc.2. (For the
latter, it may be necessary to replace $\,M'\!,N\,$ by smaller open
submanifolds.) Depending on whether $\,\si=0\,$ identically or $\,\si\ne0\,$
everywhere in $\,M'$ (see Lemma \a\cm.5), let us define
$\,\fe:M'\to(0,\infty)\,$ by $\,\fe=1\,$ or, respectively, $\,\fe=Q/|\si|\,$
with $\,Q=g(\navp,\navp)$, and set, in both cases,
$\,\ve=\,\sgn\hs\hs\si\in\{-\hs1,0,1\}$, so that $\,\fe\si=\ve Q$. Thus, in
the latter case, $\,\fe=2|\vp-\y|=2\ve\hs(\vp-\y)$, with $\,\y\,$ as in Lemma
\a\cm.5.

There exists a unique Riemannian metric $\,h\,$ on $\,N\,$ such that
$\,g=\fe\hs\proj^*\nh h\,$ on $\,\Cal H$, where $\,\proj:\Cal L\to N\,$ is the
bundle projection. In fact, as $\,\proj(M')=N$, we may define $\,h\,$ by
$\,h(w,w')\,=\,g(w,w')/\!\fe$, with $\,w,w'$ standing for any two $\,C^\infty$
vector fields on $\,N$, as well as for their horizontal lifts to
$\,M'\subset\Cal L$. Namely, Corollary \a\rr.2 (when $\,\si=0$) or
\f{\rr.2.ii} (when $\,\si\ne0$) shows that $\,g(w,w')/\!\fe:M'\to\bbR\,$ is
constant in the vertical directions, and hence it
is well-defined as (i.e., descends to) a function $\,N\to\bbR\hs$. (Since the
flows of $\,v,u\,$ in $\,\Cal L\,$ consist of fibrewise multiplications by
complex constants, the assumptions about $\,w,w'$ in \f{\rr.2.ii} mean that
they are the horizontal lifts of some vector fields in $\,N$.) Furthermore,
denoting $\,r\hs$ the norm function of $\,\langle\,,\rangle$, we have
$\,\langle v,v\rangle=\langle u,u\rangle=a^2r^2$ and
$\,\,\text{\rm Re}\hskip1pt\langle v,u\rangle=0\,$ (see Remark \a\cc.2), and
so $\,g=\theta\,\text{\rm Re}\hskip1pt\langle\,,\rangle\,$ on $\,\Cal V$,
with $\,\theta=Q/(ar)^2$. Since $\,g(\Cal H,\Cal V)=\{0\}\,$ (i.e.,
$\,\Cal H=\Cal V^\perp$, cf. \f{\sr.3}), our $\,g\,$ thus is given by the
same formulae as in \S\xm, and all we need to verify is that the
data \f{\xm.1} just obtained (with $\,\jy\,$ standing for the range of
$\,r:M'\to\bbR$) have all the properties required in \S\xm.

First, making $\,M'\!,N\,$ even smaller if necessary, we may assume that
$\,\vp\,$ is a function of $\,r:M'\to\bbR\,$ (as $\,v=\navp\,$ and
$\,\nabla r\,$ are both $\,g$-or\-thog\-o\-nal to $\,u\,$ and $\,\Cal H$, cf.
\f{\sr.5}, \f{\sr.3}, \f{\cc.6}); hence so are $\,Q\,$ (see Lemma \a\qs.1(a))
and $\,\theta,\fe$. Also, $\,d\vp/dr=\hs d_v\vp/d_vr=Q/(ar)=ar\theta\,$ by
\f{\pr.6.i} and \f{\cc.6}.

Next, $\,\proj:M'\to N\,$ is holomorphic and $\,\Cal H\,$ is $\,J$-invariant
(since so is $\,\Cal V$, cf. \f{\sr.3}). Thus, if we use the same symbol
$\,J\,$ for the complex structure tensors of both $\,M'$ and $\,N$, no
ambiguity will arise as to the meaning of $\,Jw$, where vectors $\,w\,$
tangent to $\,N\,$ are identified with their horizontal lifts to
$\,M'\subset\Cal L$. The metric $\,h\,$ on $\,N\,$ now is clearly Hermitian,
since so is $\,g$, on $\,M'$. Applying Lemma \a\hr.1 to $\,\proj:M'\to N$, we
see that, for $\,w\,$ as above, $\,\,\text{\rm D}_w$ commutes with $\,J\,$ in
$\,N$, as $\,\nabla_{\!w}$ and $\,...\hrz$ both commute with $\,J\,$ in
$\,M'$ (due to the fact that $\,g\,$ is K\"ahler, and $\,J$-invariance of
$\,\Cal H$). Consequently, $\,J\,$ is $\,\,\text{\rm D}$-parallel in $\,N$,
i.e., $\,h\,$ is a K\"ahler metric on the complex manifold $\,N$. Moreover,
$\,\proj:M'\to N\,$ and our $\,\fe,\ve,\hs v=\navp$, along with $\,\si,\ta\,$
as in \f{\sr.4}, satisfy the hypotheses of Corollary \a\hr.3 (provided that
$\,m\ge2$): \f{\hr.1} is our definition of $\,h\,$ (with
$\,[\nabla\!\fe]\hrz=0\,$ as $\,v=\navp\,$ is vertical, cf. \f{\sr.3}, and
$\,\fe\,$ is a linear function of $\,\vp$), while \f{\hr.3} and \f{\sh.3.c}
(for $\,\ba=-\hs2\ve\hskip1pt\om/\!\fe$, with $\,\fe,\ve\,$ as above) follow
from Lemma \a\sr.5(ii), \f{\rr.2.i} and \f{\rr.1.b}. Relation
$\,\,\ri\,=\la\hs g\,$ on $\,\Cal H\,$ in \f{\sr.4} now yields the premise of
Corollary \a\hr.3(ii), so that, by Corollary \a\hr.3, $\,(N,h)\,$ is an
Einstein manifold unless $\,m=2$. This is also trivially true when $\,m=1$,
i.e., $\,N\,$ consists of a single point.

Finally, as $\,g=\fe\hs\proj^*\nh h\,$ on $\,\Cal H$, \f{\kr.1} gives
$\,\hs\om=\fe\hs\proj^*\omh\,$ on $\,\Cal H\,$ for the K\"ahler forms
$\,\omh$ and $\,\hs\om\,$ of $\,h\,$ and $\,g$. Since $\,\fe\si=\ve Q\,$ (see
above), combining \f{\cc.8}, \f{\cc.5} and \f{\rr.1.b} we conclude that our
connection in $\,\Cal L\,$ has the curvature form
$\,\varOmega=-\hs2\ve a\,\omh$. This completes the proof.{\hfill\qd}
\enddemo

\head\S\sz. Solutions to \ \f{\de.1} -- \f{\de.3} \ with $\,\si=0$\endhead
Given an integer $\,m\ge2$, a septuple $\,(Q,Y,\hs\sc\hs,\si,\ta,\la,\my)\,$
of $\,C^\infty$ functions of a variable $\,\vp\,$ satisfies conditions
\f{\de.1} -- \f{\de.3} on some interval on which $\,\si=0\,$ identically if
and only if, for some constants $\,K,\ah,\ts\,$ with
$\,|K|+|\ah|+|\ts|>0$,
$$Q\,=\,-\hs K\vpsq\,+\,(2m-1)^{-1}\hs[\hs\ah\hs\vp^{2m-1}\,
-\,\ts/m\hs]\,,\ff\sz.1$$
while $\,Y=-\hs2K\nh\vp+\ah\hs\vp^{2m-2}$,
$\,\,\sc\,=-\hs(2m-1)(2m-4)\hskip.6ptK-2(m-1)\ah\hs\vp^{2m-3}$,
$\,\ta=-\hs K\nh\vp+\ah\hs\vp^{2m-2}/2$, $\,\la=(3-2m)K$,
$\,\my=K-(m-1)\ah\hs\vp^{2m-3}\,$ and, of course, $\,\si=0$. Note that
positivity of $\,|K|+|\ah|+|\ts|\,$ amounts to \f{\de.3.i} on {\it some\/}
interval.

In fact, a septuple just defined is easily seen to satisfy \f{\de.1} --
\f{\de.3} on a suitable interval. Conversely, let us assume \f{\de.1} --
\f{\de.3} with $\,\si=0\,$ on an interval $\,\iy$. By \f{\de.1}, $\,\la\,$
then is constant, i.e., $\,\la=(3-2m)K\,$ for some $\,K\in\bbR\hs$. If
$\,\ta=\si=0\,$ everywhere in $\,\iy$, \f{\de.3.iii} gives $\,\la=\my$, while
$\,Y=0\,$ and $\,\my=0\,$ identically by \f{\de.3.ii} and
\f{\de.1}, i.e. (cf. \f{\de.1}, \f{\de.3.ii}), $\,Q\,$ is constant and
$\,Y=\hs\sc\hs=\si=\ta=\la=\my=0\,$ on $\,\iy$. (This is a special case of the
above formulae, with $\,K=\ah=0$, $\,\ts\ne0$.) On the other hand, if
$\,\ta\ne\si\,$ on a dense subset of $\,\iy$, \f{\de.3.ii} gives
$\,2(2m-3)\my=(2m-3)[\hs\sc\,+(2-2m)\la]\,$ and, by \f{\de.2.ii},
\f{\de.3.iii}, $\,\vp\hs d\hs\sc/d\vp=(2m-3)\hs\sc\,-(2m-1)(2m-4)\la$.
Replacing $\,\la\,$ in this equality with $\,(3-2m)K\,$ and multiplying by
$\,\vp^{2-2m}$, we obtain the required formula for $\,\,\sc\hs$, with
$\,\ah\,$ representing a constant of integration. As $\,\ta-\si=\ta=Y/2\,$
(due to \f{\de.3.ii} with $\,\si=0$) and $\,2(\la-\my)=\,2m\la-\,\sc\,\,$ (by
\f{\de.3.ii}), with $\,\la=(3-2m)K$, we may rewrite \f{\de.3.iii} as
$\,(m-1)Y=[m(3-2m)K-\,\sc/2]\hs\vp\,$ which, combined with the formula for
$\,\,\sc\hs$, gives the expressions for $\,Y\,$ and $\,\ta$. The first two
equations in \f{\de.2} now yield the formula for $\,\my\,$ and \f{\sz.1} for
some $\,\ts\in\bbR\hs$.

Finally, let $\,\ta\ne\si\,$ somewhere in $\,\iy$, and let $\,\iyp$ be a
maximal subinterval such that $\,\ta\ne\si\,$ on a dense subset of $\,\iyp$.
If one of the endpoints of $\,\iyp$ were an interior point of $\,\iy$, our
discussion of the case $\,\ta=\si\,$ would imply vanishing, at that point, of
$\,\la\,$ and all derivatives of $\,\my$, so that on $\,\iyp$ we would have
the above formulae with $\,K=\ah=0$, contradicting the assumption that
$\,\ta\ne\si=0$.

\head\S\as. A single equation, equivalent to \ \f{\de.1} -- \f{\de.4}\endhead
As in Lemma \a\qs.1, we will refer to individual equations in the
multi-formula expressions \f{\de.1} and \f{\as.1} by labelling them
with lower-case Roman numerals: \f{\de.1.i} -- \f{\de.1.v}, \f{\as.1.i} --
\f{\as.1.v}.

For a fixed integer $\,m\ge1\,$ and any solution
$\,(Q,Y,\hs\sc\hs,\si,\ta,\la,\my)\,$ to \f{\de.1} -- \f{\de.4},
$$\aligned
Q\,&=\,2(\vp-\y)\hs\si\,,\quad\ta\,=\,\si\,+\,(\vp-\y)\hs\si'\,,\quad
Y\,=\,2m\hs\si\,+\,2(\vp-\y)\hs\si'\,,\\
\my\,&=\,-\hs(m+1)\hs\si'\,-\,(\vp-\y)\hs\si''\,,\qquad
\la\,-\,\my\,=\,2(m-1)(\vp-\y)\hskip1pt\si'\ns/\vp\,,\endaligned\ff\as.1$$
where $\,\y\,$ is the constant in Lemma \a\de.1 and $\,\si'=\,d\si/d\vp$. In
fact, \f{\as.1.i} is obvious; differentiating \f{\as.1.i} and using
\f{\de.1.i}, we obtain \f{\as.1.ii}; \f{\as.1.iii} is clear from \f{\de.3.ii}
and \f{\as.1.ii}; differentiating \f{\as.1.iii} we get \f{\as.1.iv} (cf.
\f{\de.1.ii}); while \f{\de.3.iii} gives \f{\as.1.v}, as
$\,\ta-\si=(\vp-\y)\hs\si'$ by \f{\as.1.ii}.
\remark{Remark \a\as.1}A solution $\,(Q,Y,\hs\sc\hs,\si,\ta,\la,\my)\,$
to \f{\de.1} -- \f{\de.4} with a fixed integer $\,m\ge1\,$ is uniquely
determined by its constituent function $\,\si\,$ and the constant $\,\y\,$
defined in Lemma \a\de.1. Explicitly, formulae for $\,Q,Y,\ta,\la,\my\,$ (or,
$\,\,\sc$) in terms of $\,\si\,$ and $\,\y\,$ are provided by \f{\as.1} (or,
by \f{\de.3.ii} and \f{\as.1.iv} -- \f{\as.1.v}).
\endremark
\proclaim{Lemma \a\as.2}Given an integer\/ $\,m\ge2\,$ and a\/ $\,C^1$
solution\/ $\,(Q,Y,\hs\sc\hs,\si,\ta,\la,\my)\,$ to\/ \f{\de.1} -- \f{\de.4}
on an interval, let\/ $\,\y\,$ be the constant in\/ {\rm Lemma \a\de.1}.
Then
$$\vpsq(\vp-\y)\,{\dcu\si\over d\vp^3}\,
=\,[(m-4)\vpsq-2(m-1)\hs\y\hs\vp]\,{\dsq\si\over d\vpsq}\,
+\,2(m-1)(\vp+\y)\,{d\si\over d\vp}\,.\ff\as.2$$
\endproclaim
\demo{Proof}We have \f{\as.2} on every subinterval where
$\,\si'=\hs d\si/d\vp\,$ vanishes identically, and so we may restrict our
consideration to a fixed subinterval in which $\,\si'\ne0\,$ everywhere and
$\,0\ne\vp\ne\y\hs$. (The union of all subintervals of one or the
other type is dense in the original interval of the variable $\,\vp$.) By
\f{\as.1.ii}, $\,\ta-\si=(\vp-\y)\hs\si'$, and so, from \f{\as.1.i},
$\,(\ta-\si)\hs\si=Q\si'\ns/2$. Thus, if we replace $\,2(m-1)(\ta-\si)\,$ on
the left-hand side of \f{\de.2.i} with $\,(\la-\my)\vp\,$ (cf. \f{\de.3.iii})
and then divide by $\,(\la-\my)Q$, we obtain
$\,\vp\my'\,=\,\la-\my\,+\,2(m-1)\hs\my\,+\,2(m-1)^2\si'$. (Note that
$\,\la\ne\my\,$ everywhere, due to our choice of the subinterval and
\f{\as.1.v}.) Replacing $\,\la-\my\,$ and then $\,2(m-1)\hs\my\,$ with the
expressions provided by \f{\as.1.iv} -- \f{\as.1.v}, we now have
$\,\vp\my'=-\hs4(m-1)\hs\si'-2(m-1)(\vp-\y)\hs\si''+2(m-1)(\vp-\y)\hs\si'\ns/\vp$. On
the other hand, differentiating \f{\as.1.iv} we get
$\,\vp\my'=-\,(m+2)\vp\si''\,-\,(\vp-\y)\vp\si'''$. Equating the last two
formulae for $\,\vp\my'$ we arrive at \f{\as.2}, which completes the
proof.{\hfill\qd}
\enddemo
\proclaim{Proposition \a\as.3}Given an integer\/ $\,m\ge2$, a real
number\/ $\,\y\,$ and an interval of the variable\/ $\,\vp$, not containing\/
$\,0\,$ or\/ $\,\y\hs$, there exists a bijective correspondence between\/
$\,C^1$ solutions\/ $\,(Q,Y,\hs\sc\hs,\si,\ta,\la,\my)\,$ to\/ \f{\de.1} --
\f{\de.4} defined on our interval, for which\/ $\,\y\,$ is the constant in\/
{\rm Lemma \a\de.1}, and\/ $\,C^3$ functions\/ $\,\si\,$ on this interval
which are nonzero everywhere and satisfy equation\/ \f{\as.2}. The
correspondence in question associates with\/
$\,(Q,Y,\hs\sc\hs,\si,\ta,\la,\my)\,$ the constituent function\/ $\,\si$.
\endproclaim
\demo{Proof}Any solution $\,(Q,Y,\hs\sc\hs,\si,\ta,\la,\my)\,$ to \f{\de.1} --
\f{\de.4} has $\,\si\ne0\,$ everywhere, cf. \f{\de.4}, and the assignment
$\,(Q,Y,\hs\sc\hs,\si,\ta,\la,\my)\mapsto\si\,$ is injective (Remark \a\as.1).
It is also surjective since, whenever \f{\as.2} and \f{\de.4} hold, the
functions $\,Q,Y,\si,\ta,\la,\my\,$ and $\,\,\sc\hs$, determined by $\,\si\,$
as in Remark \a\as.1, must satisfy \f{\de.1} -- \f{\de.3}. In fact,
\f{\de.1.i} follows
from \f{\as.1.i} and \f{\as.1.ii}; \f{\de.1.ii} from \f{\as.1.iii} and
\f{\as.1.iv}; \f{\de.1.iv} from \f{\as.1.ii}, \f{\as.1.iv}, \f{\as.1.i} and,
again, \f{\as.1.ii}; while \f{\de.1.iii} is clear since \f{\as.1.i} and
\f{\as.1.ii} give $\,Q\si'=2(\vp-\y)\hs\si\si'=2(\ta-\si)\hs\si$. Also, writing both
$\,\la'-\my'$ and $\,\my'$ in terms of $\,\si\,$ and its derivatives (via
\f{\as.1.iv}, \f{\as.1.v}), then replacing $\,\si'''$ with the expression
provided by \f{\as.2}, and noting that the occurrences of $\,\si''$ in the
resulting formula for $\,\la'=(\la'-\my')+\my'$ cancel one another, we easily
obtain \f{\de.1.v} from \f{\as.1.i} and \f{\as.1.v}. Next, \f{\de.3.i}
follows from our assumptions $\,\si\ne0\,$ and $\,\vp\ne\y\hs$, combined with
\f{\as.1.i}; the first formula in \f{\de.3.ii} is immediate from \f{\as.1.iii}
and \f{\as.1.ii}, while the second one is the equality defining $\,\,\sc\,\,$
in Remark \a\as.1; and \f{\de.3.iii} is clear from \f{\as.1.ii} and
\f{\as.1.v}. As \f{\de.1} and \f{\de.3} now imply \f{\de.2} (Remark \a\de.2),
this completes the proof.{\hfill\qd}
\enddemo
\remark{Remark \a\as.4}A $\,C^3$ function $\,\si\,$ of the variable $\,\vp$,
defined on an interval $\,\iy\,$ such that $\,0\notin\iy$, satisfies \f{\as.2}
(with a constant $\,\y\,$ and a fixed integer $\,m\ge2$) if and only if the
function
$$m\hs\si\,-\,(\vp-\y)^2\si''\,-\,\hs m(\vp-\y)\hs\si'\,
+\,\hs2(m-1)(\vp-\y)^2\hs\si'\ns/\vp\,,\ff\as.3$$
with $\,\si'=\hs d\si/d\vp$, is constant on $\,\iy$. (To see this, just
differentiate \f{\as.3}.)
\endremark
\remark{Remark \a\as.5}Proposition \a\as.3 and Remark \a\as.1 lead to an
explicit description of those solutions $\,(Q,Y,\hs\sc\hs,\si,\ta,\la,\my)\,$
to \f{\de.1} -- \f{\de.4} for which $\,\y\,$ defined in Lemma \a\de.1 is zero.
Namely, for a fixed integer $\,m\ge2$, \f{\as.2} with $\,\y=0\,$ is an Euler
equation, and its solution space has a basis formed by the power functions
$\,1$, $\,\vp^m$ and $\,1/\vp$. As $\,Q=2\si\vp\,$ (cf. \f{\as.1}), there
exist constants $\,\nk,\ah,\ts\,$ with
$$\alignedat2
&\text{\rm\ptmi i)}\quad&&
2m\hs\si\,=\,\nk\,+\,m\ah\hs\vp^m\,-\,2(m+1)^{-1}\ts/\vp\,,\hskip83pt\\
&\text{\rm ii)}\quad&&
Q\,=\,m^{-1}\nk\hs\vp\,+\,\ah\hs\vp^{m+1}\,-\,2(m+1)^{-1}\ts/m\,.
\endalignedat\ff\as.4$$
\endremark

\head\S\su. Solutions with $\,\y\ne0$\endhead
\proclaim{Lemma \a\su.1}For any fixed integer\/ $\,m\ge1$, the rational
functions
$$F(\ps)\,=\,{(\ps-2)\ps^{2m-1}\over(\ps-1)^m}\,,\qquad E(\ps)\,
=\,(\ps-1)\,\sum_{k=1}^m{k\over m}{2m-k-1\choose m-1}\ps^{k-1}\ff\su.1$$
of the real variable $\,\ps\hs$ and their derivatives\/ $\,\dot F,\dot E\,$
satisfy the equations
$$\alignedat2
&\text{\rm a)}\quad&&
\ps(\ps-2)\hs\dot F(\ps)\,=\,\vx(\ps)\hs F(\ps)\,,\\
&\text{\rm b)}\quad&&
\ps(\ps-2)\dot E(\ps)\,=\,\vx(\ps)\hs E(\ps)\,+\,2(2m-1)\hs E(0)\,,\hskip9pt
\text{\rm where}\hskip39pt\\
&\text{\rm c)}\quad&&
\vx(\ps)\,=\,m(\ps-1)\,+\,m(\ps-1)^{-1}\hs-\,2(m-1)\,.
\endalignedat\ff\su.2$$
\endproclaim
\demo{Proof}Relation \f{\su.2.a} is immediate from \f{\su.1}. Using a
subscript to mark the dependence of $\,F,E\,$ on $\,m$, we have
$\,(\ps-1)\hs[E_m(\ps)-E_m(0)]=\ps^2\hs E_{m-1}(\ps)\,$ for $\,m\ge2\,$ (as
one sees expanding the difference into powers of $\,\ps\,$ via \f{\su.1}),
while \f{\su.1} gives $\,(\ps-1)\hs F_m(\ps)=\ps^2\hs F_{m-1}(\ps)\,$ and
$\,m\hs E_m(0)=2(2m-3)\hs E_{m-1}(0)$. Also,
$$\alignedat2
&\text{\rm\ptmi i)}\quad&&
[E_m-E_m(0)]/F_m\,=\,E_{m-1}/F_{m-1}\qquad\text{\rm for}\quad m\ge2\,,\\
&\text{\rm ii)}\quad&&
\ps(\ps-2)\,d\hs[E(\ps)/F(\ps)]/d\ps\,=\,2(2m-1)\hs E(0)/F(\ps)\qquad
\text{\rm for}\quad m\ge1\,.\hskip7pt\endalignedat\ff\su.3$$
In fact, \f{\su.3.i} is obvious from the preceding equalities, while
\f{\su.3.ii}, easily verified for $\,m=1$, can be established by induction:
let \f{\su.3.ii} hold when some given $\,m\ge2\,$ is replaced by $\,m-1$. The
above recursion formulae for $\,F(\ps)\,$ and $\,E(0)\,$ now yield
$\,\ps(\ps-2)\,d\hs[\hs E_{m-1}/F_{m-1}\hs]/d\ps
=m\ps^2E_m(0)/[(\ps-1)F_m]$. Using this to evaluate
$\,\ps(\ps-2)\,d\hs[\hs E_m/F_m\hs]/d\ps\,$ from \f{\su.3.i} (where
$\,d\hs[\hs E_m(0)/F_m\hs]/d\ps\,$ is computed with the aid of \f{\su.2.a}),
and noting that $\,m\ps^2(\ps-1)^{-1}-\vx(\ps)=2(2m-1)\,$ for $\,\vx\,$
as in \f{\su.2.c}, we obtain \f{\su.3.ii} for the given $\,m$, as required.
Now, as $\,E=(E/F)F$, \f{\su.2.b} follows from \f{\su.3.ii} and \f{\su.2.a},
which completes the proof.{\hfill\qd}
\enddemo
Note that applying $\,d/d\ps\,$ to \f{\su.2.a} and \f{\su.2.b} we obtain
$$\ps(\ps-2)\hs\ddot\si\,+\,2(\ps-1)\hs\dot\si\,
=\,\dot\vx\hs\si\,+\,\vx\hs\dot\si\,,\hskip8pt\text{\rm both\ for}
\hskip7pt\si=E\hskip6pt\text{\rm and}\hskip6pt\si=F\hs,\ff\su.4$$
where $\,\dot\si=\hs d\si/d\ps\,$ and $\,\si,\hs\dot\si,\vx,\dot\vx\,$
stand for
$\,\si(\ps),\hs\dot\si(\ps),\vx(\ps),\dot\vx(\ps)$.
\proclaim{Lemma \a\su.2}Given an integer\/ $\,m\ge1$, a constant\/
$\,\y\in\bbR\smallsetminus\{0\}$, and an interval\/ $\,\iy\,$ of the
variable\/ $\,\vp$, not containing\/ $\,0\,$ or\/ $\,\y\hs$, the vector
space\/ $\,V\,$ of all\/ $\,C^3$ solutions\/ $\,\si\,$ to the linear
equation\/ \f{\as.2} on\/ $\,\iy\,$ has a basis consisting of the constant
function\/ $\,1\,$ along with\/ $\,E(\ps)\,$ and\/ $\,F(\ps)$, defined by\/
\f{\su.1}, and treated as functions of the variable\/ $\,\vp\,$ via the
substitution\/ $\,\ps=\vp/\y\hs$.
\endproclaim
In fact, $\,d/d\ps\,$ applied to \f{\su.4} yields a third-order linear
differential equation imposed on $\,\si$. Using \f{\su.4} again to express
$\,\si\,$ through $\,\dot\si\,$ and $\,\ddot\si$, and
substituting that expression for $\,\si\,$ (but not its derivatives) in our
third-order equation, we obtain \f{\as.2} with $\,\vp\,$ and $\,\y\,$
replaced everywhere by $\,\ps\,$ and $\,1$. (Note that
$\,\dot\vx=m-m/(\ps-1)^2=m\ps(\ps-2)/(\ps-1)^2$,
$\,\ddot\vx=2m/(\ps-1)^3$, and
$\,\ddot\vx/\dot\vx=2/[\ps(\ps-1)(\ps-2)]$.) However, this modified
version of \f{\as.2}, satisfied by both $\,\si=E\,$ and $\,\si=F$, is
equivalent to \f{\as.2} under the variable change $\,\vp=\y\hs\ps$.{\hfill\qd}
\remark{Remark \a\su.3}Given an integer $\,m\ge2$, let
$\,(Q,Y,\hs\sc\hs,\si,\ta,\la,\my)\,$ be a $\,C^1$ solution to \f{\de.1} --
\f{\de.4} on some given interval, for which $\,\y\,$ in Lemma \a\de.1 is
nonzero. By Lemmas \a\as.2 and \a\su.2, the constituent function $\,\si\,$ is
a linear combination of $\,E(\ps)$, $\,F(\ps)\,$ and $\,1$, treated as
functions of $\,\vp=\y\hs\ps$. Since $\,Q\,$ is related to $\,\si\,$ as in
\f{\as.1}, there exist constants $\,\ax,\bx,\cx\,$ such that
$\,|\ax|+|\bx|+|\cx|>0\,$ (cf. \f{\de.3.i}) and
$$Q\hs=\hs(\ps-1)\hs[\hs\ax\hs+\hs\bx E(\ps)\hs+\hs\cx F(\ps)\hs]\hskip6pt
\text{\rm with}\hskip5pt\ps=\vp/\y\hskip5pt\text{\rm and}\hskip5ptF,E\hskip5pt
\text{\rm as\ in\ \f{\su.1}.}\ff\su.5$$
(About notations, see Remark \a\su.4 below.) Conversely, let $\,Q\,$ be a
function of the variable $\,\vp\,$ given by \f{\su.5} with constants
$\,\ax,\bx,\cx,\y\,$ such that $\,\y\ne0\,$ and $\,|\ax|+|\bx|+|\cx|>0$. Then
$\,Q$, restricted to any interval $\,\iy\,$ on which $\,(\vp-\y)\vp Q\ne0$,
arises in the manner just described from a unique $\,C^1$ solution to
\f{\de.1} -- \f{\de.4} for which $\,\y\,$ is the constant in Lemma \a\de.1;
explicitly, $\,Q,Y,\hs\sc\hs,\si,\ta,\la,\my\,$ are obtained as in Remark
\a\as.1. In fact, by Lemma \a\su.2, $\,\si=Q/[2(\vp-\y)]\,$ then satisfies
\f{\as.2}, and so such a solution exists and is unique in view of Proposition
\a\as.3.

As $\,\si=Q/[2(\vp-\y)]\,$ is a solution to \f{\as.2} on $\,\iy\,$ and
$\,0\notin\iy$, the expression \f{\as.3} is constant on $\,\iy\,$ and has a
simple expression in terms of $\,\ax,\bx,\cx,\y\,$ in \f{\su.5}. Specifically,
\f{\as.3} equals $\,m\ax/(2\y)$. In fact, according
to Remark \a\as.4, by assigning \f{\as.3} to $\,\si\,$ we define a
real-val\-ued linear functional on the three\diml\ space $\,V\,$ of solutions
mentioned in Lemma \a\su.2. Similarly, $\,m\ax/(2\y)\,$ is such a functional,
assigning to $\,\si\,$ the number $\,m/(2\y)\,$ times the coefficient
$\,\ax\,$ in the expansion \f{\su.5} of the function $\,Q\,$ related to
$\,\si\,$ via the first formula in \f{\as.1}. That these two functionals
coincide is now clear since they agree on the basis $\,1,E,F\,$ of the
solution space: namely, the former functional yields $\,m\,$ for $\,\si=1\,$
(which is obvious) and $\,0\,$ for $\,\si=E\,$ or $\,\si=F\,$ (which is
nothing else than \f{\su.4}, with $\,\dot\vx=m\ps(\ps-2)/(\ps-1)^2$,
rewritten in terms of the variable $\,\vp=\y\hs\ps$).
\endremark
\remark{Remark \a\su.4}Having run out of letters, we had to choose between
using unorthodox notations and allowing different meanings of identical or
similar symbols. We decided to do the latter; as a result, $\,\bx\,$ in
\f{\su.5} is not the same as in \f{\pr.8}, $\,\cx\,$ in \f{\su.5} is to be
distinguished from the lower-case $\,\y\hs$, and the italic $\,r,\sa\,$ for
the norm function and various curve parameters must not be confused with the
roman $\,\,\ri\hs,\hs\sc\,\,$ for the Ricci tensor and scalar curvature.
\endremark

\head\S\ty. Three types of conformally-Einstein K\"ahler manifolds\endhead
\proclaim{Proposition \a\ty.1}Let\/ $\,\mgmt\,$ satisfy\/ \f{\id.1} with\/
$\,m\ge3$, or\/ \f{\id.2} with\/ $\,m=2$, and let\/ $\,Q:M\to\bbR\,$ be given
by\/ $\,Q=g(\navp,\navp)$. Then\/ $\,Q\,$ is a rational function of\/ $\,\vp$.
More precisely, the open set\/ $\,M'\subset M\,$ on which\/ $\,d\vp\ne0\,$ is
connected and dense in\/ $\,M\,$ and, for\/ $\,\si,\hs\y\,$ as in\/ {\rm Lemma
\a\cm.5}, one of the following three cases occurs\/{\rm:}
\widestnumber\item{(c)}\roster
\item"(a)"$\si=0\,$ identically on\/ $\,M'$.
\item"(b)"$\si\ne0\,$ everywhere in\/ $\,M'$ and\/ $\,\y=0$.
\item"(c)"$\si\ne0\,$ everywhere in\/ $\,M'$ and\/ $\,\y\ne0$.
\endroster
In case\/ {\rm(a)}, {\rm(b)}, or\/ {\rm(c)}, the functions\/
$\,\vp,Q:M\to\bbR\,$ satisfy\/ \f{\sz.1}, or\/ \f{\as.4.ii} or, respectively,
\f{\su.5}, for some constants\/ $\,K,\ah,\ts$, or\/ $\,\ax,\bx,\cx$. In
case\/ {\rm(c)}, $\,\vp\ne\y\,$ everywhere in\/ $\,M\,$ unless\/ $\,\cx=0\,$
in\/ \f{\su.5}.
\endproclaim
\demo{Proof}Corollary \a\mc.3 yields \f{\sr.1}, so that $\,M'$ is connected
and dense (Corollary \a\cm.3) and, by Lemma \a\cm.5, we have (a), (b), or
(c). However, Corollary \a\qs.2 gives \f{\de.1} -- \f{\de.3}, everywhere in
$\,M'$, for the septuple $\,(Q,Y,\hs\sc\hs,\si,\ta,\la,\my)\,$ associated with
$\,\vp$. Therefore, in case (a) (or (b), or (c)), any given point of $\,M'$
has a neighborhood in which \f{\sz.1} (or \f{\as.4.ii}, or \f{\su.5}) holds
for some constants $\,K,\ah,\ts\,$ or $\,\ax,\bx,\cx$. (See \S\sz\ and
Remarks \a\as.5, \a\su.3.) Consequently, this triple of constants forms a
locally constant function $\,M'\to\rtr$ which (as $\,M'$ is connected) must be
constant. Hence \f{\sz.1}, or \f{\as.4.ii}, or \f{\su.5} holds {\it
everywhere\/} in $\,M'$, with a single triple of constants. Denseness of
$\,M'$ in $\,M\,$ implies the same everywhere in $\,M$, while \f{\su.5} yields
$\,\vp\ne\y\,$ on $\,M\,$ unless $\,\cx=0$, since $\,Q\,$ given by \f{\su.5}
with $\,\cx\ne0\,$ has a pole at $\,\ps=1$. This completes the
proof.{\hfill\qd}
\enddemo

\head\S\ec. Examples of conformally-Einstein K\"ahler metrics\endhead
Suppose that $\,m\ge2\,$ is an integer, $\,\iy\subset\bbR\,$ is an open
interval of the variable $\,\vp$, while $\,Q\,$ is a rational function of
$\,\vp$, analytic and positive everywhere in $\,\iy$, and
\widestnumber\item{(iii)}\roster
\item"(i)"$\,Q\,$ is defined by \f{\sz.1} for some constants
$\,K,\ah,\ts$, or
\item"(ii)"$\,Q\,$ has the form \f{\as.4.ii} with some constants
$\,\nk,\ah,\ts$, or, finally,
\item"(iii)"$\,Q\,$ is given by \f{\su.5} for some constants
$\,\ax,\bx,\cx,\y\,$ with $\,\y\ne0$.
\endroster
We further set $\,\y=0\,$ in case (ii); in case (i), we leave $\,\y\,$
undefined. Replacing $\,\iy\,$ with a subinterval, we also assume that
$\,\y\notin\iy\,$ in cases (ii), (iii), and define $\,\ve\in\{-\hs1,0,1\}\,$
by $\,\ve=0\,$ (case (i)) or $\,\ve=\,\text{\rm sgn}\,(\vp-\y)\,$ for all
$\,\vp\in\iy\,$ (cases (ii), (iii)).

Using our $\,Q\,$ and a fixed real constant $\,a\ne0$, let us now select a
positive $\,C^\infty$ function $\,r\hs$ on $\,\iy\,$ with $\,Q\,dr/d\vp=ar$.
Thus, $\,\vp\mapsto r\,$ is a diffeomorphism of $\,\iy\,$ onto an open
interval $\,\jy\subset(0,\infty)$. Replacing $\,\vp\,$ with $\,r$, we will
treat functions of $\,\vp\in\iy$, including $\,Q$, as functions of
$\,r\in\jy$. Further $\,C^\infty$ functions of $\,r\in\jy$, introduced in this
way, are $\,\theta=Q/(ar)^2$ and $\,\fe$, with $\,\fe=1\,$ in case (i) and
$\,\fe=\,2\hs|\vp-\y|=2\ve\hs(\vp-\y)\,$ in cases (ii), (iii).

Let $\,(N,h)\,$ be a K\"ahler-Einstein manifold of complex
dimension $\,m-1\,$ with the Ricci tensor $\,\,\rih=\kx\hs h$, where the
constant $\,\kx\,$ is defined by
$$\text{\rm$\kx=(3-2m)K\,$ in case (i), $\,\kx=\ve\nk\,$ in case (ii),
$\,\kx=\ve m\ax/\y\,$ in case (iii).}\ff\ec.1$$
Also, let there be given a holomorphic line bundle
$\,\Cal L\,$ over $\,N\,$ with a Hermitian fibre metric $\,\langle\,,\rangle$,
and a $\,C^\infty$ connection in $\,\Cal L\,$ with a $\,J$-invariant
horizontal distribution
$\,\Cal H\,$ (cf. Remark \a\lc.4), making $\,\langle\,,\rangle\,$ parallel and
having the curvature form (Remark \a\cc.1) equal to $\,-\hs2\ve a\,$ times the
K\"ahler form of $\,(N,h)\,$ (see \f{\kr.1}). Using the symbol $\,r\hs$ also
for the norm function of $\,\langle\,,\rangle\,$ (Remark \a\cc.2), we will
treat $\,\fe,\theta,\vp,Q\,$ as functions $\,M'\to\bbR\hs$, where
$\,M'\subset\Cal L\smallsetminus N\,$ is a fixed connected open subset of the
pre\-im\-age $\,r^{-1}(\jy)\,$ (notation of \f{\cc.4}). We now define a metric
$\,g\,$ on $\,M'$ by declaring $\,\Cal H\,$ to be
$\,g$-or\-thog\-o\-nal to the vertical distribution $\,\Cal V\,$ (\S\cc) and
setting $\,g=\fe\hs\proj^*\nh h\,$ on $\,\Cal H\,$ and
$\,g=\theta\,\text{\rm Re}\hskip1pt\langle\,,\rangle\,$ on $\,\Cal V$. Here
$\,\proj:\Cal L\to N\,$ is the bundle projection, and
$\,\,\text{\rm Re}\hskip1pt\langle\,,\rangle\,$ is the standard Euclidean
metric on each fibre of $\,\Cal L$.
\remark{Remark \a\ec.1}The above construction is, obviously, a special case of
that in \S\xm, with a description modified as in Remark \a\dx.1. Thus (see
\S\dx), $\,(M'\!,g)\,$ is a K\"ahler manifold and $\,\vp\,$ is a special \krp\
on $\,(M'\!,g)$. There are {\it just three additional assumptions\/} on the
data $\,\iy,\vp,Q,a,\ve,\dots\,$ that distinguish the present case from the
general situation in \S\xm: first, $\,m\ge2\,$ (while, in \S\xm, $\,m=1\,$ is
also allowed); next, $\,Q\,$ is required to satisfy (i), (ii) or (iii) above
(rather than being just any positive $\,C^\infty$ function of the variable
$\,\vp\in\iy\hs$); and, finally, the function $\,\kx:N\to\bbR\,$ with
$\,\,\rih=\kx\hs h\,$ is now assumed to be a specific constant with \f{\ec.1}:
in \S\xm\ it need not be constant when $\,m=2$, and can be {\it any\/}
constant if $\,m\ge3$.
\endremark
\remark{Remark \a\ec.2}If we relax the third ``additional assumption'' in
Remark \a\ec.1 and just require $\,\kx\,$ with $\,\,\rih=\kx\hs h\,$ to be a
function on $\,N$, constant unless $\,m=2$, our $\,\vp\,$ will still be a
special \krp\ on the K\"ahler manifold $\,(M'\!,g)\,$ (cf. \S\dx).
Condition \f{\ec.1} then holds if and only if the septuple
$\,(Q,Y,\hs\sc\hs,\si,\ta,\la,\my)\,$ associated with $\,\vp\,$ (Definition
\a\sr.2) satisfies \f{\de.3.iii}.

In fact, by (c) in \S\dx, $\,\si=0\,$ identically in
case (i) (as $\,\ve=0$), while, in cases (ii), (iii), $\,\ve=\pm\hs1$, and so
$\,2\si=Q/(\vp-\y)$, with the constant $\,\y\,$ chosen above. Since, in all
cases, $\,2\ta=\hs dQ/d\vp$, $\,Y=2\ta+2(m-1)\hs\si\,$ and
$\,2\my=-\hs dY/d\vp\,$ (Lemma \a\sr.5), $\,\ta,Y,\my\,$ are given either by
the formulae following \f{\sz.1} (in case (i)), or by those listed in
\f{\as.1}, with $\,\si'=\,d\si/d\vp\,$ (cases (ii), (iii)). Also, (b) in
\S\dx\ provides a formula for $\,\la$, with $\,f\,$ as above. One now easily
verifies that the difference $\,2(m-1)(\ta-\si)-(\la-\my)\vp\,$ (cf.
\f{\de.3.iii}) equals $\,[(3-2m)K-\kx]\hs\vp\,$ (in case (i), proving our
claim), or $\,\vp/(\vp-\y)\,$ times the difference between \f{\as.3} and
$\,\ve\kx/2\,$ (in cases (ii), (iii)). Our claim now also follows in case
(ii), since \f{\as.3} with $\,\y=0\,$ equals $\,K/2\,$ (by \f{\as.4.i}).
Finally, in case (iii), \f{\as.3} equals $\,m\ax/(2\y)\,$ (see Remark
\a\su.3), so that \f{\de.3.iii} is equivalent to $\,\kx=\ve m\ax/\y$, as
required.
\endremark
\proclaim{Proposition \a\ec.3}Let\/ $\,(M'\!,g)\,$ and\/ $\,\vp,Q:M'\to\bbR\,$
be obtained as above, for some integer\/ $\,m\ge2$. Then\/ $\,(M'\!,g)\,$ is a
K\"ahler manifold and the conformally related metric\/ $\,\tilde g=g/\vpsq$,
defined wherever\/ $\,\vp\ne0$, is Einstein, while\/ $\,Q=g(\navp,\navp)$ and
the Laplacian of\/ $\,\vp\,$ is a function of\/ $\,\vp$.
\endproclaim
In fact, according to Remark \a\ec.1, $\,(M'\!,g)\,$ is a K\"ahler manifold
and $\,\vp\,$ is a special \krp\ on $\,(M'\!,g)\,$ (as in \f{\sr.1}), with
$\,g(\navp,\navp)=Q\,$ (see Remark \a\dx.1). Since we assumed \f{\ec.1},
Remark \a\ec.2 gives \f{\de.3.iii}, and our assertion follows from Remark
\a\mc.4, the claim about the Laplacian $\,Y=\Delta\vp\,$ being obvious since
$\,\ta\,$ and $\,\si\,$ are functions of $\,\vp\,$ and
$\,Y=2\ta+2(m-1)\hs\si\,$ (see Remark \a\ec.2).{\hfill\qd}

\head\S\ls. Local structure of conformally-Einstein K\"ahler metrics\endhead
We can now prove our main local result: locally, at points in general
position, and up to biholomorphic isometries, the only conformally-Einstein
K\"ahler metrics in complex dimensions $\,m>2\,$ are those described in
Proposition \a\ec.3.
\proclaim{Theorem \a\ls.1}Let\/ $\,\vp:M\to\bbR\,$ be a nonconstant\/
$\,C^\infty$ function on a K\"ahler manifold\/ $\,(M,g)\,$ of
complex dimension\/ $\,m\ge2\,$ such that the conformally related metric\/
$\,\tilde g=g/\vpsq$, defined wherever\/ $\,\vp\ne0$, is Einstein. If\/
$\,m=2$, let us also assume that\/ $\,d\vp\wedge\hs d\Delta\vp=0\,$ everywhere
in\/ $\,M$. Then the open set\/ $\,M'\subset M\,$ given by\/ $\,d\vp\ne0\,$ is
dense in\/ $\,M$, and every point of\/ $\,M'$ has a neighborhood on which\/
$\,g\,$ and\/ $\,\vp\,$ are, up to a biholomorphic isometry, obtained as in\/
{\rm\S\ec} with a rational function\/ $\,Q\,$ of the form\/ \f{\sz.1}, or\/
\f{\as.4.ii}, or\/ \f{\su.5}.

More precisely, in cases\/ {\rm(a)}, {\rm(b)}, {\rm(c)} of\/ {\rm Proposition
\a\ty.1}, the function\/ $\,Q\,$ used in the construction of\/ $\,g\,$ and\/
$\,\vp\,$ satisfies, respectively, {\rm(i)}, {\rm(ii)}, or\/ {\rm(iii)} in\/
{\rm\S\ec}. In all cases\/ $\,Q\,$ equals\/ $\,g(\navp,\navp)\,$ if one uses
its dependence on\/ $\,\vp\,$ to treat it as a function\/ $\,M\to\bbR\hs$.
\endproclaim
In fact, Corollary \a\mc.3 yields \f{\sr.1}, so that, by Theorem \a\ct.1,
$\,M'$ is dense in $\,M\,$ and, locally in $\,M'$, our $\,g,\vp\,$ are
obtained as in \S\xm. Moreover, the objects \f{\xm.1} used in \S\xm\ may be
replaced by the data $\,\iy,\vp,Q,a,\ve,\dots$, with
$\,Q=g(\navp,\navp):M\to\bbR\,$ treated as a function of the variable
$\,\vp\in\iy\,$ (see Remark \a\dx.1). All we now need to show is that the
latter data satisfy the ``three additional assumptions'' in Remark \a\ec.1. To
this end, note that $\,m\ge2$, while $\,Q\,$ has the form (i), (ii) or (iii)
in \S\ec\ (by Proposition \a\ty.1). Finally, $\,\kx:N\to\bbR\,$ with
$\,\,\rih=\kx\hs h\,$ is given by \f{\ec.1} according to Remark \a\ec.2, as
\f{\de.3.iii} holds according to Remark \a\mc.4.{\hfill\qd}
\remark{Remark \a\ls.2}For any given integer $\,m\ge2$, the local
biholomorphic-isometry types of quadruples $\,\mgmt\,$ satisfying \f{\id.1}
(when $\,m\ge3$), or \f{\id.2} (when $\,m=2$), depend on three real constants
along with an arbitrary local biholomorphic-isometry type of a
K\"ahler-Einstein metric $\,h\,$ in complex dimension $\,m-1$. In fact,
Theorem \a\ls.1 states that locally, at points with $\,d\vp\ne0$, some data
analogous to \f{\xm.1} with $\,\jy,r,\theta,\vp,\fe\,$ replaced by
$\,\iy,\vp,Q\,$ lead to $\,g\,$ and $\,\vp\,$ via an explicit construction. Of
these data, $\,a\,$ is not an invariant (as it may be arbitrarily chosen for
any given $\,g,\vp$, cf. Remark \a\lc.3), while
$\,\ve\in\{-\hs1,0,1\}\,$ is a discrete parameter and $\,m\,$ is fixed, which
leaves, aside from $\,h$, just the parameters $\,\ax,\bx,\cx,\y\,$ or
$\,\nk,\ah,\ts$, representing $\,Q\,$ in (i) -- (iii) of \S\ec). Formula
\f{\ec.1} for the Einstein constant $\,\kx\,$ of $\,h\,$ now reduces the
number of independent parameters by one, and, according to Remark \a\ec.1,
this is the only constraint (cf. Remark \a\ls.3 below).

The parameters $\,\ax,\bx,\cx,\y\,$ or $\,\nk,\ah,\ts\,$ explicitly appear
in the geometry of $\,(M,g)\,$ and $\,\vp$, namely, through a functional
relation between $\,\vp\,$ and $\,Q=g(\navp,\navp)\,$ (Proposition \a\ty.1).
Similarly, $\,h\,$ may be explicitly constructed as a metric on a local leaf
space for the distribution $\,\Cal V\,$ with \f{\sr.3} (see the proof of
Theorem \a\ct.1).
\endremark
\medskip
The concept of a local-isometry type in Remark \a\ls.2 is the one
traditionally used in local differential geometry, where the restrictions of
a real-analytic metric to two disjoint open subsets of the underlying manifold
are regarded as belonging to the same type. The number of independent real
parameters would, however, increase by one (for $\,m=2$) or by $\,2m-1\,$ (for
$\,m\ge3$) if, instead, our ``types'' involved a base point (at which
$\,d\vp\ne0$) and base-point preserving isometries.
\remark{Remark \a\ls.3}Neither the fibre metric $\,\langle\,,\rangle\,$ nor
the connection in $\,\Cal L$, subject to the conditions listed in \S\ec, give
rise to any additional parameters in Remark \a\ls.2: locally, they form a
single equivalence class. In fact, let us fix a holomorphic local trivializing
section $\,w\,$ of $\,\Cal L$, defined on a contractible open set
$\,N'\subset N$, and use this fixed $\,w\,$ to represent connections in
$\,\Cal L\,$ through complex-val\-ued $\,1$-forms $\,\varGamma$, as in Remark
\a\cc.1. The conditions listed in \S\ec: $\,J$-invariance of the horizontal
distribution $\,\Cal H$, and a specific choice of the curvature form
$\,\varOmega$, now mean that $\,\varGamma\,$ is of type $\,(1,0)\,$ (i.e.,
$\,\varGamma(y):T_yN\to\bbC\,$ is complex-linear for every $\,y\in N'$, cf.
\f{\cc.7}), and $\,d\varGamma=-\hs i\hs\varOmega\,$ for a prescribed
real-val\-ued closed $\,2$-form $\,\varOmega\,$ of type $\,(1,1)$.

That a $\,(1,0)\,$ form $\,\varGamma\,$ with
$\,d\varGamma=-\hs i\hs\varOmega\,$ exists on $\,N'$ is clear as we may
choose a function $\,\varphi:N'\to\bbR\,$ with
$\,i\hs\varOmega=\partial\overline{\partial}\varphi$, and set
$\,\varGamma=\partial\varphi$. Any other such form $\,\tilde\varGamma\,$
obviously equals $\,\varGamma-\hs d\varPhi$, where $\,\varPhi:N'\to\bbC\,$ is
holomorphic. The connections represented by $\,\varGamma\,$ and
$\,\tilde\varGamma\,$ thus are holomorphically equivalent, as the latter is
the image of the former under the bundle automorphism (over $\,N'$), sending
the section $\,w\,$ (chosen above) to $\,\tilde w=e^\varPhi w$. Finally, a
connection with the curvature form $\,\varOmega\,$ admits a parallel fibre
metric $\,\langle\,,\rangle$, unique up to a constant factor, since
$\,\langle w,w\rangle\,$ is a positive $\,C^\infty$ function subject to the
sole requirement that
$\,d\log\hs\langle w,w\rangle=2\,\text{\rm Re}\,\varGamma$, the form
$\,\,\text{\rm Re}\,\varGamma\,$ being closed due to real-val\-ued\-ness of
$\,i\hskip1ptd\varGamma=\varOmega$. A holomorphic bundle automorphism as
above, multiplied by a suitable real scale factor, will not only identify the
two connections, but also send one fibre metric onto the other.
\endremark

\head\S\tg. Appendix: Compact product manifolds with \ \f{\id.1}\endhead
The examples described next are well-known, at least for $\,m=2\,$ (cf.
\cite{\pet, pages 343 and 345 -- 350} and \cite{\hem}), and have an immediate
generalization (see \cite{\dmr}) to the case where $\,M=N\times S^2\,$ is
replaced by a suitable flat $\,S^2$ bundle over $\,N$.

Given an integer $\,m\ge2\,$ and a real number
$\,K>0$, let $\,(M,g)\,$ be a product K\"ahler manifold whose factors are an
oriented $\,2$-sphere $\,S^2$ with a metric $\,\gx\,$ of constant Gaussian
curvature $\,K$, and any compact K\"ahler-Einstein manifold $\,(N,h)\,$ of
complex dimension $\,m-1\,$ with the negative-definite Ricci tensor
$\,\,\rih=\kx\hs h$, where $\,\kx=(3-2m)K$. Treating our $\,S^2$ as the sphere
of radius $\,1/\sqrt K\,$ about $\,0\,$ in a Euclidean $\,3$-space $\,V$, let
us choose $\,\vp:S^2\to\bbR\,$ to be the restriction to $\,S^2$ of any nonzero
linear homogeneous function $\,V\to\bbR\hs$. In other words, $\,\vp\,$ is an
eigenfunction of the spherical Laplacian for the first nonzero eigenvalue
$\,-\hs2K\,$ (see below). The quadruple $\,\mgmt$, with $\,\vp\,$ viewed as a
function on $\,M\,$ constant in the direction of the $\,N\,$ factor, then
satisfies \f{\id.2}. Specifically, the tensor field $\,\bz\,$ on $\,M\,$ given
by $\,\bz=2(m-1)\hs\nabla d\vp+\vp\hskip1pt\ri\,\,$ equals
$\,(3-2m)K\nh\vp\hs g$, which implies \f{\ck.2} with $\,n=2m\,$ (and hence
\f{\id.1}), while $\,\Delta\vp=-\hs2K\nh\vp\,$ (which gives \f{\id.2}).

In fact, let $\,\dot\vp=\hs d\hs[\vp(x(\sa))]/d\sa\,$ for a fixed $\,C^2$
curve $\,\sa\mapsto x(\sa)\,$ in a Riemannian manifold $\,(M,g)\,$ and a
$\,C^2$ function $\,\vp:M\to\bbR\hs$. By \f{\pr.4}, $\,\dot\vp=g(v,\dot x)$,
where $\,\dot x=\hs dx/d\sa\,$ and $\,v=\navp$. Applying $\,d/d\sa\,$ again,
we obtain $\,\ddot\vp=g(\nabla_{\dot x}v,\dot x)+g(v,\nabla_{\dot x}\dot x)$.
If the curve is a geodesic, i.e., $\,\nabla_{\!\dot x}\dot x=0$, this becomes
$\,\ddot\vp=(\nabla d\vp)(\dot x,\dot x)$, by \f{\pr.5} with $\,u=w=\dot x$.
The last relation, applied to our linear function $\,\vp:S^2\to\bbR\,$ and all
great circles in $\,S^2$ with standard sine/cosine parameterizations
$\,\sa\mapsto x(\sa)$, yields $\,\nabla d\vp=-\hs K\nh\vp\gx\,$ in
$\,(S^2\!,\gx)$, and hence also in $\,(M,g)\,$ (where $\,\gx\,$ is identified
with its pullback to $\,M$). Thus, $\,\Delta\vp=-\hs2K\nh\vp\,$ in both
$\,(S^2\!,\gx)\,$ and $\,(M,g)$.

Finally, to see that $\,\bz=(3-2m)K\nh\vp g$, note that the factor
distributions in $\,M\,$ are $\,\,\ri$-or\-thog\-o\-nal and
$\,\nabla d\vp$-or\-thog\-o\-nal, and hence $\,\bz$-or\-thog\-o\-nal to each
other, while, as $\,\nabla d\vp=-\hs K\nh\vp\gx$, the restrictions of
$\,\nabla d\vp$, $\,\,\ri\,\,$ and $\,\bz\,$ to the $\,N\,$ (or, $\,S^2$)
factor distribution are $\,0$, $\,\,\rih=(3-2m)Kh\,$ and
$\,(3-2m)K\nh\vp h\,$ (or, respectively, $\,-\hs K\nh\vp\gx$,
$\,K\gx\,$ and $\,(3-2m)K\nh\vp\gx$), which proves our claim.

According to Corollaries \a\mc.3 and \a\rr.2(iii), the
quadruples $\,\mgmt$, just described all satisfy condition (a) in Proposition
\a\ty.1, and so, by Theorem \a\ls.1, they arise (locally, at points with
$\,d\vp\ne0$) from the construction of \S\ec, with $\,Q\,$ as in (i) of \S\ec.
More precisely, $\,Q=g(\navp,\navp)\,$ is given by \f{\sz.1} with $\,K\,$ as
above, $\,\ah=0$, and some $\,\ts\in(-\infty,0)\,$ depending on the
choice of the linear function $\,V\to\bbR$. In fact, we also have
$\,Q=\gx(\navp,\navp)\,$ (the $\,\gx$-gradient in $\,(S^2\!,\gx)$) and, since
$\,\nabla d\vp=-\hs K\nh\vp\gx\,$ in $\,(S^2\!,\gx)$, \f{\pr.6.ii} and
\f{\pr.4} applied to $\,(S^2\!,\gx)\,$ (rather than $\,(M,g)$) give
$\,dQ=-\hs2K\nh\vp\hs d\vp$, so that $\,Q+K\nh\vpsq$ is a positive constant,
as required.

\head\S\bb. Appendix: B\'erard Bergery's and Page's examples\endhead
Section 8 of Lionel B\'erard Bergery's 1982 paper \cite{\ber} describes a
family of non-K\"ahler, Einstein metrics on holomorphic $\,\bbCP^1$ bundle
spaces $\,M\,$ in all complex dimensions $\,m\ge2$, which includes the Page
metric (see \cite{\pag}) with $\,m=2$. Every metric $\,\tilde g\,$ in that
family is globally conformal to a K\"ahler metric $\,g\,$ (see \cite{\ber}),
i.e., $\,\tilde g=g/\vpsq$ for some nonconstant $\,C^\infty$ function
$\,\vp:M\to\bbR\smallsetminus\{0\}$. Our Theorem \a\ls.1 now implies that
locally, at points with $\,d\vp\ne0$, those $\,g,\vp\,$ may be obtained via
the construction in \S\ec\ (also for $\,m=2$, since the additional assumption
in \f{\id.2} happens to hold as well).

Instead of using Theorem \a\ls.1, we will now verify the last statement
directly, by explicitly showing that B\'erard Bergery's construction is a
special case of ours.

By saying `special case' we do not claim that the existence of the {\it
compact\/} manifolds found by B\'erard Bergery easily follows from our
results. Our local approach ignores the boundary conditions
necessary for compactness which, in B\'erard Bergery's construction, amount to
a very careful choice of the pair $\,\lambda,x\,$ in the table below; in the
exposition that follows, $\,\lambda=\ts/m\,$ is arbitrary, and $\,x\,$ does
not appear at all. (B\'erard Bergery's choice of $\,\lambda,x\,$ leads, when
$\,N\,$ is compact, to a compactified version $\,(M,\tilde g)\,$ of the
Einstein manifold $\,(M'\!,\tilde g)\,$ described below.)  For compact
manifolds that generalize B\'erard Bergery's examples, see \cite{\dmg}.

Since the original paper \cite{\ber} is difficult to obtain, our account of
B\'erard Bergery's examples follows their presentation in \cite{\bes, Section
K of Chapter 9, pp. 273 -- 275}. It is described there how an Einstein metric
$\,g\,$ is constructed using a real variable $\,t$, functions
$\,f,h,\varphi\,$ of $\,t$, and $\,P\,$ of $\,\varphi$, an integer $\,n$, a
rational constant written as $\,s/q$, real constants $\,l,\lambda,x$, a
K\"ahler-Einstein manifold $\,(B,\check g)$, and a $\,\,\text{\rm U}\hs(1)\,$
bundle over $\,B\,$ with the projection mapping $\,p$. The symbols just listed
are literally quoted from \cite{\bes}; our outline of the discussion in
\cite{\bes} employs different notations, consistent with the rest of the
present paper, and related to those in \cite{\bes} as follows:
\vskip5pt
\centerline{
\vbox{\offinterlineskip
\hrule
\halign
{&\vrule#&\strut
\hskip2pt\hfil#\hfil\hskip2pt\cr
&\hfil\phantom{$^{1^1}$}in \cite{\bes}:\phantom{$_j$}\hfil&
&$\,g\,$&&$\,t\,$&&$\,f\,$&&$\,h\,$&&$\,\varphi\,$&&$\,P\,$&&$\,n\,$&
&$\,s/q\,$&&$\,l\,$&&$\,\lambda\,$&&$\,x\,$&&$\,B\,$&&$\,\check g\,$&
&$\,p\,$&\cr
\noalign{\hrule}
&\hfil\phantom{$^{1^1}$}below:\phantom{$_j$}\hfil&&$\,\tilde g\,$&&$\,\sa\,$&
&$\,\ly\,$&&$\,(1-\varphi^2)^{1/2}\,$&&$\,\varphi\,$&&$\,P\,$&&$\,2m\,$&&$\,q\,$&
&$\,\ell\,$&&$\,\ts/m\,$&&$\,$(not used)$\,$&&$\,N\,$&&$\,h\,$&&$\,\proj\,$&\cr
\noalign{\hrule}
}}    }
\vskip4pt
Specifically, given an integer $\,m\ge2$, let $\,P(\varphi)\,$ be the unique
even polynomial in the real variable $\,\varphi\,$ with $\,P(0)=1\,$ and
$\,\dsq P/d\varphi^2=2m(1-\varphi^2)^{m-1}$. It follows now that, for
$\,E,F\,$ as in \f{\su.1} and any $\,\ps\in\bbR\smallsetminus\{0\hs,1\}$,
$$2(1-2m)E(0)\hs P(\varphi)\,=\,
(1-\varphi^2){}^m[\hs2E(\ps)-F(\ps)\hs]\hs,\hskip5pt\text{\rm where}\hskip4pt
\varphi=(2-\ps)/\ps\hs.\ff\bb.1$$
Namely, $\,(1-\varphi^2){}^mF(\ps)=-\hs4^m\varphi$, and so \f{\bb.1} amounts
to $\,2(1-\varphi^2){}^mE(\ps)+4^m\varphi=2(1-2m)E(0)\hs P(\varphi)$, which
one easily verifies by induction on $\,m\ge1$. In fact, for any $\,m\ge2$,
marking the dependence of $\,E,P\,$ on $\,m\,$ with a subscript, we have
$\,E_m(\ps)=4(1-\varphi^2)^{-1}E_{m-1}(\ps)+2(2-3/m)\hs E_{m-1}(0)\,$ and
$\,mE_m(0)/2=(2m-3)\hs E_{m-1}(0)\,$ (which are simply the recursion formulae
preceding \f{\su.3}), while
$\,(1-2m)P_m(\varphi)=(1-\varphi^2){}^m-2mP_{m-1}(\varphi)$, which is
immediate from the obvious expansion
$\,P_m(\varphi)=1+2m\hs\sum_{k=1}^m(-1)^{k-1}{m-1\choose k-1}(2k)^{-1}
(2k-1)^{-1}\varphi^{2k}$.

Next, let us choose a complex line bundle $\,\Cal L\,$ over a
K\"ahler-Einstein manifold $\,(N,h)\,$ of complex dimension $\,m-1$, with the
Ricci tensor $\,\,\rih=2mh$, and a $\,C^\infty$ connection in $\,\Cal L$,
with a horizontal distribution denoted $\,\Cal H$, which makes some Hermitian
fibre metric $\,\langle\,,\rangle\,$ parallel and has the curvature form
(Remark \a\cc.1) equal to $\,2mq\,$ times the K\"ahler form of $\,(N,h)\,$
(see \f{\kr.1}), where $\,q\,$ is a fixed rational number with $\,0<q<1$.

It follows that $\,\Cal L\,$ admits a structure of a holomorphic line
bundle for which $\,\Cal H\,$ is $\,J$-invariant, $\,J\,$ being the complex
structure tensor on the total space. In fact, one defines an integrable almost
complex structure $\,J\,$ on $\,\Cal L\,$ by requiring that $\,\Cal H\,$ be
$\,J$-invariant, and the projection $\,\proj:\Cal L\to N\,$ as well as all the
fibre inclusions $\,\Cal L_y\to\Cal L$, $\,y\in N$, be
pseu\-do\-hol\-o\-mor\-phic in the sense that their differentials at all
points are complex-linear; cf. \cite{\bes, Remark 9.126(b), p. 275}.

Finally, let $\,\varphi\,$ be a function $\,(0,\ell)\to(-\hs1,1)\,$ of the
real variable $\,\sa\in(0,\ell)$, with some $\,\ell>0$, for which, writing
$\,\varphi\hs'=\hs d\varphi/d\sa\,$ we have, for a suitable constant $\,\ts$,
$$(1-\varphi^2)^{m-1}[\varphi\hs'\hs]^2\,
=\,(1-\varphi^2)^m\,+\,(2m-1-\ts/m)\,P(\varphi)\,>\,0\,.\ff\bb.2$$
We also fix a positive function $\,r\hs$ of the variable $\,\sa\,$ such that
$\,\log\hs r\,$ is an antiderivative of $\,1/\ly\,$ or $\,-\hs1/\ly$, i.e.,
$\,dr/d\sa=\pm\hs r/\ly(\sa)$, where $\,\ly=\varphi\hs'\nh/(mq)$. Thus,
$\,\sa\mapsto r\,$ is a $\,C^\infty$ diffeomorphism of $\,(0,\ell)\,$ onto a
subinterval of $\,(0,\infty)$, which allows us to treat $\,\varphi\,$ as a
function of $\,r$. The symbol $\,r\hs$ also stands for the
norm function of $\,\langle\,,\rangle$, cf. Remark \a\cc.2, so that
$\,\varphi\,$ now becomes a function $\,M'\to\bbR\,$ on a suitable open
connected subset $\,M'\subset\Cal L\smallsetminus N\,$ of the total space of
$\,\Cal L\,$ (cf. \f{\cc.4}).

Any of B\'erard Bergery's Einstein metrics $\,\tilde g\,$ is defined, on this
$\,M'$, so that $\,\tilde g=(1-\varphi^2)\hs\proj^*\nh h\,$ on $\,\Cal H\,$
while $\,\tilde g
=[\varphi\hs'\nh/(mqr)]^2\,\text{\rm Re}\hskip1pt\langle\,,\rangle\,$ on the
vertical distribution $\,\Cal V$, and $\,\tilde g(\Cal H,\Cal V)=\{0\}$.
According to Remark \a\bb.1 below, our description of $\,\tilde g\,$ on
$\,\Cal V\,$ is equivalent to that in \cite{\bes}, where the restriction
$\,\tilde\gx\,$ of $\,\tilde g\,$ to any fibre of $\,\Cal L\,$ is defined as
in Remark \a\bb.1, with $\,\ly\,$ depending on $\,\sa\,$ as above, that is,
$\,\ly=\varphi\hs'\nh/(mq)$.

Using a metric $\,\tilde g\,$ on $\,M'$ thus obtained, let us now define real
constants $\,a,\ve,\y,\kx$, functions $\,\vp,Q:M'\to\bbR$, and a metric
$\,g\,$ on $\,M'$, by setting $\,a=-\hs mq$, $\,\ve=1$, $\,\y=1/2$,
$\,\kx=2m$, as well as $\,Q=\vpsq[\varphi\hs'\hs]^2$ and $\,g=\vpsq\tilde g$,
where $\,\vp=\y\hs\ps\,$ and $\,\ps:M'\to\bbR\,$ is related to $\,\varphi\,$
by $\,\varphi=(2-\ps)/\ps$. Thus, $\,1-\varphi^2=2\ve\hs(\vp-\y)/\vpsq$,
and so \f{\bb.1}, \f{\bb.2} give \f{\su.5} for $\,\ax=1\,$ and some specific
$\,\bx,\cx$, namely
$${Q\over\ps-1}\,=\,1\,+\,{(2m-1)^{-1}m^{-1}\ts\,-\,1\over2\hs E(0)}\,
[\hs2E(\ps)-F(\ps)\hs]\hs.\ff\bb.3$$
B\'erard Bergery's Einstein metric now has the form $\,\tilde g=g/\vpsq$, with
$\,(M'\!,g)\,$ and
$\,\vp,Q:M'\to\bbR\,$ constructed as in \S\ec, cf. Proposition \a\ec.3; the
function $\,Q\,$ of the variable $\,\vp\,$ used here is of type (iii) in
\S\ec, with $\,a,\ve,\y,\kx\,$ defined above.
\remark{Remark \a\bb.1}In \cite{\bes} the restriction of the metric
$\,\tilde g\,$ to any fibre of $\,\Cal L\,$ appears as a metric on a cylinder,
while in this presentation we explicitly define it on an open annulus in the
fibre. Equivalence of the two approaches can be seen as follows.

Let $\,S^1$ be the unit circle centered at $\,0\,$ in
a one\diml\ complex vector space $\,V\,$ with a Hermitian inner product
$\,\langle\,,\rangle$. Any positive $\,C^\infty$ function $\,\ly\,$ of the
real variable $\,\sa$, defined on an interval of the form $\,(0,\ell)$, gives
rise to a metric $\,\tilde\gx\,$ on the cylinder $\,(0,\ell)\times S^1$ such
that for each $\,\sa\in(0,\ell)\,$ (or, $\,\zx\in S^1$),
$\,\sa\mapsto(\sa,\zx)\,$ (or, $\,\zx\mapsto(\sa,\zx)$) is a
$\,\tilde\gx$-isometric embedding (or, a closed curve of $\,\tilde\gx$-length
$\,2\pi\hs\ly(\sa)$, with $\,\tilde\gx$-constant speed).

The Riemannian surface $\,((0,\ell)\times S^1,\tilde\gx)\,$ then admits a
canonical conformal diffeomorphism onto an annulus
$\,V'\subset V\smallsetminus\{0\}$, centered at $\,0$. In fact, any fixed
norm-preserving isomorphic identification $\,V=\bbC\,$ gives rise to the polar
coordinates $\,r,\theta\,$ in $\,V\,$ and the corresponding Cartesian
coordinates $\,x=r\cos\theta$, $\,y=r\sin\theta$. (This $\,\theta\,$ is not
related to $\,\theta\,$ in the preceding sections.) The standard
Euclidean metric $\,\,\text{\rm Re}\hskip1pt\langle\,,\rangle\,$
on $\,V\,$ then equals $\,dx^2+\hs dy^2=\hs dr^2+r^2d\theta^2$, while
$\,\tilde\gx=\hs d\sa^2+\ly^2\hs d\theta^2$, where $\,\theta$, restricted to
$\,S^1$, serves as a local coordinate for $\,S^1$. If we now choose a
diffeomorphism $\,\sa\mapsto r\,$ of $\,(0,\ell)\,$ onto a subinterval of
$\,(0,\infty)\,$ such that $\,dr/d\sa=\pm\hs r/\ly(\sa)$, then
$\,\tilde\gx=(\ly/r)^2[\hs dr^2+r^2d\theta^2]$. In other words, the
push-forward of $\,\tilde\gx\,$ under the diffeomorphism
$\,(0,\ell)\times S^1\to V'$ given by $\,(\sa,\zx)\mapsto r\zx\,$ (with
$\,r\hs$ depending on $\,\sa\,$ as above) equals $\,\ly^2\!/r^2$ times the
Euclidean metric.
\endremark

\head\S\is. Appendix: Further integrals of the system \ \f{\de.1} --
\f{\de.4}\endhead
The system \ \f{\de.1} -- \f{\de.4} has further interesting integrals:
$\,\kx=(\la Q\,+\,\si Y)/|\si|\,$ and
$\,\ts=[\my+(m-1)\la]Z^2+2(2m-1)[\ta+(m-1)\hs\si]Z-m(2m-1)Q\,$ (with $\,Z\,$
as in Lemma \a\de.1), in addition to those listed in Lemma \a\de.1. In fact,
wherever $\,\si\ne0$, \f{\de.1} gives $\,\,d\hs[Y\,+\,\la Q/\si]/d\vp=0\,$
(i.e., $\,d\kx/d\vp=0$). Similarly, for $\,m\ge2$, \f{\de.1} -- \f{\de.3.i}
yield $\,d\hskip.2pt\ts/d\vp=0\,$ wherever $\,\la\ne\my$, as one sees
evaluating $\,d\hskip.2pt\ts/d\vp\,$ from \f{$*$} in \S\de\ and \f{\de.1}, and
then replacing $\,Z\,$ with $\,2(m-1)(\ta-\si)/(\la-\my)$.

Geometrically, constancy of $\,\kx\,$ and $\,\ts\,$ is related to Schur's
theorem for the Ricci curvature, since they are proportional to the scalar
curvatures of some Einstein metrics: for $\,\kx$, it is the metric $\,h\,$
defined in the proof of Theorem \a\ct.1, under the assumption that $\,\vp\,$
satisfies \f{\sr.1} on a K\"ahler manifold of complex dimension $\,m\ge3\,$
and $\,\si\ne0\,$ everywhere in $\,M'$ (see Lemma \a\cm.5, with
$\,Q,Y,\hs\sc\hs,\si,\ta,\la,\my\,$ as in Lemma \a\qs.1); while, in the case
of $\,\ts$, it is the Einstein metric $\,\tilde g=g/\vpsq$ for $\,\mgmt\,$
satisfying \f{\id.1} with $\,m\ge3$, or \f{\id.2} with $\,m=2$, since the
scalar curvature $\,\,\tilde{\sc}\,\,$ of $\,\tilde g\,$ equals $\,2\ts\,$ by
\f{\ck.1} for $\,n=2m\,$ combined with Lemma \a\sr.5(i), where
$\,Y=\Delta\vp$, and the equality $\,Z=\vp\,$ (cf. \f{\de.3.iii} and Corollary
\a\qs.2).

>From now on we assume that $\,Q,Y,\hs\sc\hs,\si,\ta,\la,\my\,$ are $\,C^1$
functions on an interval $\,\iy\,$ of the variable $\,\vp$, satisfying, except
for one subcase in (d), conditions \f{\de.1} -- \f{\de.4} with a fixed integer
$\,m\ge2$. Aside from that subcase, we set $\,\ve=\,\sgn\hs\hs\si\,$ on
$\,\iy\,$ and let $\,\y\,$ be the constant defined in Lemma \a\de.1. Then, for
$\,\kx,\ts\,$ as above,
\widestnumber\item{(e)}\roster
\item"(a)"$\hs\ve\kx/2\,$ equals \f{\as.3}, as one sees writing
$\,\ve\kx=Y+\la Q/\si\,$ and then replacing $\,Y,\la,Q\,$ by the expressions
in \f{\as.1}.
\item"(b)"If $\,\y=0$, \f{\as.4.i} and (a) give $\,\kx=\ve\nk$, with $\,\nk\,$
appearing in \f{\as.4}.
\item"(c)"$\,\ts/m
=-\hs\vpsq(\vp-\y)\hs\si''+[(m-1)\vp-2m\y]\hs\vp\si'+2(2m-1)\y\si$, with
$\,\si'=\hs d\si/d\vp$, as one sees replacing $\,Z\,$ in $\,\ts\,$ by
$\,\vp\,$ (cf. \f{\de.3.iii}), and then using \f{\as.1} to
express $\,\my,\la,\ta,Q\,$ through $\,\si\hs,\si'\nh,\si''$.
\item"(d)"Under the assumptions \f{\de.1} -- \f{\de.4} and $\,\y=0\,$ (or,
\f{\de.1} -- \f{\de.3} and $\,\si=0$), our $\,\ts\,$ coincides with
the constant $\,\ts\,$ in \f{\as.4} (or, \f{\sz.1}), which is immediate from
(c) with $\,\y=0\,$ (or, respectively, the relation $\,Z=\vp$, cf.
\f{\de.3.iii}, and the expressions for $\,\my,\la,\ta,\si,Q\,$ in \S\sz).
\item"(e)"Let $\,\y\ne0\,$ and let $\,\ax,\bx,\cx\,$ be as in \f{\su.5}. Then,
by (a) and the final part of
Remark \a\su.3, $\,\kx=\ve m\ax/\y$. On the other hand,
$\,\ts=(2m-1)m\hs[\ax+\bx E(0)]$, with $\,E\,$ as in \f{\su.1}. In fact, by
(c), $\,\ts/m=\y\hskip1pt\Cal P[\si]$, where $\,\si\,$ is treated as a $\,C^2$
function of the variable $\,\ps=\vp/\y\,$ and, for any such function $\,\si$,
we set $\,\dot\si=\hs d\si/d\ps=\y\hs\si'$ and
$\,\Cal P[\si]=-\hs\ps^2(\ps-1)\hs\ddot\si+[(m-1)\ps-2m]\hs\ps\dot\si
+2(2m-1)\hs\si$. What we need to show is
$\,\y\hskip1pt\Cal P[\si]=(2m-1)\hs[\ax+\bx E(0)]\,$ for every $\,\si\,$ in
the three\diml\ solution space $\,V\,$ mentioned in Remark \a\su.3, with
$\,\ax,\bx\,$ corresponding as in \f{\su.5} to $\,Q=2(\vp-\y)\hs\si$, or,
equivalently, $\,\Cal P[\si]=2(2m-1)\hs\si(0)\,$ for all $\,\si\in V$. This
amounts to equality between two linear functionals on $\,V$, which we can
establish using the basis $\,1,E,F$. The functionals obviously agree on
$\,\si=1$. Also, \f{\su.2} with $\,F(0)=0\,$ (see \f{\su.1}) yields
$\,\Cal R[\si]=0\,$ both for $\,\si=E\,$ and $\,\si=F$, where
$\,\Cal R[\si]=\ps(\ps-2)\hs\dot\si-\vx\hs\si-2(2m-1)\hs\si(0)$. Multiplying
the equality $\,\Cal R[\si]=0\,$ by $\,\ps(\ps-1)(\ps-2)^{-1}$, then applying
$\,d/d\ps$, and noting that $\,\ps(\ps-1)(\ps-2)^{-1}\varXi(\ps)
=m\ps^2-2(m-1)\ps+2+4/(\ps-2)$ and $\,\ps(\ps-1)(\ps-2)^{-1}=\ps+1+2/(\ps-2)$,
we obtain a formula expressing $\,-\hs\ps^2(\ps-1)\hs\ddot\si\,$ through
$\,\si\,$ and $\,\dot\si$. That formula alone now easily gives
$\,\Cal P[\si]-\Cal R[\si]-2(2m-1)\hs\si(0)=[1-2/(\ps-2)^2]\hs\Cal R[\si]$, so
that $\,\Cal P[\si]=2(2m-1)\hs\si(0)\,$ if $\,\si=E\,$ or $\,\si=F$, as
required.
\endroster

\enddocument